\documentclass[12pt]{amsart}

\usepackage{amsfonts,amssymb,stmaryrd,amscd,amsmath,latexsym,amsbsy}

\usepackage{amssymb}
\usepackage{amsfonts}
\usepackage{latexsym}

\newtheorem{theorem}{Theorem}[section]
\newtheorem{lemma}[theorem]{Lemma}
\newtheorem{proposition}[theorem]{Proposition}
\newtheorem{corollary}[theorem]{Corollary}
\theoremstyle{definition}
\newtheorem{definition}[theorem]{Definition}

\newtheorem{example}[theorem]{Example}
\newtheorem{exercise}[theorem]{Exercise}

\newtheorem{remark}[theorem]{Remark}

\newcommand{\e}{{\bold e}}
\newcommand{\Ext}{\text{Ext}}
\newcommand{\Tr}{\text{Tr}}

\newcommand{\diag}{\text{diag}}

\renewcommand{\O}{\mathcal{O}}

\newcommand{\g}{\mathfrak{g}}
\newcommand{\h}{\mathfrak{h}}
\newcommand{\n}{\mathfrak{n}}

\newcommand{\ben}{\begin{enumerate}}
\newcommand{\een}{\end{enumerate}}

\newcommand{\mC}{{\mathcal C}}

\newcommand{\CC}{{\mathbb{C}}}

\theoremstyle{plain}

\newtheorem*{sol}{Solution}

\theoremstyle{definition}

\theoremstyle{remark}

\newcommand{\solu}[1]{\begin{sol}{\bf (\ref{#1})}}

\newcommand{\p}{{\bold p}}

\begin{document}

\title{Lectures on Calogero-Moser systems}

\author{Pavel Etingof}
\address{Department of Mathematics, Massachusetts Institute of Technology,
Cambridge, MA 02139, USA}
\email{etingof@math.mit.edu}

\maketitle

\centerline{\bf To my mother Yelena Etingof on her 75-th birthday, 
with admiration}

\vskip .1in
\centerline{\bf Introduction}
\vskip .1in

Calogero-Moser systems, which were originally discovered 
by specialists in integrable systems, are currently at the crossroads
of many areas of mathematics and within the scope of interests of
many mathematicians. More specifically, these systems and their generalizations 
turned out to have intrinsic connections with such fields 
as algebraic geometry (Hilbert schemes of surfaces),
representation theory (double affine Hecke algebras, Lie groups, quantum
groups), deformation theory (symplectic reflection algebras), 
homological algebra (Koszul algebras), Poisson
geometry, etc. The goal of the present lecture notes is to give an introduction 
to the theory of Calogero-Moser systems, highlighting their
interplay with these fields. Since these lectures are designed
for non-experts, we give short introductions to each of the
subjects involved, and provide a number of exercises. 

We now describe the contents of the lectures in more
detail. 

In Lecture 1, we give an introduction to Poisson geometry 
and to the process of classical Hamiltonian reduction. More
specifically, we define Poisson manifolds (smooth, analytic, and
algebraic), momemt maps and their main properties, and then 
describe the procedure of (classical) Hamiltonian reduction.
We give an example of computation of Hamiltonian reduction 
in algebraic geometry
(the commuting variety). Finally, we define Hamiltonian
reduction along a coadjoint orbit, and give the example which plays a
central role in these lectures -- the Calogero-Moser space 
of Kazhdan, Kostant, and Sternberg. 

In Lecture 2, we give an introduction to classical Hamiltonian 
mechanics and the theory of integrable systems. Then we explain
how integrable systems may sometimes be constructed using
Hamiltonian reduction. After this we define the classical  
Calogero-Moser integrable system using Hamiltonian reduction
along a coadjoint orbit (the Kazhdan-Kostant-Sternberg construction),  
and find its solutions. Then, by introducing coordinates on the 
Calogero-Moser space, we write both the system and the solutions
explicitly, thus recovering the standard results 
about the Calogero-Moser system. Finally, we generalize these
results to construct the trigonometric Calogero-Moser system. 

Lecture 3 is an introduction to deformation theory. 
This lecture is designed, in particular, to enable us 
to discuss quantum-mechanical versions of 
the notions and results of Lectures 1 and 2 in a manner 
parallel to the classical case. Specifically, we 
develop the theory of formal and algebraic deformations of assosiative
algebras, introduce Hochschild cohomology and discuss 
its role in studying deformations, and define universal
deformations. Then we discuss the basics of the theory of 
deformation quantization of Poisson (in particular, symplectic) manifolds, 
and state the Kontsevich quantization theorem. 

Lecture 4 is dedicated to the quantum-mechanical generalization
of the material of Lecture 1. Specifically, we define the notions
of quantum moment map and quantum Hamiltonian reduction. Then we 
give an example of computation of quantum reduction (the
Levasseur-Stafford theorem), which is the quantum analog 
of the example of commuting variety given in Lecture 1.
Finally, we define the notion of quantum reduction with respect
to an ideal in the enveloping algebra, which is the quantum
version of reduction along a coadjoint orbit, and give an example
of this reduction -- the construction of the spherical subalgebra
of the rational Cherednik algebra. Being a quantization of the
Calogero-Moser space, this algebra is to play a central role 
in subsequent lectures. 

Lecture 5 contains the quantum-mechanical version
of the material of Lecture 1. Namely, after recalling the basics 
of quantum Hamiltonian mechanics, we introduce the notion of 
a quantum integrable system. Then we explain how to construct
quantum integrable systems by means of quantum reduction (with
respect to an ideal), and give an example of this which is 
central to our exposition -- the quantum Calogero-Moser system. 

In Lecture 6, we define and study more general classical and quantum
Calogero-Moser systems, which are associated to finite Coxeter
groups (they were introduced by Olshanetsky and Perelomov). 
The systems defined in previous lectures correspond 
to the case of the symmetric group. In general, these integrable systems are
not known (or expected) to have a simple construction using reduction; in their
construction and study, Dunkl operators are an indispensible
tool. We introduce the Dunkl operators (both classical and
quantum), and explain how the Olshanetsky-Perelomov Hamiltonians
are constructed from them. 

Lecture 7 is dedicated to the study of the rational Cherednik
algebra, which naturally arises from Dunkl operators (namely, it is
generated by Dunkl operators, coordinates, and
reflections). Using the Dunkl operator representation, we prove
the Poincare-Birkhoff-Witt theorem for this algebra, and 
study its spherical subalgebra and center. 

In Lecture 8, we consider symplectic reflection algebras, 
associated to a finite group $G$ of automorphisms of a symplectic
vector space $V$. These algebras are natural generalizations of
rational Cherednik algebras (although in general they are not 
related to any integrable system). It turns out that 
the PBW theorem does generalize to these algebras, but its proof does not, since 
Dunkl operators don't have a counterpart. Instead, the proof is
based on the theory of deformations of Koszul algebras, 
due to Drinfeld, Braverman-Gaitsgory, Polishchuk-Positselski, 
and Beilinson-Ginzburg-Soergel. We also study the spherical 
subalgebra of the symplectic reflection algebra,
and show by deformation-theoretic arguments 
that it is commutative if the Planck constant is equal
to zero. 

In Lecture 9, we describe the deformation-theoretic
interpretation of symplectic reflection algebras. 
Namely, we show that they are universal deformations of
semidirect products of $G$ with the Weyl algebra of $V$. 

In Lecture 10, we study the center of the symplectic reflection
algebra, in the case when the Planck constant equals zero. 
Namely, we consider the spectrum of the center, which is an
algebraic variety analogous to the Calogero-Moser space, and show
that the smooth locus of this variety is exactly 
the set of points where the symplectic
reflection algebra is an Azumaya algebra; this requires 
some tools from homological algebra, such as the Cohen-Macaulay
property and homological dimension, which we briefly introduce. 
We also study finite dimensional representations of symplectic
reflection algebras with the zero value of the Planck constant. 
In particular, we show that for $G$ being the symmetric group $S_n$
(i.e. in the case of rational Cherednik algebras of type $A$),
every irreducible representation has dimension $n!$, and 
irreducible representations are parametrized by the
Calogero-Moser space defined in Lecture 1. A similar theorem
is valid if $G=S_n\ltimes \Gamma^n$, where $\Gamma$ is a finite
subgroup of $SL_2(\Bbb C)$. 

Lecture 11 is dedicated to representation theory of rational 
Cherednik algebras with a nonzero Planck constant. 
Namely, by analogy with semisimple Lie algebras, we develop 
the theory of category $\mathcal O$. In particular, we introduce
Verma modules, irreducible highest weight modules, which are
labeled by representations of $G$, and compute the
characters of the Verma modules. The main challenge is to compute
the characters of irreducible modules, and find out which of them
are finite dimensional. We do some of this in the case when 
$G=S_n\ltimes \Gamma^n$, where $\Gamma$ is a cyclic group. 
In particular, we construct and compute the characters of all 
the finite dimensional simple modules in the case $G=S_n$
(rational Cherednik algebra of type $A$). It turns out that 
a finite dimensional simple module exists if and only if the
parameter $k$ of the Cherednik algebra equals $r/n$, where 
$r$ is an integer relatively prime to $n$. For such values 
of $k$, such representation is unique, its dimension is
$|r|^{n-1}$, and it has no self-extensions. 


At the end of each lecture, we provide remarks and references, 
designed to put the material of the lecture in a broader
prospective, and link it with the existing literature. 
However, due to a limited size and scope of these lectures, 
we were, unfortunately, unable to give an exhaustive list of references 
on Calogero-Moser systems; such a list would have been truly enormous. 

{\bf Acknowledgments.} 
These lecture notes are dedicated to my mother Yelena Etingof on the
occasion of her 75-th birthday. Her continuous care from my early
childhood to this day has shaped me both as a person and as a
mathematician, and there are no words that are sufficient to
express my gratitude and admiration. 

I gave these lectures in the spring and summer of 
2005 at ETH (Zurich). I am greatly indebted to
Professor Giovanni Felder and the participants of his seminar for
being active listeners; they are responsible for the very
existence of these notes, as well as for improvement of their
quality. I am also very grateful to ETH for its hospitality
and financial support.  

This work was partially supported by the NSF grant DMS-0504847
and the CRDF grant RM1-2545-MO-03.

\section{Poisson manifolds and Hamiltonian reduction}

\subsection{Poisson manifolds}

Let $A$ be a commutative algebra over a field $k$. 

\begin{definition} We say that
$A$ is a Poisson algebra if it is equipped with a Lie bracket
$\lbrace{,\rbrace}$ such that $\lbrace{a,bc\rbrace}=
\lbrace{a,b\rbrace}c+b\lbrace{a,c\rbrace}$. 
\end{definition}

Let $I$ be an ideal in $A$. 

\begin{definition}
We say that $I$ is a Poisson ideal if
$\lbrace{A,I\rbrace}\subset I$. 
\end{definition}

In this case $A/I$ is a Poisson
algebra. 

Let $M$ be a smooth manifold. 

\begin{definition}
We say that $M$ is a Poisson
manifold if its structure algebra $C^\infty(M)$ is equipped  
with a Poisson bracket. 
\end{definition}

The same definition can be applied to 
complex analytic and algebraic varieties: a Poisson structure on
them is just a Poisson structure on the structure sheaf. 
Note that this definition may be used even for singular
varieties. 

\begin{definition}
A morphism of Poisson manifolds (=Poisson map) is a regular map 
$M\to N$ that induces a homomorphism of Poisson algebras
$C^\infty(N)\to C^\infty(M)$, i.e. a map that preserves Poisson
structure. 
\end{definition}

If $M$ is a smooth variety ($C^\infty$, analytic, or algebraic), 
then a Poisson structure on $M$ is defined by a Poisson
bivector $\Pi\in \Gamma(M,\wedge^2 TM)$ such that its 
Schouten bracket with itself is zero: $[\Pi,\Pi]=0$. 
Namely, $\lbrace{f,g\rbrace}:=(df\otimes dg)(\Pi)$
(the condition that $[\Pi,\Pi]=0$ is equivalent to the Jacobi
identity for $\lbrace{,\rbrace}$). 
In particular, if $M$ is symplectic (i.e. equipped with a closed
nondegenerate 2-form $\omega$) then it is Poisson with
$\Pi=\omega^{-1}$, and conversely, a Poisson manifold 
with nondegenerate $\Pi$ is symplectic with $\omega=\Pi^{-1}$. 

For any Poisson manifold $M$, we have a homomorphism of Lie
algebras $v: C^\infty(M)\to Vect_\Pi(M)$ from the Lie algebra 
of functions on $M$ to the Lie algebra of vector fields on $M$
preserving the Poisson structure, 
given by the formula $f\mapsto \lbrace{f,?\rbrace}$. 
In classical mechanics, one says that $v(f)$ is the Hamiltonian
vector field corresponding to the Hamiltonian $f$. 

\begin{exercise}\label{kvc} If $M$ is a connected symplectic manifold, then 
${\rm Ker}(v)$ consists of constant functions. 
If in addition $H^1(M,\Bbb C)=0$ then the map 
$v$ is surjective. 
\end{exercise}

\begin{example} $M=T^*X$, where $X$ is a smooth manifold. 
Define the Liouville 1-form $\eta$ on $T^*X$ as follows.
Let $\pi: T^*X\to X$ be the projection map. 
Then given $v\in T_{(x,p)}(T^*X)$, we set 
$\eta(v)=(d\pi\cdot v,p)$. Thus if $x_i$ are local coordinates on $X$ and
$p_i$ are the linear coordinates in the fibers of $T^*X$ with
respect to the basis $dx_i$ then $\eta=\sum p_idx_i$. 

Let $\omega=d\eta$. Then $\omega$ is a symplectic 
structure on $M$. In local coordinates, $\omega=\sum dp_i\wedge dx_i$. 
\end{example}

\begin{example}
 Let $\g$ be a finite dimensional Lie algebra. 
Let $\Pi: \g^*\to \wedge^2\g^*$ be the dual map to the Lie
bracket. Then $\Pi$ is a Poisson bivector on $\g^*$ 
(whose coefficients are linear). This Poisson structure 
on $\g^*$ is called the Lie Poisson structure. 

Let $\mathcal O$ be an orbit of the coadjoint action 
in $\g^*$. Then it is easy to check that 
the restriction of $\Pi$ to $\mathcal O$ is a 
section of $\wedge^2T\mathcal O$, which is nondegenerate. 
Thus $\mathcal O$ is a symplectic manifold. 
The symplectic structure on $\mathcal O$ is called the
Kirillov-Kostant structure. 
\end{example}

\subsection{Moment maps} 

Let $M$ be a Poisson manifold and $G$ a Lie group acting on $M$
by Poisson automorphisms. 
Let $\g$ be the Lie algebra of $G$. Then we have a homomorphism
of Lie algebras $\phi: \g\to Vect_\Pi(M)$. 

\begin{definition}
A $G$-equivariant regular map 
$\mu: M\to \g^*$ is said to be a moment map for the $G$-action on
$M$ if the pullback map $\mu^*: \g\to C^\infty(M)$ 
satisfies the equation $v(\mu^*(a))=\phi(a)$. 
\end{definition}

It is easy to see
that in this case $\mu^*$ is a homomorphism of Lie algebras, so 
$\mu$ is a Poisson map. Moreover, it is easy to show that 
if $G$ is connected then 
the condition of $G$-equivariance in the above definition can be
replaced by the condition that $\mu$ is a Poisson map.

A moment map does not always exist, and if it does, it is not 
always unique. However, if $M$ is a simply connected 
symplectic manifold, then the homomorphism 
$\phi: \g\to Vect_\Pi(M)$ can be lifted to a homomorphism 
$\hat\g\to C^\infty(M)$, where $M$ is a 1-dimensional 
central extension of $\g$. Thus there exists a moment map 
for the action on $M$ of the simply connected Lie group $\hat G$ 
corresponding to the Lie algebra $\hat\g$. In particular, if 
in addition the action of $G$ on $M$ is transitive, then $M$ is a
coadjoint orbit of $\hat G$. 

We also see that if $M$ is a connected
symplectic manifold then any two moment maps $M\to \g^*$ differ by 
shift by a character of $\g$. 
 
\begin{exercise} Show that if $M=\Bbb R^2$ with symplectic form
$dp\wedge dx$ and 
$G=\Bbb R^2$ acting by translations, 
then there is no moment map $M\to \g^*$. 
What is $\hat G$ in this case? 
\end{exercise}

\begin{exercise}
Show that if $M$ is simply connected and symplectic and 
$G$ is compact then there is a moment map $M\to \g^*$. 
\end{exercise}

\begin{exercise}\label{subme}
Show that if $M$ is symplectic then $\mu$ is a submersion near $x$ 
(i.e., the differential $d\mu_x: T_xM\to \g^*$ is surjective)
if and only if the stabilizer $G_x$ of $G$ is a discrete subgroup
of $G$ (i.e. the action is locally free near $x$). 
\end{exercise}

\begin{example}\label{mome}  Let $M=T^*X$, and let $G$ act on $X$. 
Define $\mu: T^*X\to \g^*$ by $\mu(x,p)(a)=p(\psi(a))$, $a\in
\g$, where $\psi: \g\to Vect(X)$ is the map defined by the action. 
 Then $\mu$ is a moment map. 
\end{example}

\subsection{Hamiltonian reduction}

Let $M$ be a Poisson manifold with an action of a Lie group $G$
preserving the Poisson structure, and with a moment map $\mu$. 
Then the algebra of $G$-invariants $C^\infty(M)^G$ is a Poisson algebra. 

Let $J$ be the ideal in $C^\infty(M)$ generated by 
$\mu^*(a)$, $a\in \g$. It is easy to see that 
$J$ is invariant under Poisson bracket with $C^\infty(M)^G$. 
Therefore, the ideal $J^G$ in $C^\infty(M)^G$ 
is a Poisson ideal, and hence the algebra $A:=C^\infty(M)^G/J^G$ is a Poisson
algebra. 

The geometric meaning of the algebra $A$ is as follows. 
Assume that the action of $G$ on $M$ is proper, i.e. 
for any two compact sets $K_1$ and $K_2$, the 
set of elements $g\in G$ such that $gK_1\cap K_2\ne \emptyset$ is 
compact. Assume also that the action of $G$ is free. 
In this case, the quotient $M/G$ is a manifold, and 
$C^\infty(M)^G=C^\infty(M/G)$. Moreover, as we mentioned 
in Exercise \ref{subme}, 
the map $\mu$ is a submersion (so $\mu^{-1}(0)$ is a smooth
submanifold of $M$), and the ideal $J^G$ corresponds to the 
submanifold $M//G:=\mu^{-1}(0)/G$ in $M$. Thus,
$A=C^\infty(M//G)$, and so $M//G$ is a Poisson manifold. 

\begin{definition}
The manifold $M//G$ is called the Hamiltonian reduction 
of $M$ with respect to $G$ using the moment map $\mu$. 
\end{definition}

\begin{exercise} Show that in this setting, if $M$ is symplectic, so is
$M//G$. 
\end{exercise} 

This geometric setting can be generalized in various directions. 
First of all, for $M//G$ to be a manifold, it suffices to require
that the action of $G$ be free only near $\mu^{-1}(0)$. Second, one can
consider a locally free action which is not necessarily free. 
In this case, $M//G$ is a Poisson orbifold. 

Finally, we can consider a purely algebraic setting 
which will be most convenient for us: $M$ is a scheme of
finite type over $\Bbb C$ (for example, a variety), and 
$G$ is an affine algebraic group. In this case, 
we do not need to assume that the action of $G$ is locally free 
(which allows us to consider many more examples). Still, 
some requirements are needed to ensure the existence of quotients. For example, 
a sufficient condition that often applies is that $M$ is an affine scheme
and $G$ is a reductive group. Then $M//G$ is an affine
Poisson scheme (possibly nonreduced and singular even if $M$ was
smooth). 

\begin{example}  Let $G$ act properly and freely on a manifold
$X$, and $M=T^*X$. Then $M//G$ (for the moment map as in Example 
\ref{mome}) 
is isomorphic to $T^*(X/G)$. 
\end{example}

On the other hand, the following example shows that when the
action of $G$ on an algebraic variety $X$ is
not free, the computation of the reduction $T^*X//G$ (as a
scheme) may be rather difficult. 

\begin{example}\label{twomat} Let $M=T^*{\rm Mat}_n(\Bbb C)$, and $G=PGL_n(\Bbb
C)$ (so $\g={\frak sl}_n(\Bbb C)$). 
Using the trace form we can identify $\g^*$ with $\g$, and 
$M$ with ${\rm Mat}_n(\Bbb C)\oplus {\rm Mat}_n(\Bbb C)$.
Then a moment map is given by the formula $\mu(X,Y)=[X,Y]$, for
$X,Y\in {\rm Mat}_n(\Bbb C)$. 
Thus $\mu^{-1}(0)$ is the {\bf commuting scheme} ${\rm Comm}(n)$ defined by the
equations $[X,Y]=0$, and the quotient $M//G$
is the quotient ${\rm Comm}(n)/G$, whose ring of functions is 
$A=\Bbb C[{\rm Comm}(n)]^G$. 

It is not known whether the commuting scheme is reduced (i.e. 
whether the corresponding ideal is a radical ideal);
this is a well known open problem. The underlying variety is
irreducible (as was shown by Gerstenhaber \cite{Ge1}), 
but very singular, and has a very complicated
structure. However, we have the following result. 

\begin{theorem} \label{gg} 
(Gan, Ginzburg,\cite{GG}) ${\rm Comm}(n)/G$ is reduced,
and isomorphic to $\Bbb C^{2n}/S_n$. Thus $A=\Bbb
C[x_1,...,x_n,y_1,...,y_n]^{S_n}$. 
The Poisson bracket on this algebra is induced from the standard
symplectic structure on $\Bbb C^{2n}$.  
\end{theorem}
\end{example}

{\bf Remark.} The hard part of this theorem is to show that 
${\rm Comm}(n)/G$ is reduced (i.e. $A$ has no nonzero nilpotent
elements).  

{\bf Remark.} Let $\g$ be a simple complex Lie algebra, and $G$
the corresponding group. The commuting scheme ${\rm Comm}(\g)$ 
is the subscheme of $\g\oplus \g$ defined by the equation
$[X,Y]=0$. Similarly to the above discussion, 
${\rm Comm}(\g)/G=T^*\g//G$. It is conjectured
that ${\rm Comm}(\g)$ and in particular 
${\rm Comm}(\g)/G$ is a reduced scheme; the latter is known for
$\g=sl(n)$ thanks to the Gan-Ginzburg theorem.  
It is also known that the underlying variety $\overline{{\rm
Comm}(\g)}$ is irreducible (as was shown by Richardson), and 
$\overline{{\rm Comm}(\g)}/G=(\h\oplus \h)/W$, where $\h$ is a
Cartan subalgebra of $\g$, and $W$ is the Weyl group of $\g$
(as was shown by Joseph \cite{J}).  

\subsection{Hamiltonian reduction along an orbit.} 

Hamiltonian reduction along an orbit is a generalization of the
usual Hamiltonian reduction. For simplicity let us describe it in the situation
when $M$ is an affine algebraic variety and $G$ a reductive group. 
Let $\O$ be a closed coadjoint orbit
of $G$,  $I_\O$ be the ideal in $S\g$ 
corresponding to $\O$, and let $J_\O$ be the ideal in $\Bbb C[M]$
generated by $\mu^*(I_\O)$. Then $J_\O^G$ is a Poisson ideal in 
$\Bbb C[M]^G$, and $A=\Bbb C[M]^G/J_\O^G$ is a Poisson algebra. 

Geometrically, ${\rm Spec}(A)=\mu^{-1}(\O)/G$ (categorical quotient). 
It can also be written as $\mu^{-1}(z)/G_z$, where $z\in \O$ and
$G_z$ is the stabilizer of $z$ in $G$. 

\begin{definition}
The scheme $\mu^{-1}(\O)/G$ is
called the Hamiltonian reduction of $M$ with respect to $G$ 
along $\mathcal O$. We will denote by $R(M,G,\O)$. 
\end{definition}

\begin{exercise}\label{symp} Show that if the action of $G$ on
$\mu^{-1}(\O)$ is free, and $M$ is a symplectic variety, then 
$R(M,G,\O)$ is a symplectic variety, of dimension 
$\dim(M)-2\dim(G)+\dim(\O)$. 
\end{exercise} 

\begin{exercise}\label{commu} (The Duflo-Vergne theorem,
\cite{DV}) Let $z\in \g^*$
be a generic element. Show that the connected component of the
identity of the group $G_z$ is commutative.

{\it Hint}: The orbit $\O$ of $z$ is described locally near $z$ by 
equations $f_1=...=f_m=0$, where $f_i$ are Casimirs of the
Poisson structure;
the Lie algebra $Lie(G_z)$ has basis $df_i(z)$, $i=1,...,m$.  
\end{exercise}

In a similar way, one can define Hamiltonian reduction along any Zariski
closed $G$-invariant subset of $\g^*$, for example the closure of
a non-closed coadjoint orbit. 

\subsection{Calogero-Moser space} \label{cm}

Let $M$ and $G$ be as in Example \ref{twomat}, and 
$\O$ be the orbit of the matrix ${\rm diag}(-1,-1,...,-1,n-1)$, 
i.e. the set of traceless matrices $T$ such that $T+1$ has rank $1$. 

\begin{definition} (Kazhdan, Kostant, Sternberg, \cite{KKS})
The scheme $\mC_n:=R(M,G,\O)$ is called the Calogero-Moser space. 
\end{definition}

Thus, $\mC_n$ is the space of conjugacy classes of pairs of
$n\times n$ matrices $(X,Y)$ such that the matrix $XY-YX+1$ has rank $1$. 

\begin{theorem}
The action of $G$ on $\mu^{-1}(\O)$ is free, and thus 
(by Exercise \ref{symp}) $\mC_n$ is a smooth 
symplectic variety (of dimension $2n$). 
\end{theorem} 

\begin{proof}
It suffices to show that if $X,Y$ 
are such that $XY-YX+1$ has rank 1, then 
$(X,Y)$ is an irreducible set of matrices. 
Indeed, in this case, by Schur's lemma,
if $B\in GL(n)$ is such that $BX=XB$ and $BY=YB$ 
then $B$ is a scalar, so the stabilizer of $(X,Y)$ in $PGL_n$ is
trivial. 

To show this, assume that $W\ne 0$ is an invariant subspace 
of $X,Y$. In this case, the eigenvalues 
of $[X,Y]$ on $W$ are a subcollection of the collection of $n-1$
copies of $-1$ and one copy of $n-1$. The sum of the elements of
this subcollection must be zero, since it is the trace of $[X,Y]$
on $W$. But then the subcollection must be the entire collection,
so $W=\Bbb C^n$, as desired. 
\end{proof} 

In fact, one also has the following more complicated theorem. 

\begin{theorem} (G. Wilson, \cite{Wi}) The Calogero-Moser space is
connected. 
\end{theorem}

\subsection{Notes} 1. For generalities on Poisson algebras and
Poisson manifolds, see e.g. the book \cite{Va}.
The basic material on symplectic manifolds can be found in the
classical book \cite{Ar}. Moment maps appeared already in the
works of S. Lie (1890). Classical Hamiltonian reduction (for
symplectic manifolds) was introduced by Marsden and Weinstein
\cite{MW}; reduction along an orbit is introduced by
Kazhdan-Kostant-Sternberg \cite{KKS}.  
For basics on moment maps and Hamiltonian 
reduction see e.g. the book \cite{OR}. 

2. The commuting variety and the Calogero-Moser space 
are very special cases of much more general Poisson varieties 
obtained by reduction, called quiver varieties; they play an
important role in geometric representation theory, and were
recently studied by many authors, notably W. Crawley-Boevey, G. Lusztig,
H. Nakajima. Wilson's connectedness theorem for the Calogero-Moser space 
can be generalized to quiver varieties; see \cite{CB}. 

\section{Classical mechanics and Integrable systems}

\subsection{Classical mechanics}

The basic setting of Hamiltonian classical mechanics is as
follows. The phase space of a mechanical system 
is a Poisson manifold $M$. The manifold $M$ is
usually symplectic and often equals $T^*X$, where $X$ is another
manifold called the configuration space. The dynamics of the
system is defined by its Hamiltonian (or energy function) $H\in C^\infty(M)$. 
Namely, the Hamiltonian flow attached to $H$ is the flow 
corresponding to the vector field $v(H)$. If $y_i$ are
coordinates on $M$, then the differential equations defining the
flow (Hamilton's equations) are written as
$$
\frac{dy_i}{dt}=\lbrace{H,y_i\rbrace},
$$
where $t$ is the time. 
The main mathematical problem in classical mechanics 
is to find and study the solutions of these equations. 

Assume from now on that $M$ is symplectic. Then by Darboux theorem 
we can locally choose coordinates $x_j,p_j$ on $M$ 
such that the symplectic form is $\omega=\sum dp_j\wedge dx_j$. 
Such coordinates are said to be canonical. For them one has 
$\lbrace{p_i,x_j\rbrace}=\delta_{ij},\lbrace{p_i,p_j\rbrace}=
\lbrace{x_i,x_j\rbrace}=0$. In canonical coordinates 
Hamilton's equations are written as 
$$
\frac{dx_i}{dt}=\frac{\partial H}{\partial p_i},\quad
\frac{dp_i}{dt}=-\frac{\partial H}{\partial x_i}.
$$

\begin{example}
Let $M=T^*X$, where $X$ is a Riemannian manifold.
Then we can identify $TX$ with $T^*X$ using the Riemannian
metric. Let $H=\frac{p^2}{2}+U(x)$, where $U(x)$ is a smooth function
on $X$ called the potential. 
The flow defined by this Hamiltonian describes the  
motion of a particle on $X$ in the potential field $U(x)$. 

Hamilton's equations for 
the Hamiltonian $H$ are 
$$
\dot{x}=p,\ \dot{p}=-\frac{\partial U}{\partial x},
$$
where $\dot{p}$ is the covariant time derivative of $p$
with respect to the Levi-Civita connection.
(If $X=\Bbb R^n$ with the usual metric, it is the usual time derivative). 
These equations reduce to Newton's equation 
$$
\ddot{x}=-\frac{\partial U}{\partial x}
$$
(again, $\ddot{x}$ is the covariant time derivative of $\dot{x}$).
If $U(x)$ is constant, the equation has the form 
$\ddot{x}=0$, which defines the so called geodesic flow. 
Under this flow, the particle moves along geodesics in $X$
(lines if $X=\Bbb R^n$) with constant speed. 
\end{example}

The law of conservation of energy, which trivially follows from 
Hamilton's equations, says that the Hamiltonian $H$ is constant
along the trajectories of the system. This means that if $M$ is
2-dimensional then the Hamiltonian flow is essentially a
1-dimensional flow along level curves of $H$, and thus 
the system can be solved explicitly in quadratures. 
In higher dimensions this is not the case, and a generic 
Hamiltonian system on a symplectic manifold of dimension $2n$,
$n>1$, has a very complicated behavior. 

\subsection{Symmetries in classical mechanics} 

It is well known that if a classical mechanical system has a
symmetry, then this symmetry can be used to reduce its order 
(i.e., the dimension of the phase space). For instance, suppose
that $F\in C^\infty(M)$ is a first integral of the system,
i.e., $F$ and $H$ are ``in involution'': 
$\lbrace{F,H\rbrace}=0$. Then $F$ gives rise to a (local) symmetry
of the system under the group $\Bbb R$ defined by the Hamiltonian
flow with Hamiltonian $F$\footnote{The symmetry is local 
because the solutions of the system defined by $F$ may not exist
for all values of $t$; this is, however, not essential for our 
considerations.}, and if $F$ is functionally independent
of $H$, then we can use $F$ to reduce the order of the system by 2.  
Another example is motion in a rotationally symmetric field 
(e.g., motion of planets around the sun), which can be 
completely solved using the rotational symmetry. 

The mathematical mechanism of using symmetry to reduce the order 
of the system is that of Hamiltonian reduction. 
Namely, let a Lie group $G$ act on $M$ preserving the
Hamiltonian $H$. Let $\mu: M\to \g^*$ be a moment map for this
action. For simplicity assume that $G$ acts freely on $M$
(this is not an essential assumption). 
We also assume that generically the Hamiltonian
vector field $v(H)$ is transversal to the $G$-orbits. 

It is easy to see that if $y=y(t)$ is a solution 
of Hamilton's equations, then $\mu(y(t))$ is constant. 
\footnote{This fact is a source of many conservation laws in physics,
e.g. conservation of momentum for translational symmetry, or
conservation of angular momentum for rotational symmetry;
it also explains the origin of the term ``moment map''.}

Therefore, the flow descends to the symplectic manifolds
$R(M,G,\O)=\mu^{-1}(\O)/G$, where $\O$ runs over orbits of
the coadjoint representation of $G$, with the same Hamiltonian. 
These manifolds have a smaller dimensions than that of $M$. 

On the other hand, suppose that $\O$ is a generic coadjoint
orbit, and we know the image $y_*(t)$ in 
$M_*:=R(M,G,\O)=\mu^{-1}(z)/G_z$ 
of a trajectory  $y(t)\in \mu^{-1}(z)$. Then we can find
$y(t)$ explicitly. To do so, let us denote 
the trajectory $y_*(t)$ by $T$, and let 
$T'$ be the preimage of $T$ in $\mu^{-1}(z)$ (so that $T=T'/G_z$).
Let us locally pick a cross-section 
$i: T\to T'$. Then we have (locally) an identification 
$T'=i(T)\times G_z=\Bbb R\times G_z$, and 
the vector field of the flow on $T'$ is given by the formula 
$\frac{\partial}{\partial t}+\gamma(t)$, where $\gamma(t)$ is a 
left-invariant vector field on $G_z$ (i.e. 
$\gamma(t)\in {\rm Lie}(G_z)$), which depends on the choice of
the cross section. Since the Lie algebra ${\rm Lie}(G_z)$
is abelian (Exercise \ref{commu}), the flow generated by this
vector field (and hence the trajectory $y(t)$) 
can be explicitly computed using Euler's formula
for solving linear first order differential equations. 

\subsection{Integrable systems}

\begin{definition} 
An integrable system on a symplectic manifold $M$ of dimension
$2n$ is a collection of smooth functions $H_1,...,H_n$ on $M$
such that they are in involution (i.e., 
$\lbrace{H_i,H_j\rbrace}=0$), and the differentials 
$dH_i$ are linearly independent on a dense open set in $M$. 
\end{definition}

Let us explain the motivation for this definition. 
Suppose that we have a Hamiltonian flow on $M$ with 
Hamiltonian $H$, and assume that $H$ can be included in an
integrable system: $H=H_1,H_2,...,H_n$. In this case, 
$H_1,H_2,...,H_n$ are first integrals of the flow. 
In particular, we can use $H_n$ to reduce the order of the flow 
from $2n$ to $2n-2$ as explained in the previous subsection. 
Since $\lbrace{H_i,H_n\rbrace}=0$, we get a Hamiltonian system 
on a phase space of dimension $2n-2$ with Hamiltonian $H$ and 
first integrals $H=H_1,...,H_{n-1}$, which form an integrable
system. Now we can use $H_{n-1}$ to further reduce the flow 
to a phase space of dimension $2n-4$, etc. Continuing in this
way, we will reduce the flow to the 2-dimensional phase space, 
As we mentioned before, such a flow can be integrated in
quadratures, hence so can the flow on the original manifold $M$. 
This motivates the terminology ``integrable system''. 

{\bf Remark.}  In a similar way one can define integrable
systems on complex analytic and algebraic varieties. 

\begin{exercise} Show that the condition that $dH_i$ are linearly independent 
on a dense open set is equivalent to the requirement 
that $H_i$ are functionally independent, i.e. 
there does not exist a nonempty open set $U\subset M$ 
such that the points $(H_1(u),...,H_n(u))$, $u\in U$, are
contained in a smooth hypersurface in $\Bbb R^n$. 
\end{exercise}

\begin{exercise} Let $H_1,...,H_n, H_{n+1}$ be a 
system of functions on a symplectic manifold $M$, such that
$\lbrace{H_i,H_j\rbrace}=0$. Let $x_0$ be a point of $M$ 
such that $dH_i$, $i=1,...,n$, are linearly independent at
$x_0$. Show that there exists neighborhood $U$ of $x_0$ and a
smooth function $F$ on a neighborhood of $y_0:=(H_1(x_0),...,H_n(x_0))\in \Bbb
R^n$ such that $H_{n+1}(x)=F(H_1(x),...,H_n(x))$ for $x\in U$. 
\end{exercise} 

\begin{theorem} (Liouville's theorem) Let $M$ be a symplectic
manifold, with an integrable system $H_1,...,H_n$. 
Let $c=(c_1,..,c_n)\in \Bbb R^n$, and $M_c$ be the 
set of points $x\in M$ such that $H_i(x)=c_i$. Assume that 
$dH_i$ are independent at every point of $M_c$, and that $M_c$ is
compact. In this case every connected component $N$ of $M_c$ is a
torus $T^n$, and the $n$ commuting flows corresponding to the 
Hamiltonians $H_i$ define an action of the group $\Bbb R^n$ on 
$N$, which identifies $N$ with $\Bbb R^n/\Gamma$, where $\Gamma$
is a lattice in $\Bbb R^n$.  
\end{theorem}

The proof of this theorem can be found, for instance, in 
\cite{Ar}, Chapter 10. 

{\bf Remark.} The compactness assumption of Liouville's theorem
is often satisfied. For example, it holds if $M=T^*X$, where $X$
is a compact manifold, and $H=\frac{p^2}{2}+U(x)$.  

{\bf Remark.} One may define a Hamiltonian $H$ on $M$ to be
{\it completely integrable} if it can be included in an integrable
system $H=H_1,H_2,...,H_n$. Unfortunately, 
this definition is not always satisfactory: for example,
it is easy to show that any Hamiltonian is completely integrable 
in a sufficiently small neighborhood of a point $x\in M$ where
$dH(x)\ne 0$. Thus when one says that a Hamiltonian is completely
integrable, one usually means not the above definition, but
rather that one can explicitly produce $H_2,...,H_n$ such that
$H=H_1,H_2,...,H_n$ is an integrable system (a property that is
hard to define precisely). However, there are situations where 
the above precise definition is adequate. For example, this is
the case in the algebro-geometric situation 
($M$ is an algebraic variety, Hamiltonians are polynomial
functions), or in the $C^\infty$-situation in the case when 
the hypersurfaces $H={\rm const}$ are compact. 
In both of these cases integrability is a strong restriction
for $n>1$. For instance, in the second case integrability implies that 
the closure of a generic trajectory has dimension $n$, 
while it is known that for ``generic'' $H$
with compact level hypersurfaces this closure is the entire 
hypersurface $H={\rm const}$, i.e., has dimension $2n-1$. 

\subsection{Action-angle variables of integrable systems}

To express solutions of Hamilton's equations corresponding to
$H$, it is convenient to use the so called ``action-angle
variables'', which are certain corrdinates on $M$ near a point
$P\in M$ where $dH_i$ are linearly independent. Namely, 
the action variables are simply the functions $H_1,...,H_n$. 
Note that since $H_i$ are in involution, they define 
commuting flows on $M$. 
To define the angle variables, let us fix (locally near $P$) a Lagrangian 
submanifold $L$ of $M$ (i.e., one of dimension $n$ on 
which the symplectic form vanishes), which is  
transversal to the joint level surfaces of $H_1,...,H_n$.  
Then for a point $z$ of $M$ sufficiently close to $P$
one can define the ``angle variables'' $\phi_1(z),...,\phi_m(z)$,
 which are the times one needs to move along the flows 
defined by $H_1,...,H_n$ starting from $L$ to reach the point
$z$. By the above explanations, 
$\phi_i$ can be found in quadratures if $H_i$ are known. 

The reason action-angle variables are useful is that 
in them both the symplectic form and Hamilton's equations 
have an extremely simple expression: the symplectic form is
$\omega=\sum dH_i\wedge d\phi_i$, and 
Hamilton's equations are given by the formula 
$\dot{H_i}=0$, $\dot{\phi_i}=\delta_{1i}$. 
Thus finding the action-angle variables 
is sufficient to solve Hamilton's equations 
for $H$. In fact, by producing the action-angle variables,
we also solve Hamilton's equations with Hamiltonians $H_i$, or, 
more generally, with Hamiltonian $F(H_1,...,H_n)$, where $F$ is
any given smooth function. 

\subsection{Constructing integrable systems by Hamiltonian
reduction.}\label{concl}

A powerful method of constructing integrable systems is
Hamiltonian reduction. Namely, let $M$ be a symplectic manifold,
and let $H_1,...,H_n$ be smooth functions on $M$ such that
$\lbrace{H_i,H_j\rbrace}=0$ and $dH_i$ are linearly independent
everywhere. Assume that $M$ carries a symplectic action of a Lie
group $G$ with moment map $\mu: M\to
\g^*$, which preserves the functions $H_i$, and let 
$\O$ be a coadjoint orbit of $G$. Assume that $G$ acts freely on
$\mu^{-1}(\O)$, and 
that the distribution ${\rm Span}(v(H_1),...,v(H_n))$ is transversal to
$G$-orbits along $\mu^{-1}(\O)$. 
In this case, the symplectic manifold $R(M,G,\O)$ carries a
collection of functions
$H_1,...,H_n$ such that $\lbrace{H_i,H_j\rbrace}=0$ and $dH_i$ 
are linearly independent everywhere. It sometimes happens that
while 
$n<\frac{1}{2}\dim M$ (so $H_1,...,H_n$ is {\bf not} an
integrable system on $M$), 
one has $n=\frac{1}{2}\dim R(M,G,\O)$, so that $H_1,...,H_n$ is
an integrable system 
on $R(M,G,\O)$. In this case the Hamiltonian flow defined by any
Hamiltonian
of the form $F(H_1,...,H_n)$ on $R(M,G,\O)$ can be solved in
quadratures. 

{\bf Remark.} For this construction to work, it suffices to
require 
that $dH_i$ are linearly independent on a dense open set of
$\mu^{-1}(\O)$. 

\subsection{The Calogero-Moser system} 

A vivid example of constructing an integrable system by
Hamiltonian reduction is the Kazhdan-Kostant-Sternberg
construction of the Calogero-Moser system. In this case
$M=T^*{\rm Mat}_n(\Bbb C)$ (regarded as the set of pairs of
matrices $(X,Y)$ as in Subsection \ref{cm}), with the usual
symplectic form $\omega={\rm Tr}(dY\wedge dX)$. Let
$H_i=\Tr(Y^i)$, $i=1,...n$. These functions are obviously in
involution, but they don't form an integrable system, because
there are too few of them. However,
let $G=PGL_n(\Bbb C)$ act on $M$ by conjugation, and let $\O$ be
the coadjoint orbit of $G$ considered in Subsection \ref{cm}
(consisting of traceless matrices $T$ such that $T+1$ has rank
$1$). Then the system $H_1,...,H_n$ descends to a system of
functions in involution on $R(M,G,\O)$, which is the
Calogero-Moser space ${\mathcal C}_n$. Since this space is
$2n$-dimensional, $H_1,...,H_n$ is an integrable system on
${\mathcal C}_n$. It is called the (rational) Calogero-Moser
system.

The Calogero-Moser flow is, by definition, the Hamiltonian flow
on ${\mathcal C}_n$ defined by the Hamiltonian
$H=H_2=\Tr(Y^2)$. Thus this flow is integrable, in the sense that
it can be included in an integrable system. In particular, its
solutions can be found in quadratures using the inductive
procedure of reduction of order. However (as often happens with
systems obtained by reduction), solutions can also be found by a
much simpler procedure, since they can be found already on the
``non-reduced'' space $M$: indeed, on $M$ the Calogero-Moser flow
is just the motion of a free particle in the space of matrices,
so it has the form $g_t(X,Y)=(X+2Yt,Y)$. The same formula is
valid on ${\mathcal C}_n$. In fact, we can use the same method to
compute the flows corresponding to all the Hamiltonians
$H_i=\Tr(Y^i)$, $i\in \Bbb N$: these flows are given by the
formulas $$ g_t^{(i)}(X,Y)=(X+iY^{i-1}t,Y).  $$

\subsection{Coordinates on ${\mathcal C}_n$ and the explicit form
of the Calogero-Moser system}

It seems that the result of our consideratons is trivial and
we've gained nothing. To see that this is, in fact, not so, let
us write the Calogero-Moser system explicitly in coordinates. To
do so, we first need to introduce local coordinates on the
Calogero-Moser space ${\mathcal C_n}$.

To this end, let us restrict our attention to the open set
$U_n\subset {\mathcal C}_n$ which consists of conjugacy classes
of those pairs $(X,Y)$ for which the matrix $X$ is
diagonalizable, with distinct eigenvalues; by Wilson's theorem,
this open set is dense in ${\mathcal C}_n$.

A point $P\in U_n$ may be represented by a pair $(X,Y)$ such that
$X={\rm diag}(x_1,...,x_n)$, $x_i\ne x_j$. In this case, the
entries of $T:=XY-YX$ are $(x_i-x_j)y_{ij}$. In particular, the
diagonal entries are zero. Since the matrix $T+1$ has rank $1$,
its entries $\kappa_{ij}$ have the form $a_ib_j$ for some numbers
$a_i,b_j$. On the other hand, $\kappa_{ii}=1$, so $b_j=a_j^{-1}$
and hence $\kappa_{ij}=a_ia_j^{-1}$. By conjugating $(X,Y)$ by
the matrix ${\rm diag}(a_1,...,a_n)$, we can reduce to the
situation when $a_i=1$, so $\kappa_{ij}=1$. Hence the matrix $T$
has entries $1-\delta_{ij}$ (zeros on the diagonal, ones off the
diagonal). Moreover, the representative of $P$ with diagonal $X$
and $T$ as above is unique up to the action of the symmetric
group $S_n$. Finally, we have $(x_i-x_j)y_{ij}=1$ for $i\ne j$,
so the entries of the matrix $Y$ are 
$y_{ij}=\frac{1}{x_i-x_j}$ if $i\ne j$. On the other hand, the
diagonal entries $y_{ii}$ 
of $Y$ are unconstrained. Thus
we have obtained the following result. 

\begin{proposition}
Let ${\Bbb C}^n_{\rm reg}$ be the open set 
of $(x_1,...,x_n)\in {\Bbb C}^n$ such that $x_i\ne x_j$ for $i\ne j$. 
Then there exists an isomorphism of algebraic varieties 
$\xi: T^*({\Bbb C}^n_{\rm reg}/S_n)\to U_n$ given by the formula 
$(x_1,...x_n,p_1,...,p_n)\to (X,Y)$, where $X={\rm
diag}(x_1,...,x_n)$, and $Y=Y(x,p):=(y_{ij})$, 
$$
y_{ij}=\frac{1}{x_i-x_j}, i\ne j,\ y_{ii}=p_i.
$$   
\end{proposition}

In fact, we have a stronger result: 

\begin{proposition}\label{symvar} $\xi$ is an isomorphism  
of symplectic varieties (where the cotangent bundle 
is equipped with the usual symplectic structure). 
\end{proposition}

For the proof of Proposition \ref{symvar}, 
we will need the following general and important but easy
theorem. 

\begin{theorem}\label{nbf} (The necklace bracket formula) 
Let $a_1,...,a_r$ and $b_1,...,b_s$ be either $X$ or $Y$. 
Then on $M$ we have 
\begin{equation}
\begin{aligned}
\lbrace{\Tr(a_1...a_r),\Tr(b_1...b_s)\rbrace}=\\
\sum_{(i,j): a_i=Y, b_j=X}
\Tr(a_{i+1}...a_ra_1...a_{i-1}b_{j+1}...b_sb_1...b_{j-1})-\\
\sum_{(i,j): a_i=X, b_j=Y}
\Tr(a_{i+1}...a_ra_1...a_{i-1}b_{j+1}...b_sb_1...b_{j-1}). 
\end{aligned}
\end{equation}
\end{theorem}

\begin{exercise} Prove Theorem \ref{nbf}.
\end{exercise} 

\begin{proof} (of Proposition \ref{symvar})
Let $a_k=\Tr(X^k)$, $b_k=\Tr(X^kY)$.
It is easy to check using the Necklace bracket formula 
that on $M$ we have 
$$
\lbrace{a_m,a_k\rbrace}=0,
\lbrace{b_m,a_k\rbrace}=ka_{m+k-1},
\lbrace{b_m,b_k\rbrace}=(k-m)b_{m+k-1}.
$$
On the other hand, $\xi^*a_k=\sum x_i^k$, 
$\xi^*b_k=\sum x_i^kp_i$. Thus we see that 
$$
\lbrace{f,g\rbrace}=\lbrace{\xi^*f,\xi^*g\rbrace},
$$ 
where $f,g$ are either $a_k$ or $b_k$. 
On the other hand, the functions $a_k,b_k$, $k=0,...,n-1$, form a local
coordinate system near a generic point of $U_n$, so we are done. 
\end{proof} 

Now let us write the Hamiltonian of the Calogero-Moser
system in coordinates. It has the form 
\begin{equation}\label{cmeq}
H=\Tr(Y(x,p)^2)=\sum_i p_i^2-\sum_{i\ne j}\frac{1}{(x_i-x_j)^2}.
\end{equation}
Thus the Calogero-Moser Hamiltonian describes the motion 
of a system of $n$ particles on the line with interaction
potential $-1/x^2$. 
This is the form of the Calogero-Moser Hamiltonian in which it
originally occured in the work of F. Calogero. 

Now we finally see the usefulness of the Hamiltonian reduction
procedure. The point is that it is not clear at all 
from formula (\ref{cmeq}) why the Calogero-Moser Hamiltonian
should be completely integrable. However, our reduction procedure
implies the complete integrability of $H$, and gives an explicit
formula for the first integrals: 
$$
H_i=\Tr(Y(x,p)^i).
$$
Moreover, this procedure immediately 
gives us an explicit solution of the system.
Namely, assume that $x(t),p(t)$ is the solution 
with initial condition $x(0),p(0)$.
Let $(X_0,Y_0)=\xi(x(0),p(0))$.  
Then $x_i(t)$ are the eigenvalues of the matrix 
$X_t:=X_0+2tY_0$, and $p_i(t)=x_i'(t)/2$. 

{\bf Remark.} In fact, the reduction procedure not only allows 
us to solve the Calogero-Moser system, but also 
provides a ``partial compactification'' of its phase 
space $T^*{\Bbb C}^n_{\rm reg}$, namely the Calogero-Moser space 
${\mathcal C}_n$, to which the Calogero-Moser flow smoothly
extends. I think it is fair to say that ${\mathcal C}_n$ is ``the
right'' phase space for the Calogero-Moser flow. 

\begin{exercise}
Compute explicitly the integral $H_3$. 
\end{exercise}

\subsection{The trigonometric Calogero-Moser system}

Another integrable system which can be obtained by a similar
reduction procedure is the trigonometric Calogero-Moser system. 
To obtain it, take the same $M,G,\O$ as in the case of the
rational Calogero-Moser system, but define 
$H_i^*:=\Tr(Y_*^i)$, where $Y_*=XY$. 

\begin{proposition}\label{invo}
The functions $H_i^*$ are in involution. 
\end{proposition}

\begin{exercise}
Deduce Proposition \ref{invo} from  
the necklace bracket formula. 
\end{exercise}

Thus the functions $H_1^*,...,H_n^*$ define an integrable system
on the Calogero-Moser space ${\mathcal C}_n$. 

The trigonometric Calogero-Moser system is defined by the
Hamiltonian $H=H_2^*$. 
Let us compute it more explicitly on the open set $U_n$ using the
coordinates $x_i,p_i$. We get $Y_*=(y_{*ij})$, 
where 
$$
y_{*ij}=\frac{x_i}{x_i-x_j}, i\ne j,\quad y_{*ii}=x_ip_i. 
$$
Thus we have 
$$
H=\sum_i (x_ip_i)^2-\sum_{i\ne j}\frac{x_ix_j}{(x_i-x_j)^2}.
$$

Let us introduce ``additive'' coordinates $x_{i*}=\log x_i$, 
$p_{i*}=x_ip_i$. It is easy to check that these coordinates are
canonical. In them, the trigonometric Calogero-Moser Hamiltonian 
looks like 
$$
H=\sum_i p_{i*}^2-\sum_{i\ne j}\frac{4}{{\rm sinh}^2((x_{i*}-x_{j*})/2)},
$$
This Hamiltonian describes the motion of a system of 
$n$ particles on the line with interaction potential 
$-4/{\rm sinh}^2(x/2)$. 

{\bf Remark.} By replacing $x_{j*}$ with $ix_{j*}$, we can 
also integrate the system with Hamiltonian 
$$
H=\sum_i p_{i*}^2+\sum_{i\ne j}\frac{4}{{\rm sin}^2((x_{i*}-x_{j*})/2)},
$$
which corresponds to a system of 
$n$ particles on the circle of length $2\pi$  with interaction potential 
$+4/{\rm sin}^2(x/2)$. 

\subsection{Notes} 1. For generalities on classical mechanics,
symmetries of a mechanical system, reduction of order using
symmetries, integrable systems, action-angle variables we refer
the reader to \cite{Ar}. Classical Calogero-Moser systems go back to the
papers \cite{Ca}, \cite{Mo}; their construction using reduction along orbit 
is due to Kazhdan, Kostant, and Sternberg, \cite{KKS}. 

2. The necklace bracket formula is a starting point of
noncommutative symplectic geometry; it appears in \cite{Ko3}. 
This formula was generalized to the case of quivers in \cite{BLB},\cite{Gi}. 

\section{Deformation theory} 

Before developing the quantum analogs of the 
notions and results of lectures 1 and 2, 
we need to discuss the general theory of quantization 
of Poisson manifolds. We start with an even more general
discussion -- the deformation theory of associative algebras. 

\subsection{Formal deformations of associative algebras}

Let $k$ be a field, and $K=k[[\hbar_1,...,\hbar_n]]$. 
Let ${\mathfrak{m}}=(\hbar_1,...,\hbar_n)$ be the maximal ideal in
$K$. 

\begin{definition} A topologically free 
$K$-module is a $K$-module isomorphic to $V[[\hbar_1,...,\hbar_n]]$ for some
vector space $V$ over $k$. 
\end{definition}

Let $A_0$ be an algebra\footnote{In these lectures, by an algebra
we always mean an associative algebra with unit} over $k$. 

\begin{definition} A (flat) formal
$n$-parameter deformation of $A_0$ is an  algebra $A$ 
over $K$ which is topologically free as a $K$-module, together
with an algebra isomorphism $\eta_0: A/{\mathfrak{m}} A\to A_0$. 
\end{definition}

When no confusion is possible, we will call $A$ a
deformation of $A_0$ (omitting ``formal''). 

Let us restrict ourselves to one-parameter deformations with
parameter $\hbar$. Let us choose an identification $\eta: A\to A_0[[\hbar]]$ as 
$K$-modules, such that $\eta=\eta_0$ modulo $\hbar$. Then the
product in $A$ is completely determined by the product of
elements of $A_0$, which has the form of a ``star-product''
$$
\mu(a,b)=a\ast b=\mu_0(a,b)+\hbar \mu_1(a,b)+\hbar^2 \mu_2(a,b)+...,
$$
where $\mu_i: A_0\otimes A_0\to A_0$ are linear maps, and
$\mu_0(a,b)=ab$. 

\subsection{Hochschild cohomology}
The main tool in deformation theory of associative algebras is Hochschild cohomology.
Let us recall its definition. 

Let $A$ be an associative algebra. 
Let $M$ be a bimodule over $A$.  
A Hochschild $n$-cochain of $A$ with coefficients in $M$ 
is a linear map $A^{\otimes n}\to M$. 
The space of such cochains is denoted by $C^n(A,M)$. 
The differential $d:C^n(A,M)\to C^{n+1}(A,M)$ is defined 
by the formula 
$$
df(a_1,...,a_{n+1})=f(a_1,...,a_n)a_{n+1}-f(a_1,...,a_na_{n+1})
$$
$$
+f(a_1,...,a_{n-1}a_n,a_{n+1})-...+(-1)^nf(a_1a_2,...,a_{n+1})+
$$
$$
(-1)^{n+1}a_1f(a_2,...,a_{n+1}).
$$
It is easy to show that $d^2=0$. 

\begin{definition}
The Hochschild cohomology $H^\bullet(A,M)$ is defined  
 to be the cohomology of the complex 
$(C^\bullet(A,M),d)$. 
\end{definition}

\begin{remark} The Hochschild cohomology is often denoted by 
$HH^\bullet(A,M)$, but we'll use the shorter notation $H^\bullet$. 
\end{remark}

\begin{proposition}
One has a natural isomorphism 
$$
H^i(A,M)\to {\rm Ext}^i_{A-{\rm bimod}}(A,M),
$$ 
where $A-{\rm bimod}$ denotes the category of $A$-bimodules. 
\end{proposition}

\begin{proof}
The proof is obtained immediately by considering the bar resolution 
of the bimodule $A$:
$$
...\to A\otimes A\otimes A\to A\otimes A\to A,
$$
where the bimodule structure on $A^{\otimes n}$ is given by 
$$
b(a_1\otimes a_2\otimes ...\otimes a_n)c=ba_1\otimes
a_2\otimes...\otimes a_nc,
$$
and the map $\partial_n: A^{\otimes n}\to A^{\otimes {n-1}}$ is given by the formula 
$$
\partial_n(a_1\otimes a_2\otimes...\otimes a_n)=a_1a_2\otimes...\otimes a_n-...
+(-1)^{n}a_1\otimes... \otimes a_{n-1}a_n.
$$ 
\end{proof}

Note that we have the associative 
Yoneda product 
$$
H^i(A,M)\otimes H^j(A,N)\to H^{i+j}(A,M\otimes_A N),
$$ 
induced by tensoring of cochains.   

If $M=A$, the algebra itself, then we will 
denote $H^\bullet(A,M)$ by $H^\bullet(A)$. 
For example, $H^0(A)$ is the center of $A$, and 
$H^1(A)$ is the quotient of the Lie algebra of derivations of $A$ by inner derivations.  
The Yoneda product induces a graded algebra 
structure on $H^\bullet(A)$; it can be shown that 
this algebra is supercommutative. 

\subsection{Hochschild cohomology and deformations}

Let $A_0$ be an algebra, and let us look 
for 1-parameter deformations $A=A_0[[\hbar]]$ of $A_0$. 
Thus, we look for such series $\mu$ which satisfy the associativity 
equation, modulo the automorphisms of the $k[[\hbar]]$-module $A_0[[\hbar]]$
which are the identity modulo $\hbar$. 
\footnote{Note that we don't have to worry about the existence of a unit in 
$A$ since a formal deformation of an algebra with unit always has a 
unit.}  
 
The associativity equation 
$\mu\circ(\mu\otimes Id)=\mu\circ (Id\otimes \mu)$ 
reduces to a hierarchy of linear equations: 
\begin{equation}\label{asso}
\sum_{s=0}^N \mu_s(\mu_{N-s}(a,b),c)=
\sum_{s=0}^N \mu_s(a,\mu_{N-s}(b,c)).
\end{equation}
(These equations are linear in $\mu_N$ if $\mu_i$, $i<N$, are known).

To study these equations, one can use Hochschild cohomology. 
Namely, we have the following are standard facts
(due to Gerstenhaber, \cite{Ge2}), which can be
checked directly.  

1.  The linear equation for $\mu_1$ says that $\mu_1$ is a Hochschild 
2-cocycle. Thus algebra structures on $A_0[\hbar]/\hbar^2$ 
deforming $\mu_0$ are parametrized by the space $Z^2(A_0)$ of 
Hochschild 2-cocycles of $A_0$ with values in $M=A_0$.  

2. If $\mu_1,\mu_1'$ are two 2-cocycles such that 
$\mu_1-\mu_1'$ is a coboundary, then the algebra structures 
on $A_0[\hbar]/\hbar^2$ corresponding to $\mu_1$ and $\mu_1'$ are 
equivalent by a transformation of $A_0[\hbar]/\hbar^2$ 
that equals the identity 
modulo $\hbar$, and vice versa. Thus equivalence classes 
of multiplications on $A_0[\hbar]/\hbar^2$ deforming $\mu_0$
are parametrized by the cohomology $H^2(A_0)$. 

3. The linear equation for $\mu_N$ says that $d\mu_N$ 
is a certain quadratic expression $b_N$ in $\mu_0,\mu_1,...,\mu_{N-1}$. 
This expression is always a Hochschild 3-cocycle, and the equation is 
solvable iff it is a coboundary. Thus the cohomology class 
of $b_N$ in $H^3(A_0)$ is the only obstruction to solving this equation. 

\subsection{Universal deformation}

In particular, if $H^3(A_0)=0$ then the equation for $\mu_n$ can be solved 
for all $n$, and for each $n$ the freedom in choosing the solution, 
modulo equivalences, is the space $H:=H^2(A_0)$. Thus there exists 
an algebra structure over $k[[H]]$ on the space $A_u:=A_0[[H]]$ 
of formal functions from $H$ to $A_0$, $a,b\mapsto \mu_u(a,b)\in A_0[[H]]$,
($a,b\in A_0$), such that $\mu_u(a,b)(0)=ab\in A_0$, and 
every 1-parameter flat formal deformation $A$ of $A_0$ is given 
by the formula $\mu(a,b)(\hbar)=\mu_u(a,b)(\gamma(\hbar))$
for a unique formal series $\gamma\in \hbar H[[\hbar]]$,
with the property that $\gamma'(0)$ is the cohomology class 
of the cocycle $\mu_1$. 

Such an algebra $A_u$ is called a {\bf universal deformation} of $A_0$. 
It is unique up to an isomorphism (which may involve an
automorphism of 
$k[[H]]$). 

Thus in the case $H^3(A_0)=0$, deformation theory allows us 
to completely classify 1-parameter flat formal deformations of $A_0$.
In particular, we see that the ``moduli space'' parametrizing formal deformations of $A_0$
is a smooth space -- it is the formal neighborhood of zero in $H$.

If $H^3(A_0)$ is nonzero then in general the universal 
deformation parametrized by $H$ does not exist, as there are obstructions
to deformations. In this case, the moduli space of deformations 
will be a closed subscheme of the formal neighborhood of zero in 
$H$, which is often singular. 
On the other hand, even when $H^3(A_0)\ne 0$, the universal 
deformation parametrized by $H$ may exist (although 
it may be more difficult to prove than in the vanishing case).
In this case one says that the deformations of $A_0$ are {\bf unobstructed}
(since all obstructions vanish even though the space of obstructions doesn't).

\subsection{Quantization of Poisson algebras and manifolds} 

An example when there are obstructions to deformations is the
theory of quantization of Poisson algebras and manifolds. 

Let $M$ be a smooth $C^\infty$-manifold or a 
smooth affine algebraic variety over $\Bbb C$, 
and $A_0={\mathcal O}(M)$ the structure algebra of $M$.

{\bf Remark.} In the $C^\infty$-case, we will consider only local maps 
$A_0^{\otimes n}\to A_0$, i.e. those given by polydifferential operators, and 
all deformations and the Hochschild cohomology 
is defined using local, rather than general, cochains. 

\begin{theorem}(Hochschild-Kostant-Rosenberg) \cite{HKR} One has
$H^i(A_0)=\Gamma(M,\wedge^iTM)$ as a module over $A_0=H^0(A_0)$. 
\end{theorem}

In particular, $H^2$ is the space of bivector fields, and $H^3$
the space of trivector fields. So the cohomology class of $\mu_1$
is a bivector field; in fact, it is $\frac{1}{2}\pi$, where
$\pi(a,b):=\mu_1(a,b)-\mu_1(b,a)$, since any 2-coboundary in this
case is symmetric. The equation for $\mu_2$ says that $d\mu_2$ is
a certain expression that depends quadratically on $\mu_1$. It
is easy to show that the cohomology class of 
this expression is the Schouten bracket
$\frac{1}{4}[\pi,\pi]$. Thus, for the existence of $\mu_2$ it is necessary
that $[\pi,\pi]=0$, i.e. that $\pi$ be a {\bf Poisson bracket}.
In other words, the trivector field $[\pi,\pi]\in H^3(A_0)$ is an obstruction
to extending the first order deformation $a\ast b=ab+\hbar
\mu_1(a,b)$ to higher orders. 

More generally, let $A_0$ be any commutative algebra, 
and $A=A_0[[\hbar]]$ be a not necessarily commutative (but
associative) deformation of $A_0$. In this case, 
$A_0$ has a natural Poisson structure, given by the formula 
$\lbrace{a,b\rbrace}=[a',b']/\hbar\text{ mod }\hbar$, where
$a',b'$ are any lifts of $a,b$ to $A$. It is easy to check that
this expression is independent on the choice of the lifts. 
In terms of the star-product, this bracket is given by the
formula $\lbrace{a,b\rbrace}=\mu_1(a,b)-\mu_1(b,a)$.  

\begin{definition} In this situation, $(A,\mu)$ is said to be {\bf a quantization} of
$(A_0,\lbrace{,\rbrace})$, and 
$(A_0,\lbrace{,\rbrace})$ 
is said to be the {\bf quasiclassical limit} of
$(A,\mu)$. 
\end{definition} 

{\bf Remark.} If $A_0={\mathcal O}(M)$, one says that $A$ is a
quantization of $M$, and $M$ the quasiclassical limit of $A$. 

This raises the following important question. 
Suppose $A_0={\mathcal O}(M)$. Given a Poisson bracket $\pi$ on $M$, 
is it always possible to construct its quantization?

By the above arguments, $\mu_2$ exists (and a choice of $\mu_2$
is unique up to adding an arbitrary bivector). So there arises
the question of existence of $\mu_3$ etc., i.e. the question
whether there are other obstructions.

The answer to this question is yes and no. Namely, if you don't
pick $\mu_2$ carefully, you may be unable to find $\mu_3$, but
you can always pick $\mu_2$ so that $\mu_3$ exists, and there is
a similar situation in higher orders. This subtle fact is a
consequence of the following deep theorem of Kontsevich:

\begin{theorem}\label{kon} 
\cite{Ko1,Ko2} Any Poisson structure $\pi$ on $A_0$ can be quantized. 
Moreover, there is a natural bijection between products $\mu$ 
up to an isomorphism equal to $1$ modulo $\hbar$, 
and Poisson brackets $\pi_0+\hbar \pi_1+\hbar^2\pi_2+...$ up to a
formal diffeomorphism equal to $1$ modulo $\hbar$, 
such that the quasiclassical limit of $\mu$ is $\pi_0$.
\end{theorem}

{\bf Remarks.} 1. The Kontsevich deformation quantization has an
additional property called locality: the maps $\mu_i(f,g)$ are
differential operators with respect to both $f$ and $g$. 

2. Note that, as was shown by O. Mathieu, 
a Poisson bracket on a general commutative $\Bbb C$-algebra may fail 
to admit a quantization. 

Let us consider the special case of symplectic manifolds, 
i.e. the case when $\pi$ is a nondegenerate bivector. 
In this case we can consider $\pi^{-1}=\omega$, which is a closed, nondegenerate 2-form 
(=symplectic structure) on $M$. In this case, Kontsevich's theorem  
is easier, and was proved by De Wilde - 
Lecomte, and later Deligne and Fedosov.
Moreover, in this case there is the 
following additional result, also due to Kontsevich, \cite{Ko1,Ko2}. 

\begin{theorem} If $M$ is symplectic and $A$ is a quantization of $M$, then 
the Hochschild cohomology $H^i(A[\hbar^{-1}])$ is isomorphic to 
$H^i(M,\Bbb C((\hbar)))$.  
\end{theorem} 

{\bf Remark.} Here the algebra $A[\hbar^{-1}]$ is regarded as a (topological) 
algebra over the field of Laurent series 
$\Bbb C((\hbar))$, so Hochschild cochains are, by definition, 
linear maps $A_0^{\otimes n}\to A_0((\hbar))$.

\begin{example}\label{quansym}
The algebra $B=A[\hbar^{-1}]$ provides an example of an algebra 
with possibly nontrivial $H^3(B)$, for which the universal 
deformation parametrized by $H=H^2(B)$ exists. 
Namely, this deformation is attached through the correspondence of Theorem \ref{kon}
(and inversion of $\hbar$)
to the Poisson bracket
$\pi=(\omega+t_1\omega_1+...+t_r\omega_r)^{-1}$, 
where $\omega_1,...,\omega_r$ 
are closed 2-forms on $M$ which represent a basis of $H^2(M,\Bbb
C)$, and $t_1,...,t_r$ are the 
coordinates on $H$ corresponding to this basis. 
\end{example}

\subsection{Algebraic deformations}

Formal deformations of algebras often arise from algebraic
deformations. The most naive definition of an algebraic deformation 
is as follows. 

Let $\Sigma$ be a smooth affine algebraic curve over $k$
(often $\Sigma=k$ or $\Sigma=k^*$), and let $0\in \Sigma(k)$. 
Let $I_0$ be the maximal ideal corresponding to $0$. 

\begin{definition} An algebraic deformation over $B:=k[\Sigma]$ of an algebra $A_0$
is a $B$-algebra $A$ which is a free
$B$-module, together with the identification $\eta_0: A/I_0A\to A_0$
of the zero-fiber of $A$ with $A_0$ as algebras.   
\end{definition}

Any algebraic deformation defines a formal deformation. 
Indeed, let $\hbar$ be a formal parameter of $\Sigma$ 
near $0$, and let $\widehat A$ be the completion of $A$ with respect
to $I_0$ (i.e., $\widehat A={\rm limproj}_{n\to
\infty}A/I_0^nA$). Then $\widehat A$ is a topologically free 
$k[[\hbar]]$-module which is a deformation of $A_0$. 

\begin{definition}
An algebraic quantization of a Poisson algebra $A_0$ is 
an algebraic deformation $A$ of $A_0$ such that 
the completion $\widehat{A}$ is a deformation quantization 
of $A_0$. 
\end{definition}

\begin{example} 1. (Weyl algebra) 
$A_0=k[x,p]$ with the usual Poisson
structure, $\Sigma=k$, 
$A=k[\hbar,x,\hbar\partial_x]$.

2. (Generalization of 1) Let $\overline{A}$ be a filtered
algebra: $k=F^0\overline{A}\subset F^1\overline{A}\subset...$,
$\cup_i F^i\overline{A}=\overline{A}$. Assume that ${\rm
gr}\overline{A}=A_0$. Assume that $A_0$ is commutative. 
Then $A_0$ has a natural Poisson structure of degree $-1$
(why?). In this case, an algebraic quantization of $A_0$ is 
given by the so called {\it Rees algebra} $A$ of $\overline{A}$. 
Namely, the algebra $A={\rm Rees}(\overline{A})$ is defined by
the formula $A=\oplus_{n=0}^\infty F^n\overline{A}$. 
This is an algebra over $k[\hbar]$, where $\hbar$ is the element 1 
of the summand $F^1\overline{A}$. It is easy to see 
that $A$ is an algebraic deformation of $A_0$ (with 
$\Sigma=k$). 

An important example of this: $X$ is a manifold, 
$\overline{A}=D(X)$, the algebra of differential operators 
on $X$, $F^\bullet$ is the filtration by order. 
In this case $A_0=C^\infty_{pol}(T^*X)$, the space of 
smooth functions on $T^*X$ which are polynomial 
along fibers (of uniformly bounded degree), with the usual Poisson
structure. 

Another important example: $\g$ is a Lie algebra, 
$\overline{A}=U(\g)$ is the universal enveloping algebra of $\g$,
$F^\bullet$ is the natural filtration on $U(\g)$ which is defined
by the condition that $\deg(x)=1$ for $x\in \g$. 
In this case $A_0$ is the Poisson algebra $S\g$ with the Lie
Poisson structure, and $A=U(\g_\hbar)$, where $\g_\hbar$ 
is the Lie algebra over $\Bbb C[\hbar]$ which is equal to 
$\g[\hbar]$ as a vector space, with bracket
$[a,b]_\hbar:=\hbar [a,b]$. 

3. (Quantum torus)
$A_0=k[x^{\pm 1},y^{\pm 1}]$, with Poisson bracket
$\lbrace{x,y\rbrace}=xy$, $A=k[q,q^{-1}]<x^{\pm 1},y^{\pm
1}>/(xy=qyx)$, $\Sigma=k^*$. 
\end{example}

\subsection{Notes} 1. For generalities on Hochschild cohomology, 
see the book \cite{Lo}. The basics of deformation theory 
of algebras are due to Gerstenhaber \cite{Ge2}. 

2. The notion of deformation quantization was proposed 
in the classical paper \cite{BFFLS}; in this paper the authors
ask the question whether every Poisson manifold admits a
quantization, which was solved positively by Kontsevich.  

\section{Quantum moment maps, quantum 
Hamiltonian reduction, and the Levasseur-Stafford theorem} 

\subsection{Quantum moment maps and quantum Hamiltonian reduction}

Now we would like to quantize the notion of a moment map. 
For simplicity let us work over $\Bbb C$.
Recall that a classical moment map is defined for a Poisson manifold 
$M$ with a Poisson action of a Lie group $G$, as a Poisson map $M\to
\g^*$ whose components are Hamiltonians defining the action of
$\g=Lie(G)$ on $M$. In algebraic terms, a Poisson manifold $M$ with 
a $G$-action defines a Poisson algebra $A_0$ (namely,
$C^\infty(M)$) together with a Lie algebra map $\phi_0: \g\to {\rm
Der}_{\Pi}(A_0)$ from $\g$ to the Lie algebra of derivations of 
$A_0$ preserving its Poisson bracket. A classical moment map is then a
homomorphism of Poisson algebras $\mu_0: S\g\to A_0$ such that 
for any $a\in S\g$, $b\in A_0$ one has
$\lbrace{\mu_0(a),b\rbrace}=\phi_0(a)b$. 

This algebraic reformulation makes it perfectly clear how one
should define a quantum moment map. 

\begin{definition} Let $\g$ be a Lie algebra, and 
$A$ be an associative algebra equipped with a $\g$-action, 
i.e. a Lie algebra map $\phi: \g\to {\rm Der}A$. 

(i) A quantum moment map for $(A,\phi)$ is an associative algebra
homomorphism $\mu: U(\g)\to A$ such that for any $a\in \g$, $b\in
A$ one has $[\mu(a),b]=\phi(a)b$. 

(ii) Suppose that $A$ is a filtered associative algebra, such that 
${\rm gr}A$ is a Poisson algebra $A_0$, equipped with a $\g$-action $\phi_0$
and a classical moment map $\mu_0$. Suppose that 
${\rm gr}\phi=\phi_0$. A quantization of 
$\mu_0$ is a quantum moment map 
$\mu: U(\g)\to A$ such that ${\rm gr}\mu=\mu_0$. 

(iii) More generally, suppose that $A$  is a deformation quantization 
of a Poisson algebra $A_0$ equipped with a $\g$-action $\phi_0$
and a classical moment map $\mu_0$. Suppose that 
$\phi=\phi_0 {\rm mod }\hbar$. 
A quantization of 
$\mu_0$ is a quantum moment map 
$\mu: U(\g)\to A[\hbar^{-1}]$ such that 
for $a\in \g$ we have $\mu(a)=\hbar^{-1}\mu_0(a)+O(1)$. 
\end{definition}

Thus, a quantum moment map is essentially a homomorphism of Lie
algebras $\mu: \g\to A$. 
Note that like in the classical case, the action $\phi$ is
determined by the moment map $\mu$. 

\begin{example} Let $X$ be a manifold with an action of a Lie
group $G$, and $A=D(X)$ be the algebra of differential operators
on $X$. There is a natural homomorphism $\mu: \g\to {\rm
Vect}X\subset D(X)$ which is a quantum moment map 
for the natural action of $\g$ on $A$. It is a quantization  
of the classical moment map for the action of $G$ on $T^*X$
defined in Section 1. 
\end{example}

\subsection{Quantum hamiltonian reduction} 

The algebraic definition of Hamiltonian reduction 
given in Section 1 is easy to translate to the quantum
situation. Namely, let $A$ be an algebra with a $\g$-action and a
quantum moment map $\mu: U(\g)\to A$. The space of
$\g$-invariants $A^\g$, i.e. elements $b\in A$ such that
$[\mu(a),b]=0$ for all $a\in \g$, is a subalgebra of $A$.
Let $J\subset A$ be the left ideal generated by $\mu(a)$, $a\in
\g$. Then $J$ is not a 2-sided ideal, but $J^\g:=J\cap A^\g$ is a
2-sided ideal in $A^\g$. 

Indeed, let $c\in A^\g$, and $b\in
J^\g$, $b=\sum_i b_i\mu(a_i)$, $b_i\in A,a_i\in \g$. Then 
$bc=\sum b_i\mu(a_i)c=\sum b_ic\mu(a_i)\in J^\g$. 

Thus, the algebra $A//\g:=A^\g/J^\g$ is an associative algebra,
which is called the quantum Hamiltonian reduction of $A$
with respect to the quantum moment map $\mu$.
\footnote{If the Lie algebra $\g$ is not reductive, or its action
on $A$ is not locally finite, this definition may be too naive
to give good results. In this case, 
instead of taking $\g$-invariants (which is not an exact
functor) one may need to include all its derived functors, i.e. cohomology
($Ext^i$). Moreover, if $\g$ is infinite dimensional (Virasoro,
affine Kac-Moody), and $A$ is a vertex algebra rather than a usual
associative algebra, which is an important case in string theory, then
one needs to consider ``semi-infinite'' cohomology,
i.e. $i=\infty/2+j$, $j\in \Bbb Z$.}

An easy example of quantum Hamiltonian reduction is given by the
following exercise. 

\begin{exercise} \label{qred} Let $X$ be a smooth affine algebraic variety 
with a free action of a connected reductive algebraic group $G$. 
Let $A=D(X)$, and $\mu: \g\to {\rm Vect}X\to D(X)$ 
be the usual action map. Show that $A//\g=D(X/G)$.  
\end{exercise}

\subsection{The Levasseur-Stafford theorem}
In general, similarly to the classical case,  it is rather
difficult to compute the quantum reduction $A//\g$. 
For example, in this subsection we will 
describe $A//\g$ in the case when $A=D(\g)$ is the algebra of differential operators 
on a reductive Lie algebra $\g$, and $\g$ acts on $A$ through the adjoint 
action on itself. This description is a very nontrivial result
of Levasseur and Stafford.

Let $\h$ be a Cartan subalgebra of $\g$, and $W$ 
the Weyl group of $(\g,\h)$. 
Let $\h_{reg}$ denote the set of regular points in $\h$, 
i.e. the complement of the reflection hyperhplanes. 
To describe $D(\g)//\g$, we 
will construct a homomorphism $HC: D(\g)^\g\to D(\h)^W$, 
called the Harish-Chandra homomorpism (as it was first
constructed by Harish-Chandra). 
Recall that we have the classical Harish-Chandra isomorphism
$\zeta: \Bbb C[\g]^\g\to\Bbb C[\h]^W$, defined simply by restricting
$\g$-invariant functions on $\g$ to the Cartan subalgebra $\h$. 
Using this isomorphism, we can define an action of $D(\g)^\g$ on
$\Bbb C[\h]^W$, which is clearly given by $W$-invariant 
differential operators. 
However, these operators will, in general, have poles on the
reflection hyperplanes. Thus we get a homomorphism 
$HC': D(\g)^\g\to D(\h_{\rm reg})^W$. 

The homomorphism $HC'$ is called the radial part homomorphism, as
for example for $\g={\frak{su}}(2)$ it computes the radial parts of
rotationally invariant differential operators on $\Bbb R^3$ 
in spherical coordinates. This homomorphism is not yet
what we want, since it does not actually land in $D(\h)^W$ (the
radial parts have poles). 

Thus we define the Harish-Chandra homomorphism by twisting 
$HC'$ by the discriminant
$\delta(x)=\prod_{\alpha>0}(\alpha,x)$
($x\in \h$, and $\alpha$ runs over positive roots of $\g$):
$$
HC(D):=\delta \circ HC'(D)\circ \delta^{-1}\in D(\h_{\rm reg})^W.
$$

\begin{theorem} \label{ls} (i) (Harish-Chandra, \cite{HC}) 
For any reductive $\g$, $HC$ lands in $D(\h)^W\subset
D(\h_{reg})^W$.  

(ii) (Levasseur-Stafford \cite{LS})  
The homomorphism $HC$ defines an isomorphism $D(\g)//\g=D(\h)^W$.  
\end{theorem}

{\bf Remarks.} 1. Part (i) of the theorem says that the poles magically disappear after 
conjugation by $\delta$.

2. Both parts of this theorem are quite nontrivial. 
The first part was proved by Harish-Chandra using analytic methods, and 
the second part by Levasseur and Stafford using the theory of D-modules. 

In the case $\g={\frak{gl}}_n$, Theorem \ref{ls}
is a quantum analog of Theorem \ref{gg}. 
The remaining part of this subsection is devoted to the proof   
of Theorem \ref{ls} in this special case, using Theorem \ref{gg}. 

We start the proof with the following proposition, valid for any reductive Lie algebra. 

\begin{proposition}\label{concoe}
If $D\in (S\g)^\g$ is a differential
operator with constant coefficients, then $HC(D)$ is 
the $W$-invariant differential operator with constant
coefficients on $\h$, obtained from $D$ via the classical
Harish-Chandra isomorphism $\eta: (S\g)^\g\to (S\h)^W$. 
\end{proposition}

\begin{proof} Without loss of generality, we may assume that $\g$ is simple.

\begin{lemma}\label{lapl} 
Let $D$ be the Laplacian $\Delta_\g$ of $\g$, corresponding 
to an invariant form. Then $HC(D)$ is the Laplacian $\Delta_\h$. 
\end{lemma}

\begin{proof}
let us calculate $HC'(D)$. 
We have 
$$
D=\sum_{i=1}^r \partial_{x_i}^2+2\sum_{\alpha>0}\partial_{f_\alpha}\partial_{e_\alpha}, 
$$
where $x_i$ is an orthonormal basis of $\h$, and $e_\alpha,f_\alpha$ are root elements such that 
$(e_\alpha,f_\alpha)=1$. 
Thus if $F(x)$ is a $\g$-invariant function on $\g$, then we get 
$$
(DF)|_{\h}=\sum_{i=1}^r \partial_{x_i}^2(F|_\h)+2\sum_{\alpha>0}(\partial_{f_\alpha}\partial_{e_\alpha}F)|_{\h}. 
$$
Now let $x\in \h$, and consider 
$(\partial_{f_\alpha}\partial_{e_\alpha}F)(x)$. We have 
$$
(\partial_{f_\alpha}\partial_{e_\alpha}F)(x)=
\partial_s\partial_t|_{s=t=0}F(x+tf_\alpha+se_\alpha).
$$
On the other hand, we have 
$$
{\rm Ad}(e^{s\alpha(x)^{-1}e_\alpha})(x+tf_\alpha+se_\alpha)=
x+tf_\alpha+ts\alpha(x)^{-1}h_\alpha+...,
$$
where $h_\alpha=[e_\alpha,f_\alpha]$. Hence, 
$$
\partial_s\partial_t|_{s=t=0}F(x+tf_\alpha+se_\alpha)=
\alpha(x)^{-1}(\partial_{h_\alpha}F)(x). 
$$
This implies that 
$$
HC'(D)F(x)=\Delta_\h F(x)+2\sum_{\alpha>0}\alpha(x)^{-1}\partial_{h_\alpha}F(x).
$$
Now the statement of the Lemma reduces to the identity 
$$
\delta^{-1}\circ \Delta_\h \circ \delta=\Delta_\h+2\sum_{\alpha>0}\alpha(x)^{-1}\partial_{h_\alpha}. 
$$
This identity follows immediately from the identity
$$
\Delta_\h \delta=0.
$$
To prove the latter, it suffices to note that $\delta$ is the lowest degree 
nonzero polynomial on $\h$, which is antisymmetric under the action of $W$. 
The lemma is proved. 
\end{proof} 

Now let $D$ be any element of $(S\g)^\g\subset D(\g)^\g$ of degree $d$
(operator with constant coefficients). It is obvious that the leading 
order part of the operator $HC(D)$ is the operator $\eta(D)$ 
with constant coefficients, whose symbol is just the 
restriction of the symbol of $D$ from $\g^*$ to $\h^*$. 
Our job is to show that in fact $HC(D)=\eta(D)$. To do so, denote by 
$Y$ the difference $HC(D)-\eta(D)$. Assume $Y\ne 0$. By Lemma \ref{lapl}, the operator $HC(D)$ commutes with 
$\Delta_\h$. Therefore, so does $Y$. Also $Y$ has homogeneity degree $d$ but order $m\le d-1$.
Let $S(x,p)$ be the symbol of $Y$ ($x\in \h,p\in \h^*$). Then $S$  is a homogeneous function of homogeneity degree $d$
under the transformations $x\to t^{-1}x$, $p\to tp$, polynomial in $p$ of degree $m$. 
From these properties of $S$ it is clear that $S$ is not a polynomial (its degree in $x$ is $m-d<0$). 
On the other hand, since $Y$ commutes with $\Delta_\h$, the Poisson bracket of $S$ with $p^2$ is zero.    
Thus Proposition \ref{concoe} follows from the following lemma. 

\begin{lemma}\label{brac}
Let $S: (x,p)\mapsto S(x,p)$ 
be a rational function on 
$\h\oplus \h^*$, which is polynomial in $p\in \h^*$. 
Let $f: \h^*\to \Bbb C$ be a polynomial, 
such that there is no vector $v\in \h^*$ for which $df(p)(v)$ is
identically zero in $p$ (for example, $f=p^2$). 
Suppose that the Poisson bracket $\lbrace{f ,S\rbrace}$ 
equals to zero. Then $S$ is a polynomial: $S\in \Bbb C[\h\oplus \h^*]$.  
\end{lemma}

\begin{proof} (R. Raj). Let $x_0\in \h$ be a generic point in the divisor
of poles of $S$. Then the function $S^{-1}$ is regular and vanishes
at $(x_0,p)$ for generic $p$. Also, we have
$\lbrace{S^{-1},f\rbrace}=0$, which implies that $S^{-1}$
vanishes along the entire flowline of the Hamiltonian flow
defined by $f$. This flowline is defined by the formula 
$x(t)=x_0+t\cdot df(p), p(t)=p$, and must be contained in the
pole divisor of $S$ at $x_0$. This implies that the vector 
$df(p)$ is tangent to the pole divisor of $S$ at $x_0$ for almost
every $p$ (in the sense of Zariski topology). 
Thus if $v\in \h^*$ is the normal to the pole divisor of $S$
at $x_0$, then $df(p)(v)=0$ for almost every, hence every $p$. 
This is a contradiction, which implies that $x_0$ does not exist,
and hence $S$ is a polynomial. 
\end{proof} 

Thus Proposition \ref{concoe} is proved. 
\end{proof}

Now we continue the proof of Theorem \ref{ls}. Consider the
filtration on $D(\g)$ in which
$\deg(\g)=1, \deg(\g^*)=0$ 
(the order filtration), and the associated graded map 
${\rm gr}HC: \Bbb C[\g\times \g^*]^\g\to \Bbb C[\h_{\rm reg}\times
\h^*]^W$, which attaches to every differential operator
the symbol of its radial part. 
It is easy to see that this map is just the restriction map to
$\h\oplus \h^*\subset \g\oplus \g^*$, so it actually lands in 
$\Bbb C[\h\oplus \h^*]^W$. 

Moreover, ${\rm gr}HC$ is a map {\bf onto} $\Bbb C[\h\oplus \h^*]^W$. 
Indeed, ${\rm gr}HC$ is a Poisson map, so the surjectivity follows 
from the following Lemma. 

\begin{lemma}\label{poisgen}
$\Bbb C[\h\oplus \h^*]^W$ is generated as a
Poisson algebra by $\Bbb C[\h]^W$ and $\Bbb C[\h^*]^W$, i.e. by 
functions $f_m=\sum x_i^m$ and $f_m^*=\sum p_i^m$, $m\ge 1$. 
\end{lemma}

\begin{proof}
For the proof we need the following theorem due to H. Weyl (from his book
``Classical groups'').

\begin{theorem}\label{Wey} Let $B$ be an algebra over $\Bbb C$. 
Then the algebra $S^nB$ is generated by elements of the form 
$$
b\otimes 1\otimes...\otimes 1+
1\otimes b\otimes...\otimes 1+...+
1\otimes 1\otimes...\otimes b.
$$
\end{theorem}

\begin{proof} Since $S^nB$ is linearly spanned by elements of the
form $x\otimes...\otimes x$, $x\in B$, it suffices to prove the 
theorem in the special case $B=\Bbb C[x]$. In this case, 
the result is simply the fact that the ring of symmetric
functions is generated by power sums, which is well known. 
\end{proof} 

Now, we have $\lbrace{f_m^*,f_r\rbrace}=mr\sum x_i^{r-1}p_i^{m-1}$, and by
Weyl's theorem (applied to $B=\Bbb C[x,y]$, 
such functions  generate $\Bbb C[\h\oplus
\h^*]^W$ as an algebra). The lemma is proved. 
\end{proof} 

Let $K_0$ be the kernel of ${\rm gr}HC$. Then by 
Theorem \ref{gg}, $K_0$ is the ideal of the commuting scheme
${\rm Comm}(\g)/G$. 

Now consider the kernel $K$ of the homomorphism $HC$. 
It is easy to see that $K\supset J^\g$, so ${\rm gr}(K)\supset
{\rm gr}(J)^\g$. On the other hand, since $K_0$ is the ideal of
the commuting scheme, we clearly have ${\rm
gr}(J)^\g\supset K_0$, and $K_0\supset {\rm gr}K$. 
This implies that $K_0={\rm gr}K={\rm gr}(J)^\g$, and $K=J^\g$.

It remains to show that ${\rm Im}HC=D(\h)^W$. 
Since ${\rm gr}K=K_0$, we have 
${\rm gr}{\rm Im}HC=\Bbb C[\h\oplus \h^*]^W$. 
Therefore, to finish the proof of the Harish-Chandra and Levasseur-Stafford theorems, 
it suffices to prove the following proposition. 

\begin{proposition}\label{conta}
${\rm Im}HC\supset D(\h)^W$. 
\end{proposition}

\begin{proof} 
We will use the following Lemma. 

\begin{lemma} \label{wal} (N. Wallach, \cite{Wa})
$D(\h)^W$ is generated as an algebra 
by $W$-invariant functions and $W$-invariant differential
operators with constant coefficients. 
\end{lemma}
 
\begin{proof} 
The lemma follows by taking associated graded algebras from 
Lemma \ref{poisgen}.
\end{proof}

\begin{remark}
Levasseur and Stafford showed \cite{LS} that this lemma is valid 
for any finite group $W$ acting on a finite dimensional vector
space $\h$. However, the above proof does not apply, since, as
explained in \cite{Wa}, Lemma \ref{poisgen} fails for many groups $W$, 
including Weyl groups of exceptional Lie
algebras $E_6,E_7,E_8$ (in fact it even fails for the cyclic group of order
$>2$ acting on a 1-dimensional space!). 
The general proof is more complicated and
uses some results in noncommutative algebra. 
\end{remark}

Lemma \ref{wal} and Proposition \ref{concoe} imply Proposition \ref{conta}.
\end{proof}

Thus, Theorem \ref{ls} is proved. 

\begin{exercise}
Let $\g_{\Bbb R}$ be the compact form of $\g$, and
$\O$ a regular coadjoint orbit in $\g_{\Bbb R}^*$. Consider the function 
$$
\psi_\O(x)=\int_\O e^{(b,x)}db,\ x\in \h,
$$
where $db$ is the measure on the orbit coming from the Kirillov-Kostant
symplectic structure. Prove the Harish-Chandra formula 
$$
\psi_\O(x)=\delta^{-1}(x)\sum_{w\in
W}(-1)^{l(w)}e^{(w\lambda,x)},
$$
where $\lambda$ is the intersection of $\O$ with the 
dominant chamber in the dual Cartan subalgebra $\h_{\Bbb R}^*\subset
\g_{\Bbb R}^*$, and $l(w)$ is the length of an element $w\in W$.

Deduce from this the Kirillov character formula for finite
dimensional representations \cite{Ki}:
If $\lambda$ is a dominant integral weight, and $L_\lambda$ is 
the corresponding representation of $G$, then 
$$
{\rm Tr}_{L_\lambda}(e^x)=
\frac{\delta(x)}{\delta_{tr}(x)}\int_{\O_{\lambda+\rho}} e^{(b,x)}db,
$$
where $\delta_{tr}(x)$ is the trigonometric version of
$\delta(x)$, i.e. the Weyl denominator
$\prod_{\alpha>0}(e^{\alpha(x)/2}-e^{-\alpha(x)/2})$, $\rho$ is
the half-sum of positive roots, and $\O_\mu$ denotes the
coadjoint orbit passing through $\mu$.

{\bf Hint.} Use Proposition \ref{concoe} and the fact that $\psi_\O$ 
is an eigenfunction of $D\in (S\g)^\g$ with eigenvalue $\chi_\O(D)$, where
$\chi_\O(D)$ is the value of the invariant polynomial $D$ at
the orbit $\O$. 
\end{exercise}

\subsection{Hamiltonian reduction with respect to an ideal in $U(\g)$.}

In Section 1, we defined Hamiltonian reduction along an orbit
or, more generally, closed $G$-invariant subset of $\g^*$. Algebraically 
a closed $G$-invariant subset of $\g^*$ corresponds to a Poisson
ideal in $S\g$. Thus the quantum analog of this construction
should be quantum Hamiltonian reduction with respect 
to a two-sided ideal $I\subset U(\g)$.  

The Hamiltonian reduction with respect to $I\subset U(\g)$ is
defined as follows. Let $\mu: U(\g)\to A$ be a quantum moment
map, and $J(I)=A\mu(I)\subset A$. Then $J(I)^\g$ is a 2-sided
ideal in $A^\g$ (this is shown similarly to the case of the usual
quantum hamiltonian reduction), and we set
$R(A,\g,I):=A^\g/J(I)^\g$. 
The usual quantum
reduction described above is the special case 
of this, when $I$ is the augmentation ideal.
 
The following example is a quantization of the
Kazhdan-Kostant-Sternberg construction 
of the Calogero-Moser space given in Section 1. 

\begin{example} \label{quantumcm}
Let $\g={\frak{gl}}_n$, $A=D(\g)$ as above. 
Let $k$ be a complex number, and 
$W_k$ be the representation of ${\frak{sl}}_n$ 
on the space of functions of the form 
$(x_1...x_n)^kf(x_1,...,x_n)$, where $f$ is a Laurent polynomial
of degree $0$. We regard $W_k$ as a $\g$-module by pulling it
back to $\g$ under the natural projection 
$\g\to {\frak{sl}}_n$. Let $I_k$ be the annihilator of $W_k$ in
$U(\g)$. 
\end{example}

Let $B_k=R(A,\g,I_k)$. Then $B_k$ has a filtration 
induced from the order filtration of $D(\g)^\g$. 
Let $HC_k: D(\g)^\g\to B_k$ be the natural homomorphism, and 
$K(k)$ be the kernel of $HC_k$. 

\begin{theorem} (Etingof-Ginzburg, \cite{EG}) (i) $K(0)=K$, $B_0=D(\h)^W$,
$HC_0=HC$. 

(ii) ${\rm gr}K(k)={\rm Ker}{\rm gr}HC_k=K_0$ for all complex $k$.
Thus, $HC_k$ is a flat family of homomorphisms. 

(iii) The algebra ${\rm gr}B_k$ is commutative and isomorphic to
$\Bbb C[\h\oplus \h^*]^W$ as a Poisson algebra.  
\end{theorem}

\begin{definition} The algebra $B_k$ is called the spherical rational
Cherednik algebra. The homomorphism $HC_k$ is called the deformed
Harish-Chandra homomorphism. 
\end{definition}

The algebra $B_k$ has an important representation
on the space $\Bbb C[\h_{reg}]$ of 
regular functions on the set $\h_{reg}$ of diagonal matrices
with different eigenvalues. To construct it, note that 
the algebra $B_k$ acts naturally on the space $E_k$
of $\g$-equivariant functions with values in $W_k$ on the formal
neighborhood of $\h_{reg}$ in $\g_{reg}$ 
(where $\g_{reg}$ is the set of all elements in $\g$ conjugate to 
an element of $\h_{reg}$, i.e. the set of matrices with different 
eigenvalues). Namely, the algebra $D(\g)^\g$ acts on $E_k$, and 
the ideal $I_k$ is clearly annihilated under this action. 

Now, an equivaraint function on the formal neighborhood of $\h_{reg}$ 
with values in $W_k$ is completely determined by its values 
on $\h_{reg}$ itself, and the only restriction for these values is
that they are in $W_k[0]$, the zero-weight subspace of $W_k$. 
But the space $W_k[0]$ is 1-dimensional, spanned by the function 
$(x_1,...,x_n)^k$. Thus the space $E_k$ is naturally isomorphic to 
$\Bbb C[\h_{reg}]$. Thus we have defined an action of $B_k$ on 
$\CC[\h_{reg}]$.  

The algebra $B_k$ is one of the main objects of this course. 
It is, in an appropriate sense, a quantization of
the Calogero-Moser space. We will discuss the precise formulation
of this statement later. 

\subsection{Quantum reduction in the deformational setting}

In the deformation theoretical setting, 
the construction of the quantum Hamiltonian reduction
should be slightly modified. Namely, assume that 
$A$ is a deformation quantization of an algebra $A_0=C^\infty(M)$, 
where $M$ is a Poisson manifold, and $\mu_0:\g\to A_0$ is a classical moment
map. Suppose that $\mu$ is a quantum moment map which is 
a quantization of $\mu_0$. Let $\O\subset \g^*$ be an invariant closed set, 
and $R(M,G,\O)$ be the corresponding classical reduction. 
 
Let $\g_\hbar$ be the Lie algebra over $\Bbb C[[\hbar]]$, which is
$\g[[\hbar]]$ as a vector space, with Lie bracket 
$[a,b]_\hbar=\hbar [a,b]$ ($a,b\in \g$). Let $U(\g_\hbar)$ 
be the enveloping algebra of $\g_\hbar$; it is a deformation 
quantization of $S\g$ (which is a completion of the Rees algebra of $U(\g)$).
We have a modified quantum moment map $\mu_\hbar: U(\g_\hbar)\to A$ 
given by $\mu_\hbar(a)=\hbar\mu(a)$ for $a\in \g$. 
 
Let $I\subset U(\g_\hbar)$ be an ideal deforming the ideal $I_0\subset S\g$ 
of functions vanishing on $\O$. Note that $I$ does not necessarily 
exist, but it does exist in many important cases, e.g. 
when $\O$ is a semisimple orbit of a reductive Lie algebra. Then we define 
the quantum reduction $R(A,\g,I):=A^\g/(A\mu_\hbar(I))^\g$
(quotient by an $\hbar$-adically closed
ideal). 

It is clear that the algebra $R(A,\g,I)$ is a deformation 
of the function algebra on $R(M,G,\O)$, but this deformation may not be flat. 
If it is flat (which happens in nice cases), one says that ``quantization 
commutes with reduction''. 

\subsection{Notes} 1. Quantum moment maps and quantum reduction 
have been considered for more than 20 years by many authors, in
particular, in physics literature; this notion arises naturally
when one considers BRST quantization of gauge theories. 
A convenient reference for us is the paper \cite{Lu}, which
considers a more general setting of quantum group actions. 
 
2. The ``quantization commutes with reduction'' conjecture was
formulated by Guillemin and Sternberg in \cite{GS}. 
In the case of a compact Lie group action and deformation
quantization, it was proved by Fedosov \cite{Fe}. 

\section{Quantum mechanics, quantum integrable systems, and quantization 
of the Calogero-Moser system}

\subsection{Quantum mechanics}
Before explaining the basic setting of quantum mechanics, let us 
present the basic setting of classical mechanics in a form convenient 
for quantization. 
In classical mechanics, we have the 
algebra of observables $C^\infty(M)$, where $M$ is the phase space 
(a Poisson, usually symplectic, manifold) 
and the motion of a point $x=x(t)\in M$ is described by Hamilton's 
equation for observables:
$$
\dot{f}=\lbrace{H,f\rbrace},
$$ 
where $f=f(x(t))$ is an observable
$f\in C^\infty(M)$, evaluated at $x(t)$. 

Similarly, the basic setting of Hamiltonian quantum mechanics is as follows. 
We have a (noncommutative) algebra $A$ of quantum observables, which acts 
(faithfully) in a space of states ${\mathcal H}$ (a complex Hilbert space). 
The Hamiltonian is an element $H$ of $A$ 
(i.e. an operator on ${\mathcal H}$, self-adjoint, and usually unbounded). 
The dynamics of the system $\psi=\psi(t)$, $\psi\in H$, 
is governed by the Schr\"odinger equation 
$$
\dot{\psi}=-\frac{iH\psi}{h},
$$
where $h>0$ is the Planck constant. 
Solutions of the Schr\"odinger equation have the form 
$$
\psi(t)=U(t)\psi(0),
$$
where $U(t)$ is the evolution operator $e^{-iHt/h}$.
This means that if $F\in A$ is any observable then 
its observed value at the state $\psi(t)$ is given by the formula
$$
<\psi(t)|F|\psi(t)>=<\psi(0)|F(t)|\psi(0)>, 
$$
where $F(t):=e^{iHt/h}Fe^{-iHt/h}$. The operator valued function
$F(t)$ satisfies the Heisenberg equation
$$
\dot{F}(t)=i[H,F(t)]/h.
$$

Quantum systems considered in quantum mechanics are usually deformations 
of certain classical mechanical systems, which are recovered when the Planck 
constant goes to zero. To study this limit, 
it is convenient to introduce \footnote{This notation is 
sacrilegious, as in all textbooks on quantum mechanics $h$ and $\hbar$ 
differ by $2\pi$}$\hbar=-ih$, and assume that $\hbar$ 
is a formal parameter rather than a numerical constant. 
(In physics such approach is called ``perturbation theory''). 
More specifically, we assume that the algebra $A$ of quantum observables 
is a formal deformation of the algebra of classical observables
$A_0=C^\infty(M)$ (i.e., $A$ is a deformation quantization of the 
phase space $M$), and the Hamiltonian $H\in A$ is a deformation of 
a classical Hamiltonian $H_0\in A_0$. In this case, we see that 
 the Heisenberg equation is  
nothing but a deformation of the Hamilton's equation for observables.  
Thus this framework can indeed be used to regard quantum mechanics 
as a deformation of classical mechanics. 

In many situations, $M=T^*X$, where $X$ is a Riemannian manifold, 
$H_0=\frac{p^2}{2}+U(x)$, where $U(x)$ is a potential. 
In this case, in the classical setting we can restrict ourselves to 
fiberwise polynomial functions: $A_0=C^\infty_{pol}(X)$. 
Then in the quantum setting we can work over polynomials 
$\Bbb C[\hbar]$ rather than formal series $\Bbb C[[\hbar]]$
(which is good since we can then specialize $\hbar$ to a numerical value),
and we have: $A={\rm Rees}(D(X))$, and
$H=-\frac{h^2}{2}\Delta+U(x)$, where $\Delta$ is the Laplacian.
Then $A$ is a quantization of $A_0=C^\infty_{pol}(T^*X)$, and $H$ is a 
quantization of $H_0$ (the minus sign comes from the fact that $i$ 
has been absorbed into $\hbar$). Finally, the space ${\mathcal H}$ is 
$L^2(X)$; the algebra $A/(\hbar=-ih_0)=D(X)$, $h_0\in \Bbb R_+$ acts 
in this space and is the algebra of quantum observables.

With these conventions, the Schr\"odinger equation takes the form 
$$
ih\frac{\partial \psi}{\partial t}=-\frac{h^2}{2}\Delta\psi+U(x)\psi,
$$
which is the classical form of the Schr\"odinger equation
(for a particle of unit mass). 

\subsection{Quantum integrable systems}

The above interpretation of quantum mechanics as a deformation of classical 
mechanics motivates the following definition of a quantum integrable system. 

Let $A_0=C^\infty(M)$, where $M$ is a symplectic manifold 
(or $A_0=\Bbb C[M]$, where $M$ is a symplectic affine algebraic variety
\footnote{If $M=T^*X$, we will often take the subalgebra 
$A_0=C^\infty_{pol}(T^*X)$ 
instead of the full algebra $C^\infty(X)$.}). Assume that $M$ has dimension $2n$. Let
$A$ be a quantization of $A_0$ (formal or algebraic). 

\begin{definition} A quantum integrable system 
in $A$ is a pairwise commuting system of elements $H_1,...,H_n$ 
such that their reductions $H_{1,0},...,H_{n,0}$ to $A_0$ 
form a classical integrable system on $M$. 
\end{definition}

In this situation, we say that the quantum integrable system $H_1,...,H_n$
is a quantization of the classical integrable system 
$H_{1,0},...,H_{n,0}$, and conversely, the classical system is the 
quasiclassical limit of the quantum system. 
Also, if we have a quantum 
mechanical system defined by a Hamiltonian $H$ which is included 
in a quantum integrable system $H=H_1,...,H_n$ then 
we say that $H_1,...,H_n$ are quantum integrals of $H$. 

{\bf Remark.} It is obvious that if $H_1,...,H_n$ are a quantum 
integrable system, 
then they are algebraically independent. 

\begin{example}\label{Tstar} 
$M=T^*X$, $A_0=C^\infty_{pol}(X)$, $A={\rm Rees}(D(X))
\subset D(X)[\hbar]$. Thus a quantum integrable system in $A$, upon 
evaluation $\hbar\to 1$, is just a collection of commuting 
differential operators $H_1,...,H_n$ on $X$.  
\end{example}

The motivation for this definition is the same as in the classical case. 
Namely, if $H=H_1,H_2,...,H_n$ is an integrable system, then 
the Schr\"odinger equation $\dot{\psi}=-H\psi/\hbar$ can usually be solved
``explicitly''. 

More specifically, recall that  
to solve the Schr\"odinger equation, it is sufficient 
to find an eigenbasis $\psi_m$ for $H$. Indeed, 
if $\lambda_m$ are the eigenvalues of $H$ at $\psi_m$, and 
$\psi(0)=\sum_m a_m\psi_m$, $a_m\in \Bbb C$, then 
$\psi(t)=\sum_m a_m e^{-i\lambda_m t/h}\psi_m$
(the index $m$ may be continuous, and then instead of the sum we will 
have to use integral). 

But in the problem of finding an eigenbasis, it is very useful 
to have operators commuting with $H$: then we can replace 
the eigenvalue problem $H_1\psi=\lambda_1\psi$ with the
joint eigenvalue problem 
\begin{equation}\label{jep}
H_i\psi=\lambda_i\psi, i=1,...,n,
\end{equation}  
whose solutions are easier to find. 

To illustrate this, look at the situation of Example \ref{Tstar}.
An integrable system in this example is an algebraically independent 
collection of differential operators $H_1,...,H_n$ on $X$. 
Let $S_i(x):=S_x(H_i)\in \Bbb C[T_x^*X]=S(T_xX)$ 
be the symbols of the operators $H_i$ (homogeneous polynomials), and assume 
that for any $x\in X$ 
the algebra $\Bbb C[T_x^*X]$ is finitely generated as a module over 
$\Bbb C[S_1(x),...,S_n(x)]$ (this is satisfied in interesting cases). 
In this case by a standard theorem in commutative algebra (see
Section  10.5), $\Bbb C[T_x^*X]$ is a free module 
over $\Bbb C[S_1(x),...,S_n(x)]$ of some rank $r$. 

\begin{proposition}\label{rdim}
The system (\ref{jep}) has an $r$-dimensional space of solutions near 
each point $x_0$ of $X$. 
\end{proposition}

\begin{proof} Let $P_1,...,P_r$ be the free homogeneous generators 
of $\Bbb C[T_x^*X]$ over $\Bbb C[S_1(x_0),...,S_n(x_0)]$, and 
$D_1,...,D_r\in D(X)$ be liftings of $P_1,...,P_r$. 
Then $\psi$ is a solution of (\ref{jep}) iff
the functions $D_1\psi$,...,$D_r\psi$ satisfy a first order 
linear holonomic system of differential equations near $x_0$
(i.e., represent a horizontal section of a flat connection):
$$
d(D_i\psi)=\sum_j \omega_{ij}D_j\psi,
$$
where $\omega=(\omega_{ij})$ is a matrix of 1-forms 
on $X$ satisfying the Maurer-Cartan equation. 
Therefore, the space of solutions is $r$-dimensional (as solution
is uniquely determined by the values of $D_i\psi(x_0)$). 
\end{proof}

Now we see the main difference between integrable and 
nonintegrable Hamiltonians $H$. Namely, 
we see from the proof of Proposition \ref{rdim} 
that solutions of the eigenvalue problem (\ref{jep})
can be found by solving {\bf ordinary} differential equations
(computing holonomy of a flat connection), while in the nonitegrable situation 
$H\psi=\lambda\psi$ is a {\bf partial} differential equation, 
which in general does not reduce to ODE. In the theory of PDE, we 
always regard reduction to ODE as an explicit solution. 
This justifies the above statement that the qunatum integrable systems, 
like the classical ones, can be solved ``explicitly''.  

{\bf Remark.} It can be argued that finding solutions of the eigenvalue 
problem (\ref{jep}) is the quantum analog of finding the action-angle 
variables of the quantum integrable system. 

\subsection{Constructing quantum integrable systems by quantum 
Hamiltonian reduction}

As in the classical case, an effective way of constructing quantum integrable 
systems is quantum Hamiltonian reduction. We will describe 
this procedure in the setting of formal deformations; 
the case of algebraic deformations and deformations coming from filtrations 
is similar.

Namely, suppose we are in the setting of Subsection \ref{concl}.  
That is, let $A_0$ be the function algebra on a symplectic ($C^\infty$ or 
algebraic) manifold $M$, and suppose that 
$\mu_0: \g\to A_0$ is a classical moment map. Let 
$H_{1,0}$,..., $H_{n,0}$ be a Poisson 
commuting family of $\g$-invariant elements of $A_0$, which reduces 
to an integrable system on $R(M,G,\O)$. Suppose that $A$  
is a deformation 
quantization of $A_0$, and $\mu: \g\to h^{-1}A$ is a quantum moment map 
quantizing $\mu_0$ (so $\mu=\hbar^{-1}\mu_0+O(1)$ on $\g$). 

Let $I_0$ be the ideal in $A_0$ of functions vanishing on $\O$, and $I$ a deformation of $I_0$ to an ideal in $U(\g_\hbar)$.  
Assume that $H_1,...,H_n$ is a commuting system of 
$\g$-invariant elements of $A$ quantizing $H_{1,0},...,H_{n,0}$.
Suppose that quantization commutes with reduction, i.e. that 
the algebra $R(A,\g,I)$ is a quantization of the 
symplectic manifold $R(M,G,\O)$. 
In this case, $H_1,...,H_n$ descend to commuting elements 
in $R(A,\g,I)$, which quantize $H_{1,0},...,H_{n,0}$, so they are a quantum 
integrable system. This system is called the quantum reduction 
of the original system $H_1,...,H_n$ in $A$. 

\subsection{The quantum Calogero-Moser system}

We now turn to the situation of Example \ref{quantumcm}.
Thus $\g={\frak {gl}}_n$ and $M=T^*\g=\lbrace{(X,P)|X,P\in \g\rbrace}$
(we identify $\g$ with $\g^*$ 
using the trace form). Set $H_i$, $i=1,...,n$, to be 
the homogeneous differential operators with constant coefficients 
on $\g$ with symbols $\Tr(P^i)$. According to the above, they 
descend to a quantum integrable system in the 
Cherednik algebra $B_k$. 

\begin{definition} This system (with $H=H_2$ being the Hamiltonian)
is called the quantized Calogero-Moser system. 
\end{definition}

To make the quantized Calogero-Moser system more explicit, 
let us recall that the algebra $B_k$ naturally acts on the space 
$\Bbb C[\h_{reg}]$. Thus $H_1,...,H_n$ define 
commuting differential operators on $\h$ 
with poles along the reflection hyperplanes
(it is easy to see that these operators are $W$-invariant).  
Denote these operators by $L_1',...,L_n'$.  
By analogy with the definition of the Harish-Chandra homomorphism, 
define $L_i=\delta\circ L_i'\circ \delta^{-1}$. 
It follows from the above that for $1\le m\le n$, one has
$$
L_m=\sum_{j=1}^n\partial_j^m+\text{lower order terms},
$$
where the lower order terms are polynomial in $k$, and vanish as $k=0$ 
(as for $k=0$ one has $L_m=HC(H_m)$). 

It is easy to see that $L_1=\sum_{j=1}^n \partial_j$. Now 
let us calculate more explicitly the Hamiltonian $H=L_2$. 

\begin{theorem}\label{qcm}
$$
L_2=\sum_{j=1}^n \partial_j^2-\sum_{i\ne j}\frac{k(k+1)}{(x_i-x_j)^2}.
$$
\end{theorem}

\begin{proof} The proof is analogous to the proof 
of Lemma \ref{lapl}. Namely, let us first calculate $L_2'$. 

Recall that the Casimir $D$ of $\g$ is given by the formula 
$$
D=\sum_{i=1}^n \partial_{x_i}^2+2\sum_{\alpha>0}\partial_{f_\alpha}\partial_{e_\alpha}. 
$$
Thus if $F(x)$ is a $\g$-equivariant function on 
the formal neighborhood of $\h_{reg}$ in $\g_{reg}$
with values in $W_k$ then we get 
$$
(DF)|_{\h}=\sum_{i=1}^n \partial_{x_i}^2(F|_\h)+
2\sum_{\alpha>0}(\partial_{f_\alpha}\partial_{e_\alpha}F)|_{\h}. 
$$
Now let $x\in \h$, and consider 
$(\partial_{f_\alpha}\partial_{e_\alpha}F)(x)$. We have 
$$
(\partial_{f_\alpha}\partial_{e_\alpha}F)(x)=
\partial_s\partial_t|_{s=t=0}F(x+tf_\alpha+se_\alpha).
$$
On the other hand, we have 
$$
{\rm Ad}(e^{s\alpha(x)^{-1}e_\alpha})(x+tf_\alpha+se_\alpha)=
x+tf_\alpha+ts\alpha(x)^{-1}h_\alpha+...,
$$
where $h_\alpha=[e_\alpha,f_\alpha]$. Hence, 
$$
F(x+tf_\alpha+se_\alpha)=
e^{-s\alpha(x)^{-1}e_\alpha}F(x+tf_\alpha+ts\alpha(x)^{-1}h_\alpha+...)=
$$
$$
e^{-s\alpha(x)^{-1}e_\alpha}
e^{t\alpha(x)^{-1}f_\alpha}F(x+ts\alpha(x)^{-1}h_\alpha+...).
$$
Thus 
$$
\partial_s\partial_t|_{s=t=0}F(x+tf_\alpha+se_\alpha)=
\alpha(x)^{-1}(\partial_{h_\alpha}F)(x)-\alpha(x)^{-2}e_\alpha f_\alpha 
F(x). 
$$
But $e_\alpha f_\alpha|_{W_k[0]}=k(k+1)$
(this is obtained by a direct computation using that $e_\alpha=
x_i\partial_j,f_\alpha=x_j\partial_i$).
This implies that 
$$
L_2'F(x)=\Delta_\h F(x)+2\sum_{\alpha>0}(\alpha(x)^{-1}\partial_{h_\alpha}F(x)
-k(k+1)\alpha(x)^{-2}F(x)).
$$
The rest of the proof is the same as in Lemma \ref{lapl}. 
\end{proof} 

Thus we see that the quantum Calogero-Moser system indeed describes 
a system of $n$ quantum particles on the line with interaction potential 
$k(k+1)/x^2$. 

{\bf Remark.} Note that unlike the classical case, in the
quantum case the coefficient
in front of the potential is an essential parameter, and
cannot be removed by rescaling.

\subsection{Notes} 1. Quantum integrable systems have been
studied for more than twenty years; let us mention, for instance,
the paper \cite{OP} which is relevant to the subject of these
lectures. The construction of such systems by quantum reduction 
has also been known for a long time: a famous instance of such 
a construction is the Harish-Chandra-Helgason theory of radial
parts of Laplace operators on symmetric spaces, see
e.g. \cite{He}. The quantum Calogero-Moser system (more precisely, its
trigonometric deformation), appeared in \cite{Su}. 

2. The construction of the quantum Calogero-Moser system 
by quantum reduction, which is parallel to the
Kazhdan-Kostant-Strenberg
construction,  is discussed, for instance, in \cite{EG}. 
The trigonometric version of the Calogero-Moser system
(the Sutherland system) with Hamiltonian 
$$
L_2=\sum_{j=1}^n \partial_j^2-\sum_{i\ne j}\frac{k(k+1)}{\sin^2(x_i-x_j)}
$$
is also integrable, and may be obtained by performing reduction 
from $T^*G$ (rather than $T^*{\mathfrak g}$), both classically and
quantum mechanically, mimicking the Kazhdan-Kostant-Sternberg 
construction. Another class of quantum integrable systems which can be
obtained by reduction is the Toda systems; this is done in the
paper by Kostant \cite{Kos}. 

\section{Calogero-Moser systems associated to finite Coxeter groups} 

It turns out that both the classical and the quantum Calogero-Moser system
can be generalized to the case of any finite Coxeter group $W$,
so that the case we have considered corresponds to the symmetric group $S_n$.  
This is done using Dunkl operators. 

\subsection{Dunkl operators}

Let $W$ be a finite Coxeter group, and $\h$ be its reflection representation,
$\dim\h=r$. Let $S$ be the set of reflections in $W$. 
For any reflection $s\in S$, let $\alpha_s\in \h^*$ be an eigenvector
of $s$ with eigenvalue $-1$. Then the reflection hyperplane of 
$s$ is given by the equation $\alpha_s=0$. Let $\alpha_s^\vee\in \h$ be 
the $-1$-eigenvector of $s$ such that $(\alpha_s,\alpha_s^\vee)=2$.

Let $c: S\to \Bbb C$ be a function invariant with respect to 
conjugation. Let $a\in \h$. 

\begin{definition} The Dunkl operator $D_a=D_a(c)$ on $\Bbb C(\h)$  
is defined by the formula 
$$
D_a=D_a(c):=\partial_a -\sum_{s\in S}\frac{c_s\alpha_s(a)}{\alpha_s}(1-s)
$$
\end{definition}

Clearly, $D_a\in \Bbb CW\ltimes D(\h_{reg})$. 

\begin{example}
Let $W=\Bbb Z_2$, $\h=\Bbb C$. Then there is only one Dunkl operator up 
to scaling, and it equals to 
$$
D=\partial_x -\frac{c}{x}(1-s),
$$
where the operator $s$ is given by the formula $(sf)(x)=f(-x)$. 
\end{example}

\begin{exercise} Show that $D_a$ maps the space of polynomials 
$\Bbb C[\h]$ to itself.
\end{exercise}

\begin{proposition}\label{maincom}
(i) For any $x\in \h^*$, one has
$$
[D_a,x]=(a,x)-\sum_{s\in S}c_s(a,\alpha_s)(x,\alpha_s^\vee)s.
$$

(ii) If $g\in W$ then $gD_ag^{-1}=D_{ga}$.
\end{proposition}

\begin{proof} (i) The proof follows immediately from the identity 
$x-sx=(x,\alpha_s^\vee)\alpha_s$.

(ii) The identity is obvious from the invariance of the function $c$.  
\end{proof}

The main result about Dunkl operators, on which all their applications 
are based, is the following theorem. 

\begin{theorem}\label{dunkl} (C. Dunkl, \cite{Du}) The Dunkl operators commute: 
$[D_a,D_b]=0$ for any $a,b\in \Bbb \h$. 
\end{theorem}

\begin{proof} Let $x\in \h^*$. We have 
$$
[[D_a,D_b],x]=[[D_a,x],D_b]-[[D_b,x],D_a].
$$
Now, using Proposition \ref{maincom}, we obtain:
$$
[[D_a,x],D_b]=-[\sum_s c_s (a,\alpha_s)(x,\alpha_s^\vee)s,D_b]=
$$
$$
-\sum_sc_s (a,\alpha_s)(x,\alpha_s^\vee)(b,\alpha_s)sD_{\alpha_s^\vee}.
$$
Since $a$ and $b$ occur symmetrically, we obtain that $[[D_a,D_b],x]=0$. 
This means that for any $f\in \Bbb C[\h]$, $[D_a,D_b]f=f[D_a,D_b]1=0$. 
\end{proof} 

\subsection{Olshanetsky-Perelomov operators}

Assume for simplicity that $\h$ is an irreducible $W$-module. 
Fix an invariant inner product $(,)$ on $\h$. 

\begin{definition}\cite{OP}
The Olshanetsky-Perelomov operator corresponding to a $W$-invariant 
function $c: S\to \Bbb C$ is the second order differential operator
$$
L:=\Delta_\h-\sum_{s\in S}\frac{c_s(c_s+1)(\alpha_s,\alpha_s)}{\alpha_s^2}.
$$
\end{definition} 

It is obvious that $L$ is a $W$-invariant operator.

\begin{example} Suppose that $W=S_n$. Then there is only one conjugacy 
class of reflections, so the function $c$ takes only one value. 
So $L$ is the quantum Calogero-Moser Hamiltonian, with $c=k$. 
\end{example}

\begin{exercise} Write $L$ explicitly for $W$ of type $B_n$. 
\end{exercise}

It turns out that the operator $L$ defines a quantum integrable system. 
This fact was discovered by Olshanetsky and Perelomov
(in the case when $W$ is a Weyl group). We are going 
to give a simple proof of this fact, due to Heckman, based on Dunkl operators.

To do so, note that any element $B\in (\Bbb C W\ltimes D(\h_{reg}))^W$ defines 
a linear operator on $\Bbb C(\h)$. Let  
$m(B)$ be the restriction of the operator $B$ to 
the subspace of $W$-invariant functions, $\Bbb C(\h)^W$. 
It is clear that for any $B$, $m(B)$ is a differential operator. Therefore 
the assignment $B\mapsto m(B)$ defines a homomorphism 
$m: (\Bbb C W\ltimes D(\h_{reg}))^W\to 
D(\h_{reg})^W$.

Define the operator
$$
\bar L:=\Delta_\h-\sum_{s\in S}\frac{c_s(\alpha_s,\alpha_s)}
{\alpha_s}\partial_{\alpha_s^\vee}
$$

\begin{proposition} (Heckman, \cite{Hec}) Let $\{y_1,\ldots y_r\}$ be  
an orthonormal basis of $\frak h$. Then we have 
 
$$m(\sum_{i=1}^rD_{y_i}^2)=\bar L.$$ 
\end{proposition} 

\begin{proof}  
Let us extend the map $m$ to $\Bbb C W\ltimes D(\h_{reg})$ 
by defining $m(B)$, $B\in \Bbb C W\ltimes D(\h_{reg})$, to 
be the (differential) operator $\Bbb C(\h)^W\to \Bbb C(\h)$ 
corresponding to $B$. Then 
we have $m(D_y^2)=m(D_y\partial_y)$. A simple computation 
shows that 
$$
D_y\partial_y=
\partial_y^2-\sum_{s\in S}\frac{c_s\alpha_s(y)}{\alpha_s}(1-s)\partial_y=
$$
$$
\partial_y^2-
\sum_{s\in S}\frac{c_s\alpha_s(y)}{\alpha_s}(\partial_y(1-s)+\alpha_s(y)
\partial_{\alpha_s^\vee}s).
$$
This means that 
$$
m(D_y^2)=\partial_y^2-\sum_{s\in S}c_s{\alpha_s(y)^2\over  
\alpha_s}\partial_{\alpha_s^\vee}.
$$ 
So we get 
$$
m(\sum_{i=1}^rD_{y_i}^2)=\sum_i\partial_{y_i}^2- 
\sum_{s\in S}c_s{\sum_{i=1}^r\alpha_s(y_{i})^2\over  
\alpha_s}\partial_{\alpha_s^\vee}=\bar L$$  
since $\sum_{i=1}^r\alpha_s(y_i)^2=  
(\alpha_s,\alpha_s)$. 
\end{proof}  

Recall that by Chevalley's theorem,
the algebra $(S\h)^W$ is free. Let 
$P_1=x^2,P_2,...,P_r$ be homogeneous generators of $(S\h)^W$. 

\begin{corollary}
The operator $\bar L$ defines a quantum integrable system. 
Namely, the operators $\bar L_i:=m(P_i(D_{y_1},...,D_{y_r}))$ are 
commuting quantum integrals of $\bar L=\bar L_1$. 
\end{corollary}

\begin{proof} The statement follows from Theorem \ref{dunkl} and
the fact that $m(b_1b_2)=m(b_1)m(b_2)$. 
\end{proof}

To derive from this the quantum integrability of the operator $L$, we 
will prove the following proposition. 

\begin{proposition} Let $\delta_c(x):=\prod_{s\in S}\alpha_s(x)^{c_s}$. 
Then we have 
$$
\delta_c^{-1}\circ \bar L\circ \delta_c=L.
$$
\end{proposition}

\begin{proof}
We have 
$$
\sum_i \partial_{y_i}(\log\delta_c)\partial_{y_i}=
\sum_{s\in S}\frac{c_s(\alpha_s,\alpha_s)}{2\alpha_s}\partial_{\alpha_s^\vee}.
$$
Therefore, we have 
$$
\delta_c \circ L\circ 
\delta_c^{-1}=\Delta_\h-
\sum_{s\in S}\frac{c_s(\alpha_s,\alpha_s)}{\alpha_s}\partial_{\alpha_s^\vee}
+U,
$$
where 
$$
U=\delta_c(\Delta_\h\delta_c^{-1})-
\sum_{s\in S}\frac{c_s(c_s+1)(\alpha_s,\alpha_s)}{\alpha_s^2}.
$$
Let us compute $U$. We have 
$$
\delta_c(\Delta_\h \delta_c^{-1})=
\sum_{s\in S}\frac{c_s(c_s+1)(\alpha_s,\alpha_s)}{\alpha_s^2}+
\sum_{s\ne u\in S}\frac{c_sc_u(\alpha_s,\alpha_u)}{\alpha_s\alpha_u}.
$$
We claim that the last sum $\Sigma$ is actually zero. Indeed, this sum is 
invariant under the Weyl group, so $\prod_{s\in S}\alpha_s\cdot \Sigma$
is a regular anti-invariant of degree $|S|-2$. But the smallest
degree of a nonzero anti-invariant is $|S|$, so 
$\Sigma=0$, $U=0$, and we are done. 
\end{proof}

\begin{corollary} (Heckman \cite{Hec})
The Olshanetsky-Perelomov operator $L$ defines a quantum integrable system, 
namely $\lbrace{L_i, i=1,...,r\rbrace}$, where 
$L_i=\delta_c^{-1}\circ \bar L_i\circ \delta_c$. 
\end{corollary}

\subsection{Classical Dunkl operators and Olshanetsky-Perelomov hamiltonians}

Let us define the classical analog of Dunkl and Olshanetsky-Perelomov 
operators. For this purpose we need to introduce the Planck constant $\hbar$. 
Namely, let us define renormalized Dunkl operators 
$$
D_a(c,\hbar):=\hbar D_a(c/\hbar).
$$
These operators can be regarded as elements of the Rees algebra 
$A={\rm Rees}(\Bbb CW\ltimes D(\h_{reg}))$, where 
the filtration is by order of differential operators (and $W$ sits in 
degree $0$). Reducing these operators modulo $\hbar$, we get 
classical Dunkl operators $D_a^0(c)\in A_0:=A/\hbar A=\Bbb 
CW\ltimes \O(T^*\h_{reg})$. They are given by the formula 
$$
D_a^0(c)=p_a-\sum_{s\in S}\frac{c_s\alpha_s(a)}{\alpha_s}(1-s),
$$
where $p_a$ is the classical momentum (the linear function on 
$\h^*$ corresponding to $a\in \h$). 

It follows from the commutativity of the quantum Dunkl operators $D_a$ that 
the classical Dunkl operators $D_a^0$ also commute: 
$$
[D_a^0,D_b^0]=0.
$$ 

We also have the following analog of Proposition \ref{maincom}: 

\begin{proposition}\label{maincom1}
(i) For any $x\in \h^*$, one has
$$
[D_a^0,x]=-\sum_{s\in S}c_s(a,\alpha_s)(x,\alpha_s^\vee)s.
$$

(ii) If $g\in W$ then $gD_a^0g^{-1}=D_{ga}^0$.
\end{proposition}

Now let us construct the classical Olshanetsky-Perelomov hamiltonians. 
As in the quantum case, we have a homomorphism
$m: (\Bbb CW\ltimes \O(T^*\h_{reg}))^W\to \O(T^*\h_{reg})^W$, which 
is given by the formula $\sum f_g\cdot g\to \sum f_g$, $g\in W$, 
$f\in \O(T^*\h_{reg})$. We define 
the hamiltonian
$$
\bar L^0:=m(\sum_{i=1}^r(D_{y_i}^0)^2).
$$
It is easy to see by taking the limit from the quantum situation 
that 
$$
\bar L^0=
p^2-\sum_{s\in S}\frac{c_s(\alpha_s,\alpha_s)}{\alpha_s}p_{\alpha_s^\vee}.
$$

Consider the (outer) automorphism $\theta_c$ of the algebra $\Bbb CW\ltimes 
\O(T^*\h_{reg})$ defined by the formulas 
$$
\theta_c(x)=x,\ \theta_c(s)=s,\ \theta_c(p_a)=p_a+\partial_a \log\delta_c,
$$
for $x\in \h^*$, $a\in \h$, $s\in W$. 
It is easy to see that if $b_0\in A_0$
and $b\in A$ is a deformation of $b_0$ then 
$\theta_c(b_0)=\lim_{\hbar\to 0}\delta_{c/\hbar}^{-1}b\delta_{c/\hbar}$. 
Therefore, defining 
$$
L^0:=\theta_c(\bar L^0), 
$$
we find:
$$
L^0=p^2-\sum_{s\in S}\frac{c_s^2(\alpha_s,\alpha_s)}{\alpha_s^2}.
$$
This function is called the classical Olshanetsky-Perelomov hamiltonian
for $W$. 

Thus, we obtain the following result. 

\begin{theorem}
The Olshanetsky-Perelomov hamiltonian $L^0$ defines an
integrable system with integrals 
$$
L^0_i:=m(\theta_c(P_i(D_{y_1}^0,...,D_{y_r}^0)))
$$
(so that $L^0_1=L^0$).
\end{theorem}

\subsection{Notes} 1. Dunkl operators can be generalized to
complex reflection groups, see \cite{DO}. Using them in the same
manner as above, one can construct quantum integrable systems
whose Hamiltonians have order higher than 2. 

2. If $W$ is a Weyl group, then Dunkl operators for $W$ can be
extended to the trigonometric setting, see \cite{Op}. 
The trigonometric Dunkl operators can be used to
construct the first integrals of the trigonometric 
quantum Calogero-Moser system, 
with Hamiltonian
$$
L_{trig}:=\Delta_\h-\sum_{s\in S}\frac{c_s(c_s+1)(\alpha_s,\alpha_s)}{\sin^2(\alpha_s)}.
$$
Moreover, these statements
generalize to the case of the quantum elliptic Calogero-Moser
system, with Hamiltonian
$$
L_{ell}:=\Delta_\h-\sum_{s\in S}c_s(c_s+1)(\alpha_s,\alpha_s)\wp(\alpha_s,\tau),
$$
where $\wp$ is the Weierstrass elliptic function. This is done in
the paper \cite{Ch1}.

\section{The rational Cherednik algebra}

\subsection{Definition of the rational Cherednik algebra
and the Poincar\'e-Birkhoff-Witt theorem}

In the previous section we made an essential use of the 
commutation relations between operators
$x\in \h^*$, $g\in W$, and $D_a$, $a\in \h$. 
This makes it natural to consider the algebra generated 
by these operators. 

\begin{definition} The rational Cherednik algebra $H_c=H_c(W,\h)$ 
associated to $(W,\h)$ is the algebra generated inside 
$A={\rm Rees}(\Bbb CW\ltimes D(\h_{reg}))$ by 
the elements $x\in \h^*$, $g\in W$, and $D_a(c,\hbar)$, 
$a\in \h$. If $t\in \Bbb C$, then 
the algebra $H_{t,c}$ is the specialization 
of $H_c$ at $\hbar=t$. 
\end{definition}

\begin{proposition}\label{rela}
The algebra $H_c$ is the quotient of the algebra \linebreak
$\Bbb CW\ltimes \bold T(\h\oplus \h^*)[\hbar]$ 
(where $\bold T$ denotes the tensor algebra) 
by the ideal generated by the relations
$$
[x,x']=0,\ [y,y']=0,\ [y,x]=\hbar(y,x)-\sum_{s\in S}
c_s(y,\alpha_s)(x,\alpha_s^\vee)s,
$$
$x,x'\in \h^*$, $y,y'\in \h$.
\end{proposition}

\begin{proof}
Let us denote the algebra defined in the proposition by $H_c'$. 
Then according to the results of the previous section, we have a 
surjective homomorphism 
$\phi: H_c'\to H_c$ defined by the formula 
$\phi(x)=x$, $\phi(g)=g$, $\phi(y)=D_y(c,\hbar)$. 

Let us show that this homomorphism is injective. For this purpose 
assume that $y_i$ is a basis of $\h$, and $x_i$ is the dual basis
of $\h^*$. Then it is clear from the relations of $H_c'$ that 
$H_c'$ is spanned over $\Bbb C[\hbar]$ by the elements 
\begin{equation}\label{basi}
g\prod_{i=1}^r y_i^{m_i}
\prod_{i=1}^r x_i^{n_i}.
\end{equation} 

Thus it remains to show that the images of the elements 
(\ref{basi}) under the map $\phi$, 
i.e. the elements 
$$
g\prod_{i=1}^r D_{y_i}(c,\hbar)^{m_i}
\prod_{i=1}^r x_i^{n_i}.
$$
are linearly independent. 
But this follows from the obvious fact that the symbols of these elements 
in $\Bbb CW\ltimes \Bbb C[\h^*\times \h_{reg}][\hbar]$
are linearly independent. The proposition is proved.  
\end{proof}

It is more convenient to work with algebras defined by generators and 
relations than with subalgebras of a given algebra generated by a 
given set of elements. Therefore, from now on we will use 
the statement of Proposition \ref{rela} 
as a definition of the rational Cherednik algebra $H_c$. 
According to Proposition \ref{rela}, this algebra comes with 
a natural embedding $\Theta_c: H_c\to {\rm Rees}(\Bbb CW\ltimes D(\h_{reg}))$,
defined by the formula $x\to x$, $g\to g$, $y\to D_y(c,\hbar)$. 
This embedding is called the Dunkl operator embedding. 

Let us put a filtration on $H_c$ by declaring $x\in \h^*$ and 
$y\in \h$ to have degree $1$, and $g\in G$ to have degree $0$. 
Let ${\rm gr}(H_c)$ denote the associated graded algebra 
of $H_c$ under this filtration. We have a 
natural surjective homomorphism $\xi: \Bbb CW\ltimes 
\Bbb C[\h\oplus \h^*][\hbar]\to {\rm gr}(H_c)$.  

\begin{proposition}
\label{pbw} (the PBW theorem for rational Cherednik algebras) 
The map $\xi$ is an isomorphism. 
\end{proposition}

\begin{proof} The statement is equivalent to 
the claim that the elements (\ref{basi}) are a basis 
of $H_c$, which follows from the proof of Proposition \ref{rela}. 
\end{proof}

\begin{remark} 1. It follows from Proposition \ref{pbw}
that the algebra $H_c$ is a free module over $\Bbb C[\hbar]$, 
whose basis is  the collection of elements (\ref{rela})
(In particular, $H_c$ is an algebraic deformation of $H_{0,c}$).
This basis is called a PBW basis of $H_c$ (and of $H_{t,c}$). 

2. It follows from the definition of $H_{t,c}$ 
that $H_{t,c}$ is a quotient of 
\linebreak $\Bbb CW\ltimes \bold T(\h\oplus \h^*)$ 
by the relations of Proposition \ref{rela}, with $\hbar$ replaced by $t$. 
In particular, 
for any $\lambda\in \Bbb C^*$, the algebra $H_{t,c}$ is naturally 
isomorphic to $H_{\lambda t,\lambda c}$. 

3. The Dunkl operator embedding $\Theta_c$ specializes to 
embeddings $\Theta_{1,c}: H_{1,c}\to \Bbb CW\ltimes D(\h_{reg})$
given by $x\to x$, $g\to g$, $y\to D_y$, and 
$\Theta_{0,c}: H_{0,c}\to \Bbb CW\ltimes \Bbb C[\h^*\times \h_{reg}]$, 
given by $x\to x$, $g\to g$, $y\to D_y^0$. 

4. Since Dunkl operators map polynomials to polynomials, 
the map $\Theta_{1,c}$ defines a representation 
of $H_{1,c}$ on $\Bbb C[\h]$. This representation
is called the polynomial representation of $H_{1,c}$. 
\end{remark}

\begin{example} 1. Let $W=\Bbb Z_2$, 
$\h=\Bbb C$. In this case $c$ reduces to one parameter $k$, and the algebra 
$H_{t,k}$ is generated by elements $x,y,s$ with defining relations
$$
s^2=1,\ sx=-xs,\ sy=-ys,\ [y,x]=t-2ks. 
$$

2. Let $W=S_n$, $\h=\Bbb C^n$. In this case there is also only one 
complex parameter $c=k$, and the algebra 
$H_{t,k}$ is the quotient of 
\linebreak $S_n\ltimes \Bbb C<x_1,...,x_n,y_1,...,y_n>$
by the relations 
\begin{equation}\label{gln}
[x_i,x_j]=[y_i,y_j]=0,\ [y_i.x_j]=ks_{ij},\ [y_i,x_i]=t-k\sum_{j\ne i}s_{ij}.
\end{equation}
Here $\Bbb C\langle E\rangle$ denotes the free algebra on a set $E$, 
and $s_{ij}$ is the transposition of $i$ and $j$. 
\end{example}

\subsection{The spherical subalgebra}

Let $\e\in \Bbb CW$ be the symnmetrizer, 
$\e=|W|^{-1}\sum_{g\in W}g$. We have $\e^2=\e$. 

\begin{definition}
The spherical subalgebra of $H_c$ is 
the subalgebra $B_c:=\e H_{c}\e$.
The spherical subalgebra of $H_{t,c}$ 
is $B_{t,c}:=B_c/(\hbar=t)=\e H_{t,c}\e$.  
\end{definition} 

Note that $\e \Bbb CW\ltimes D(\h_{reg})\e=D(\h_{reg})^W$. 
Therefore, the map $\Theta_{t,c}$ 
restricts to an embedding $B_{t,c}\to D(\h_{reg})^W$ 
for $t\ne 0$, and $B_{0,c}\to \Bbb C[\h^*\times \h_{reg}]^W$
for $t=0$. 

\begin{proposition}\label{zerdiv} (i) The spherical subalgebra 
$B_{0,c}$ is commutative and does not have zero divisors. 

(ii) $B_c$ is an algebraic deformation of $B_{0,c}$.
\end{proposition}

\begin{proof}
(i) Follows immediately from the fact that $B_{0,c}\subset  
\Bbb C[\h^*\times \h_{reg}]^W$.

(ii) Follows since $H_c$ is an algebraic deformation of $H_{0,c}$.
\end{proof}

Proposition \ref{zerdiv} implies that the spectrum $M_c$ 
of $B_{0,c}$ is an irreducible affine algebraic variety. 
Moreover, $M_c$ has a natural Poisson structure, 
obtained from the deformation $B_c$ of $B_{0,c}$. 
In fact, it is clear that this Poisson structure is simply the restriction 
of the Poisson structure of $\Bbb C[\h^*\times \h_{reg}]^W$ 
to the subalgebra $B_{0,c}$. 

\begin{definition}
The Poisson variety $M_c$ is called the Calogero-Moser space 
of $W,\h$. 
\end{definition}

{\bf Remark.} We will later justify this terminology by showing 
that in the case $W=S_n$, $\h=\Bbb C^n$ 
the variety $M_c$ is isomorphic, as a Poisson variety,
to the Calogero-Moser space of Kazhdan, Kostant, and Sternberg. 

Thus, we may say that the algebra $B_c$ is an algebraic quantization 
of the Calogero-Moser space $M_c$. 

\subsection{The localization lemma and the basic properties of $M_c$}

Let $H_{t,c}^{\rm loc}=H_{t,c}[\delta^{-1}]$ be the localization 
of $H_{t,c}$ as a module over $\Bbb C[\h]$ with respect to
the discriminant $\delta$. 
Define also $B_{t,c}^{\rm loc}=\e H^{\rm loc}_{t,c}\e$. 

\begin{proposition}\label{loc} (i) For $t\ne 0$ 
the map $\Theta_{t,c}$ induces an isomorphism 
of algebras from $H_{t,c}^{\rm loc}\to 
\Bbb CW\ltimes D(\h_{reg})$, which restricts 
to an isomorphism $B_{t,c}^{loc}\to D(\h_{reg})^W$. 

(ii) 
The map $\Theta_{0,c}$ induces an isomorphism 
of algebras from $H_{0,c}^{\rm loc}\to 
\Bbb CW\ltimes \Bbb C[\h^*\times \h_{reg}]$,
which restricts to an isomorphism 
$B_{0,c}^{loc}\to \Bbb C[\h^*\times \h_{reg}]^W$. 
\end{proposition}

\begin{proof}
This follows immediately from the fact that
the Dunkl operators have poles only on the reflection hyperplanes.  
\end{proof}

Thus we see that the dependence on $c$ disappears upon localization.

Since ${\rm gr}(B_{0,c})=B_{0,0}=\Bbb C[\h^*\oplus \h]^W$,
we get the following geometric corollary. 

\begin{corollary} (i) The family of Poisson varieties 
$M_c$ is a flat deformation of the Poisson variety $M_0:=(\h^*\oplus \h)/W$. 
In particular, $M_c$ is smooth outside of a subset of codimension $2$. 

(ii) We have a natural map $\beta_c: M_c\to \h/W$, such that 
$\beta_c^{-1}(\h_{reg}/W)$ is isomorphic to $(\h^*\times \h_{reg})/W$.
\end{corollary}

\begin{exercise} Let $W=\Bbb Z_2$, $\h=\Bbb C$. 
Show that $M_c$ is isomorphic to the quadric 
$pq-r^2=c^2$ in the 3-dimensional space with coordinates $p,q,r$. 
In particular, $M_c$ is smooth for $c\ne 0$.   
\end{exercise}

\begin{exercise}
Show that if $t\ne 0$ then the center of $H_{t,c}$ is trivial. 

Hint. Use Proposition \ref{loc}.
\end{exercise} 

\subsection{The $SL_2$-action on $H_{t,c}$}

It is clear from the definition of $H_{t,c}$ that 
$\h$ and $\h^*$ play completely symmetric roles, 
i.e. there is a symnmetry that exchanges them. 
In fact, a stronger statement holds.  

\begin{proposition}\label{sl2}
The group $SL_2(\Bbb C)$ acts by automorphisms 
of $H_{t,c}$, via the formulas $g\to g$, 
$x_i\to px_i+qy_i$, $y_i\to rx_i+sy_i$, $g\to g$ ($g\in W$), where $ps-qr=1$,  
$x_i$ is an orthonormal basis of $\h^*$, and $y_i$ the dual basis 
of $\h$.   
\end{proposition}

\begin{proof} The invariant bilinear form on $\h$ defines an identification 
$\h\to \h^*$, under which $x_i\to y_i$, and $\alpha$ goes to a multiple of $\alpha^\vee$. This implies   
that $(\alpha^\vee,x_i)(\alpha,y_j)$ is symmetric in $i,j$. Thus $[y_i,x_j]=[y_j,x_i]$, and hence 
for any $p,q\in \Bbb  C$, the elements $px_i+qy_i$ commute with each other. 
Also, 
$$
[rx_i+sy_i,px_j+qy_j]=(ps-qr)[y_i,x_j]=[y_i,x_j].
$$
This implies the statement. 
\end{proof}

\subsection{Notes} Rational Cherednik algebras are degenerations 
of double affine Hecke algebras introduced by Cherednik in early
nineties. For introduction to double affine Hecke algebras, see
the book \cite{Ch}. The first systematic study 
of rational Cherednik algebras was undertaken in \cite{EG}. 
In particular, the results of this lecture can be found in
\cite{EG}. 

\section{Symplectic reflection algebras}

\subsection{The definition of symplectic reflection algebras}

Rational Cherednik algebras for finite Coxeter groups are 
a special case of a wider class of algebras called 
symplectic reflection algebras. To define them, let 
$V$ be a finite dimensional symplectic vector space 
over $\Bbb C$ with a symplectic form $\omega$, and $G$
be a finite group acting symplectically (linearly)
on $V$. For simplicity let us assume that $(\wedge^2V)^G=\Bbb C\omega$
and that $G$ acts faithfully on $V$ (these assumptions are not essential).  

\begin{definition}
A symplectic reflection in $G$ is an element $g$ such that the rank of the operator
$g-1$ on $V$ is $2$. 
\end{definition}

If $s$ is a symplectic reflection, then let $\omega_s(x,y)$ be 
the form $\omega$ applied to the projections of $x,y$ 
to the image of $1-s$ along the kernel of $1-s$; thus $\omega_s$ is a skewsymmetric form of rank $2$ on $V$.   

Let $S\subset G$ be the set of symplectic reflections, and 
$c: S\to \Bbb C$ be a function which is invariant under the action of $G$. 
Let $t\in \Bbb C$. 

\begin{definition} The symplectic reflection algebra $H_{t,c}=H_{t,c}[V,G]$
is the quotient of the algebra $\Bbb C[G]\ltimes \bold T(V)$ 
by the ideal generated by the relation 
\begin{equation}\label{mainrel}
[x,y]=t\omega(x,y)-2\sum_{s\in S}c_s\omega_s(x,y)s.
\end{equation}
\end{definition}

Note that if $W$ is a finite Coxeter group with reflection representation $\h$, 
then we can set $V=\h\oplus \h^*$, 
$\omega(x,x')=\omega(y,y')=0$, $\omega(y,x)=(y,x)$, 
for $x,x'\in \h^*$ and $y,y'\in \h$. In this case

1) symplectic reflections are the usual reflections in $W$;

2) $\omega_s(x,x')=\omega_s(y,y')=0$, 
$\omega_s(y,x)=\frac{1}{2}(y,\alpha_s)(\alpha_s^\vee,x)$. 

Thus, $H_{t,c}$ is the rational Cherednik algebra defined in the previous subsection. 

Note also that for any $V,G$, $H_{0,0}[V,G]=G\ltimes SV$, and $H_{1,0}[V,G]=G\ltimes {\rm Weyl}(V)$, where ${\rm Weyl}(V)$ is the 
Weyl algebra of $V$, i.e. the quotient of the tensor algebra $\Bbb T(V)$ by the relation 
$xy-yx=\omega(x,y)$, $x,y\in V$. 

\subsection{The PBW theorem for symplectic reflection algebras}

To ensure that the symplectic reflection algebras $H_{t,c}$ have good properties, 
we need to prove a PBW theorem for them, which is an analog of Proposition \ref{pbw}. 
This is done in the following theorem, which also explains the special role 
played by symplectic reflections. 

\begin{theorem}\label{pbw1}
Let $\kappa: \wedge^2V\to \Bbb C[G]$ be a linear $G$-equivariant function. 
Define the algebra $H_\kappa$ to be the quotient of the algebra 
$\Bbb C[G]\ltimes {\bold T}(V)$ by the relation 
$[x,y]=\kappa(x,y)$, $x,y\in V$. Put an increasing filtration on $H_\kappa$ 
by setting $\deg(V)=1$, $\deg(G)=0$, and define $\xi: \Bbb C G\ltimes SV\to {\rm gr}H_\kappa$ 
to be the natural surjective homomorphism. Then $\xi$ is an isomorphism if and only if 
$\kappa$ has the form 
$$
\kappa(x,y)=
t\omega(x,y)-2\sum_{s\in S}c_s\omega_s(x,y)s,
$$
for some $t\in \Bbb C$ and $G$-invariant function $c: S\to \Bbb C$. 
\end{theorem}

Unfortunately, for a general symplectic reflection algebra 
we don't have a Dunkl operator representation, so
the proof of the more difficult ``if'' part of this Theorem 
is not as easy as the proof of Proposition \ref{pbw}. 
Instead of explicit computations with Dunkl operators, it makes use 
of the deformation theory of Koszul algebras, which we will now discuss. 

\subsection{Koszul algebras}

Let $R$ be a finite dimensional semisimple algebra (over $\Bbb C$). 
Let $A$ be a $\Bbb Z_+$-graded algebra, such that $A[0]=R$. 

\begin{definition} (i) The algebra $A$ is said to be quadratic if 
it is generated over $R$ by $A[1]$, and has 
defining relations in degree 2.

(ii) $A$ is Koszul if all elements of $Ext^i(R,R)$ 
(where $R$ is the augmentation module over $A$)
have grade degree precisely $i$. 
\end{definition}

{\bf Remarks.} 1. Thus, in a quadratic algebra, 
$A[2]=A[1]\otimes_R A[1]/E$, where $E$ is the 
subspace ($R$-subbimodule) of relations. 

2. It is easy to show that 
a Koszul algebra is quadratic, since 
the condition to be quadratic is just the Koszulity condition 
for $i=1,2$.

Now let $A_0$ be a quadratic algebra, $A_0[0]=R$. 
Let $E_0$ be the space of relations for $A_0$.
Let $E\subset A_0[1]\otimes_R A_0[1][[\hbar]]$ 
be a topologically free (over $\Bbb C[[\hbar]]$) $R$-subbimodule which 
reduces to $E_0$ modulo $\hbar$ (``deformation of the relations''). 
Let $A$ be the ($\hbar$-adically complete) 
algebra generated over $R[[\hbar]]$ by $A[1]=A_0[1][[\hbar]]$ 
with the space of defining relations $E$. 
Thus $A$ is a $\Bbb Z_+$-graded algebra. 

The following very important theorem 
is due to Beilinson, Ginzburg, and Soergel, \cite{BGS} 
(less general versions appeared earlier 
in the works of Drinfeld \cite{Dr1}, Polishchuk-Positselski \cite{PP}, 
Braverman-Gaitsgory \cite{BG}). We will not give its proof.  

\begin{theorem}\label{KDP} (Koszul deformation principle)
If $A_0$ is Koszul then $A$ is a topologically free $\Bbb C[[\hbar]]$ module
if and only if so is $A[3]$.  
\end{theorem} 

{\bf Remark.} Note that $A[i]$ for $i<3$ are obviously topologically free. 

We will now apply this theorem to the proof of Theorem \ref{pbw1}. 

\subsection{Proof of Theorem \ref{pbw1}}

Let $\kappa: \wedge^2V\to \Bbb C[G]$ be a linear $G$-equivariant map. 
We write $\kappa(x,y)=\sum_{g\in G}\kappa_g(x,y)g$, where $\kappa_g(x,y)\in \wedge^2V^*$.
To apply Theorem \ref{KDP}, let us homogenize our algebras. 
Namely, let $A_0=(\Bbb CG\ltimes SV)\otimes \Bbb C[u]$ (where $u$
has degree $1$). Also 
let $\hbar$ be a formal parameter, and 
consider the deformation $A=H_{\hbar u^2\kappa}$ of 
$A_0$. That is, $A$ is the quotient of 
$\Bbb G\ltimes {\bold T}(V)[u][[\hbar]]$ by the relations $[x,y]=\hbar u^2\kappa(x,y)$. 
This is a deformation of the type considered in Theorem \ref{KDP}, and 
it is easy to see that its flatness 
in $\hbar$ is equivalent to Theorem \ref{pbw1}. 
Also, the algebra $A_0$ is Koszul, because the polynomial algebra 
$SV$ is a Koszul algebra. 
Thus by Theorem \ref{KDP}, it suffices to show that 
$A$ is flat in degree 3. 

The flatness condition in degree 3 is ``the Jacobi identity''
$$
[\kappa(x,y),z]+[\kappa(y,z),x]+[\kappa(z,x),y]=0,
$$
which must be satisfied in $\Bbb CG\ltimes V$. In components, 
this equation transforms into the system of equations
$$
\kappa_g(x,y)(z-z^g)+\kappa_g(y,z)(x-x^g)+\kappa_g(z,x)(y-y^g)=0
$$
for every $g\in G$ (here $z^g$ denotes the result of the action of 
$g$ on $z$). 

This equation, in particular, implies that if $x,y,g$ are such
that $\kappa_g(x,y)\ne 0$ then for any $z\in V$ 
$z-z^g$ is a linear combination of $x-x^g$ and $y-y^g$. 
Thus $\kappa_g(x,y)$ is identically zero unless the 
rank of $(1-g)|_V$ is at most 2, i.e. 
$g=1$ or $g$ is a symplectic reflection. 

If $g=1$ then $\kappa_g(x,y)$ has to be $G$-invariant, 
so it must be of the form $t\omega(x,y)$, where $t\in \Bbb C$. 

If $g$ is a symplectic reflection, then
$\kappa_g(x,y)$ must be zero for any $x$ such that $x-x^g=0$. 
Indeed, if for such an $x$ there had existed $y$ with 
$\kappa_g(x,y)\ne 0$ then $z-z^g$ for any $z$ would be a multiple of $y-y^g$, which is impossible since 
$Im(1-g)|_V$ is 2-dimensional. This implies that $\kappa_g(x,y)=
2c_g\omega_g(x,y)$, and $c_g$ 
must be invariant. 

Thus we have shown that if $A$ is flat (in degree 3) then $\kappa$ must have the form 
given in Theorem \ref{pbw1}. Conversely, it is easy to see that 
if $\kappa$ does have such form, then the Jacobi identity holds. 
So Theorem \ref{pbw1} is proved.

\subsection{The spherical subalgebra of the symplectic reflection algebra}

The properties of symplectic reflection algebras are similar to properties 
of rational Cherednik algebras we have studied before. 
The main difference is that we no longer have the Dunkl representation and localization results, 
so some proofs are based on different ideas and are more complicated. 

The spherical subalgebra of the symplectic reflection algebra is defined 
in the same way as in the Coxeter group case. 
Namely, let $\e=\frac{1}{|G|}\sum_{g\in G}g$, and $B_{t,c}=\e H_{t,c}\e$. 

\begin{proposition} 
$B_{t,c}$ is commutative if and only if 
$t=0$. 
\end{proposition}

\begin{proof}
Let $A$ be a $\Bbb Z_+$-filtered algebra. If $A$ is not commutative, 
then we can define a Poisson bracket on ${\rm gr}(A)$ in the following way. 
Let $m$ be the minimum of $\deg(a)+\deg(b)-\deg([a,b])$ (over $a,b\in A$ such that $[a,b]\ne 0$). 
Then for homogeneous elements $a_0,b_0\in A_0$ of degrees $p,q$, we can define 
$\lbrace{a_0,b_0\rbrace}$ to be the image in $A_0[p+q-m]$ of $[a,b]$, where $a,b$ are any lifts 
of $a_0,b_0$ to $A$. It is easy to check that $\lbrace{,\rbrace}$ is a Poisson bracket 
on $A_0$ of degree $-m$. 

Let us now apply this construction to the filtered algebra $A=B_{t,c}$.
We have ${\rm gr}(A)=A_0=(SV)^G$. 

\begin{lemma}\label{poisbr} $A_0$ has a unique, up to scaling, Poisson bracket 
of degree $-2$, and no nonzero Poisson brackets of degrees $<-2$. 
\end{lemma}

\begin{proof}
A Poisson bracket on $(SV)^G$ is the same thing as a Poisson bracket on the variety $V/G$. 
On the smooth part $(V/G)_s$ of $V/G$, it is simply a bivector field, and we can lift it 
to a bivector field on the preimage $V_s$ of $(V/G)_s$ in $V$, which 
is the set of points in $V$ with trivial stabilizers. 
But the codimension on $V\setminus V_s$ in $V$ is 2 (as $V\setminus V_s$ is a union 
of symplectic subspaces), so the bivector on $V_s$ extends to a regular bivector on $V$. 
So if this bivector is homogeneous, it must have degree $\ge -2$, and if it has
degree $-2$ then it must be with constant coefficients, so being $G$-invariant, it 
is a multiple of $\omega^{-1}$. The lemma is proved.   
\end{proof}

Now, for each $t,c$ we have a natural Poisson bracket on 
$A_0$ of degree $-2$, which depends linearly on $t,c$. 
So by the lemma, this bracket has to be of the form $f(t,c)\Pi$, 
where $\Pi$ is the unique up to scaling 
Poisson bracket of degree $-2$, and $f$ a homogeneous linear function. 
Thus the algebra $A=B_{t,c}$ is not commutative unless $f(t,c)=0$. 
On the other hand, if $f(t,c)=0$, and $B_{t,c}$ is not commutative, then, as we've shown, $A_0$ 
has a nonzero Poisson bracket of degree $<-2$. But By Lemma \ref{poisbr}, there is no such brackets. 
So $B_{t,c}$ must be commutative if $f(t,c)=0$.

It remains to show that $f(t,c)$ is in fact a nonzero multiple of $t$. 
First note that $f(1,0)\ne 0$, since $B_{1,0}$ is definitely noncommutative. 
Next, let us take a point $(t,c)$ such that $B_{t,c}$ is commutative. 
Look at the $H_{t,c}$-module $H_{t,c}\e$, which has a commuting action of $B_{t,c}$ from the right. 
Its associated graded is $SV$ as an $(\Bbb CG\ltimes SV,(SV)^G)$-bimodule, which implies that 
the generic fiber of $H_{t,c}\e$ as a $B_{t,c}$-module is the regular representation of $G$. 
So we have a family of finite dimensional representations of $H_{t,c}$ 
on the fibers of $H_{t,c}\e$, all realized in the regular representation 
of $G$. Computing the trace of the main commutation relation 
(\ref{mainrel}) of $H_{t,c}$ in this representation, we obtain
that $t=0$ (since $Tr(s)=0$ for any reflection $s$). The Proposition is proved. 
\end{proof}

Note that $B_{0,c}$ has no zero divisors, since its associated graded algebra $(SV)^G$ does not. 
Thus, like in the Cherednik algebra case, we can define 
a Poisson variety $M_c$, the spectrum of $B_{0,c}$, called the Calogero-Moser space of $G,V$. 
Moreover, the algebra $B_c:=B_{\hbar,c}$ over $\Bbb C[\hbar]$ is an algebraic quantization of $M_c$. 

\subsection{Notes} 1. Symplectic reflection algebras are 
a special case of generalized Hecke algebras defined by Drinfeld
in Section 4 of \cite{Dr2}. They were systematically studied in
\cite{EG}, which is a basic reference for the results of this lecture. 
A generalization of symplectic reflection algebras, in which the
action of $G$ on $V$ may be non-faithful and projective, is 
considered by Chmutova \cite{Chm}.  

2. The notion of a Koszul algebra is due to Priddy \cite{Pr}. 
A good reference for the theory of quadratic and Koszul algebras
is the book \cite{PP}. 

\section{Deformation-theoretic interpretation of symplectic reflection algebras}

\subsection{Hochschild cohomology of semidirect products}

In the previous lectures we saw that the symplectic reflection algebra
$H_{1,c}[V,G]$, for formal $c$, is a flat 
formal deformation of $G\ltimes {\rm Weyl}(V)$. 
It turns out that this is in fact a universal deformation. 

Namely, let $A_0=G\ltimes {\rm Weyl}(V)$, 
where $V$ is a symplectic vector space (over $\Bbb C$), 
and $G$ a finite group acting symplectically on $V$. 
Let us calculate the Hochschild 
cohomology of this algebra.

\begin{theorem} \label{afls} (Alev, Farinati, Lambre, Solotar, \cite{AFLS})
The cohomology space $H^i(G\ltimes {\rm Weyl}(V))$ is 
naturally isomorphic to the space of 
conjugation invariant functions on the set $S_i$ 
of elements $g\in G$ such that ${\rm rank}(1-g)|_V=i$. 
\end{theorem}

Since ${\rm Im}(1-g)$ is a symplectic vector space, we get the following corollary. 

\begin{corollary}
The odd cohomology of $G\ltimes {\rm Weyl}(V)$ vanishes, 
and \linebreak $H^2(G\ltimes {\rm Weyl}(V))$ is the space $\Bbb C[S]^G$ of conjugation invariant 
functions on the set of symplectic reflections. 
In particular, there exists a universal deformation $A$ of $A_0=G\ltimes {\rm Weyl}(V)$
parametrized by $\Bbb C[S]^G$.
\end{corollary}

\begin{proof} (of Theorem \ref{afls})

\begin{lemma}
Let $B$ be a $\Bbb C$-algebra together with an action of 
a finite group $G$. Then 
$$
H^*(G\ltimes B,G\ltimes B)=(\oplus_{g\in G}H^*(B,Bg))^G,
$$ 
where $Bg$ is the bimodule isomorphic to $B$ as a space where 
the left action of $B$ is the usual one and the right action is the usual action twisted by $g$. 
\end{lemma}

\begin{proof} The algebra $G\ltimes B$ is a projective $B$-module. 
Therefore, using the Shapiro lemma, we get
$$
H^*(G\ltimes B,G\ltimes B)=\Ext^*_{(G\times G)\ltimes (B\otimes B^{op})}(G\ltimes B,G\ltimes B)=
$$
$$
\Ext^*_{G_{\rm diagonal}\ltimes (B\otimes B^{op})}(B,G\ltimes B)=
\Ext^*_{B\otimes B^{op}}(B,G\ltimes B)^G=
$$
$$
(\oplus_{g\in G}\Ext^*_{B\otimes B^{op}}(B,Bg))^G=
(\oplus_{g\in G}H^*(B,Bg))^G,
$$
as desired. 
\end{proof}

Now apply the lemma to $B={\rm Weyl}(V)$. 
For this we need to calculate $H^*(B,Bg)$, where 
$g$ is any element of $G$. We may assume that 
$g$ is diagonal in some symplectic basis: 
$g=\diag(\lambda_1,\lambda_1^{-1},...\lambda_n,\lambda_n^{-1})$. 
Then by the K\"unneth formula we find that 
$H^*(B,Bg)=\otimes_{i=1}^n H^*(\Bbb A_1,\Bbb A_1g_i)$, where 
$\Bbb A_1$ is the Weyl algebra of the 2-dimensional space, 
(generated by $x,y$ with defining relation $xy-yx=1$), and 
$g_i=\diag(\lambda_i,\lambda_i^{-1})$. 

Thus we need to calculate 
$H^*(B,Bg)$, $B=\Bbb A_1$, $g=\diag(\lambda,\lambda^{-1})$. 

\begin{proposition}
$H^*(B,Bg)$ is 1-dimensional, concentrated in degree 
$0$ if $\lambda=1$ and in degree $2$ otherwise. 
\end{proposition}

\begin{proof} 
If $B=\Bbb A_1$ then $B$ has the following Koszul resolution 
as a B-bimodule:
$$
B\otimes B\to B\otimes \Bbb C^2\otimes B\to B\otimes B\to B.
$$
Here the first map is given by the formula 
$$
b_1\otimes b_2\to b_1\otimes x\otimes yb_2-b_1\otimes y\otimes xb_2
$$
$$
-b_1y\otimes x\otimes b_2+b_1x\otimes y\otimes b_2,
$$ 
the second map is given by 
$$
b_1\otimes x\otimes b_2\to b_1x\otimes b_2-b_1\otimes xb_2,\ 
b_1\otimes y\otimes b_2\to b_1y\otimes b_2-b_1\otimes yb_2,
$$
and the third map is the multiplication. 

Thus the cohomology of $B$ with coefficients in $Bg$ can be computed by mapping 
this resolution into $Bg$ and taking the cohomology. 
This yields the following complex $C^\bullet$: 
\begin{equation}\label{comple}
0\to Bg\to Bg\oplus Bg\to Bg\to 0
\end{equation}
where the first nontrivial map is given by 
$bg\to [bg,y]\otimes x-[bg,x]\otimes y$, and the second nontrivial map is 
given by $bg\otimes x\to [x,bg]$, $bg\otimes y\to [y,bg]$. 

Consider first the case $g=1$. Equip the complex $C^\bullet$ with the 
Bernstein filtration ($\deg(x)=\deg(y)=1$), starting with $0,1,2$, for $C^0,C^1,C^2$, respectively
(this makes the differential perserve the filtration). 
Consider the associated graded complex $C_{gr}^\bullet$. In this complex, brackets are replaced with 
Poisson brackets, and thus it is easy to see that $C_{gr}^\bullet$ is the De Rham complex for 
the affine plane, so its cohomology is $\Bbb C$ in degree 0 and 0 in other
degrees. Therefore, the cohomology of $C^\bullet$ is the same. 

Now consider $g\ne 1$. In this case, 
declare that $C^0,C^1,C^2$ start in degrees 2,1,0 respectively
(which makes the differential preserve the filtration), and 
again consider the graded complex $C_{gr}^\bullet$. 

The graded Euler characteristic of this complex 
is \linebreak $(t^2-2t+1)(1-t)^{-2}=1$. 

The cohomology in the $C^0_{gr}$ term is 
the set of $b\in \Bbb C[x,y]$ such that $ab=ba^g$ for all $a$. This means 
that $H^0=0$.  

The cohomology of the $C^2_{gr}$ term is the quotient of $\Bbb C[x,y]$ by the ideal generated by 
$a-a^g$, $a\in \Bbb C[x,y]$. Thus the cohomology $H^2$ of the rightmost term is 1-dimensional, in degree 0. 
By the Euler characteristic argument, this implies that $H^1=0$.
The cohomology of the filtered complex $C^\bullet$ is therefore the same, and we are done. 
\end{proof} 

The proposition implies that in the n-dimensional case 
$H^*(B,Bg)$ is 1-dimensional, concentrated in 
degree ${\rm rank}(1-g)$. It is not hard to check that 
the group $G$ acts on the sum of these 1-dimensional spaces by 
simply premuting the basis vectors. Thus the theorem is proved. 
\end{proof}

\subsection{The universal deformation of $G\ltimes {\rm Weyl}(V)$}

\begin{theorem} The algebra $H_{1,c}[G,V]$, with formal $c$, is 
the universal deformation of $H_{1,0}[G,V]=G\ltimes {\rm Weyl}(V)$. 
That is, the map $f: \Bbb C[S]^G\to H^2(G\ltimes {\rm Weyl}(V))$
induced by this deformation is an isomorphism.  
\end{theorem}

\begin{proof}
The map $f$ is a map between spaces of the same dimension, so it suffices to show 
that $f$ is injective. For this purpose it suffices to show that 
that for any $\gamma\in \Bbb C[S]^G$, the algebra $H_{1,\hbar\gamma}$ 
over $\Bbb C[\hbar]/\hbar^2$ is a nontrivial deformation of 
$H_{1,0}$. For this, in turn, it suffices to show that 
$H_{1,\hbar\gamma}$ does not have a subalgebra isomorphic to 
${\rm Weyl}(V)[\hbar]/\hbar^2$. This can be checked directly, as explained in 
\cite{EG}, Section 2. 
\end{proof}

\subsection{Notes} The results of this lecture can be found 
in \cite{EG}. For their generalizations, see \cite{E}. 

\section{The center of the symplectic reflection algebra}

\subsection{The module $H_{t,c}\e$}

Let $H_{t,c}$ be a symplectic reflection algebra. 
Consider the bimodule $H_{t,c}\e$, which has a left action of $H_{t,c}$ and a right 
commuting action of $\e H_{t,c}\e$. It is obvious that 
${\rm End}_{H_{t,c}}H_{t,c}\e=\e H_{t,c}\e$ (with opposite
product). The following theorem shows that 
the bimodule $H_{t,c}\e$ has the double centralizer property. 

Note that we have a natural map $\xi_{t,c}: H_{t,c}\to {\rm End}_{\e H_{t,c}\e}H_{t,c}\e$. 

\begin{theorem}\label{doubcen}
$\xi_{t,c}$ is an isomorphism for any $t,c$. 
\end{theorem}

\begin{proof}
The complete proof is given \cite{EG}.
We will give the main ideas of the proof skipping straightforward technical details.
The first step is to show that the result is true in the graded case, $(t,c)=(0,0)$.
To do so, note the following general fact:

\begin{exercise}\label{freeac} If $X$ is an affine complex algebraic variety with algebra of functions $\O_X$
and $G$ a finite group acting freely on $X$ then the natural map $\xi_X: G\ltimes \O_X\to 
{\rm End}_{\O_X^G}\O_X$ is an isomorphism. 
\end{exercise} 

Therefore, the map $\xi_{0,0}: G\ltimes SV\to {\rm End}_{(SV)^G}(SV)$ 
is injective, and moreover becomes an isomorphism after localization 
to the field of quotients $\Bbb C(V)^G$. To show it's surjective, take $a\in {\rm End}_{(SV)^G}(SV)$. 
There exists $a'\in G\ltimes \Bbb C(V)$ which maps to $a$. Moreover, by Exercise \ref{freeac},
$a'$ can have poles only at fixed points of $G$ on $V$. 
But these fixed points form a subset of codimension $\ge 2$, so there can be no poles and we are done
in the case $(t,c)=(0,0)$. 

Now note that the algebra ${\rm End}_{\e H_{t,c}\e}H_{t,c}\e$
has an increasing integer filtration (bounded below) induced by the filtration on $H_{t,c}$. 
This is due to the fact that $H_{t,c}\e$ is a finitely generated $\e H_{t,c}\e$-module 
(since it is true in the associated graded situation, by Hilbert's 
theorem about invariants). Also, the natural map
${\rm gr}{\rm End}_{\e H_{t,c}\e}H_{t,c}\e\to 
{\rm End}_{{\rm gr}\e H_{t,c}\e}{\rm gr}H_{t,c}\e$
is clearly injective. Therefore, our result in the case $(t,c)=(0,0)$ implies that
this map is actually an isomorphism (as so is its composition with the associated graded of $\xi_{t,c}$). 
Identifying the two algebras by this isomorphism, 
we find that ${\rm gr}(\xi_{t,c})=\xi_{0,0}$. Since $\xi_{0,0}$ is an isomorphism, 
$\xi_{t,c}$ is an isomorphism for all $t,c$, as desired.
\footnote{Here we use the fact that the filtration is bounded from below. 
In the case of an unbounded filtration, it is possible for a map not to be an isomorphism 
if its associated graded is an isomorphism. An example of this is the operator 
of multiplication by $1+t^{-1}$ in the space of Laurent polynomials in $t$, filtered 
by degree.}   
\end{proof} 

\subsection{The center of $H_{0,c}$}

It turns out that if $t\ne 0$, then the center of $H_{t,c}$ is trivial
(it was proved by Brown and Gordon \cite{BGo}). However, if $t=0$, 
$H_{t,c}$ has a nontrivial center, which we will denote $Z_c$. 
Also, as before, we denote the spherical subalgebra 
$\e H_{0,c}\e$ by $B_{0,c}$. 

We have a natural homomorphism 
$\zeta_c: Z_c\to B_{0,c}$, defined by the formula $\zeta(z)=z\e$. 
We also have a natural injection $\tau_c: {\rm gr}(Z_c)\to Z_0=SV$.
These maps are clearly isomorphisms for $c=0$. 

\begin{proposition}
The maps $\zeta_c,\tau_c$ are isomorphisms for any $c$.  
\end{proposition}

\begin{proof} It is clear that the morphism $\tau_c$ is injective. 
Identifying ${\rm gr}B_{0,c}$ with $B_{0,0}=SV$, we get 
${\rm gr}\zeta_c=\tau_c$. This implies that ${\rm gr}(\zeta_c)$ is injective.

Now we show that $\zeta_c$ is an isomorphism. 
To do so, we will construct the inverse homomorphism $\zeta_c^{-1}$. 
Namely, take an element $b\in B_{0,c}$. Since the algebra $B_{0,c}$ is commutative, 
it defines an element in ${\rm End}_{B_{0,c}}(H_{0,c}\e)$. Thus using 
Theorem \ref{doubcen}, we can define $a:=\xi_{0,c}^{-1}(b)\in H_{0,c}$. 
Note that $a$ commutes with $H_{0,c}$ as an operator on $H_{0,c}\e$ since it is defined by
right multiplication by $b$. Thus $a$ is central. It is easy to see that $b=z\e$. 

Thus $\zeta_c$ is an filtration preserving isomorphism whose associated graded
is injective. This implies that ${\rm gr}\zeta_c=\tau_c$ is an isomorphism. 
The theorem is proved.  
\end{proof}

Thus the spectrum of $Z_c$ is equal to $M_c={\rm Spec}B_{0,c}$. 

\subsection{Finite dimensional representations of $H_{0,c}$.} 

Let $\chi\in M_c$ be a central character: $\chi: Z_c\to \Bbb C$. 
Let $(\chi)$ be the ideal in $H_{0,c}$ generated by the kernel of $\chi$. 

\begin{proposition}\label{gen} If $\chi$ is generic then 
$H_{0,c}/(\chi)$ is the matrix algebra of size $|G|$. In particular, 
$H_{0,c}$ has a unique irreducible representation $V_\chi$ 
with central character $\chi$. This representation 
is isomorphic to $\Bbb C G$ as a $G$-module.
\end{proposition}

\begin{proof} It is shown by a standard argument (which we will skip)
that it is sufficient to check the statement 
in the associated graded case $c=0$. In this case, for generic $\chi$  
$G\ltimes SV/(\chi)=G\ltimes {\rm Fun}(\O_\chi)$, where 
$\O_\chi$ is the (free) orbit of $G$ consisting of the points of $V$ that map to $\chi\in V/G$,
and ${\rm Fun}(\O_\chi)$ is the algebra of functions on $\O_\chi$.  
It is easy to see that this algebra is isomorphic to a matrix algebra, and 
has a unique irreducible representation, ${\rm Fun}(\O_\chi)$, 
which is a regular representation of $G$. 
\end{proof}

\begin{corollary}\label{est}
Any irreducible representation of $H_{0,c}$ has dimension $\le |G|$.  
\end{corollary}

\begin{proof} Since for generic $\chi$ the algebra $H_{0,c}/(\chi)$ is a matrix algebra, 
the algebra $H_{0,c}$ satisfies the standard polynomial identity
(the Amitsur-Levitzki identity) for matrices $N\times N$ ($N=|G|$):
$$
\sum_{\sigma\in S_{2N}}(-1)^\sigma X_{\sigma(1)}...X_{\sigma(2N)}=0.
$$
Next, note that since $H_{0,c}$ is a finitely generated $Z_c$-module
(by passing to the associated graded and using Hilbert's theorem), every irreducible
representation of $H_{0,c}$ is finite dimensional. 
If $H_{0,c}$ had an irreducible representation $E$ of dimension
$m>|G|$, then 
by the density theorem the matrix algebra $Mat_m$ would be a quotient of 
$H_{0,c}$. But the Amitsur-Levitski identity of degree $|G|$ is not satisfied 
for matrices of bigger size than $|G|$. Contradiction. 
Thus, $\dim E\le |G|$, as desired.      
\end{proof}

In general, for special central characters there are representations 
of $H_{0,c}$ of dimension $<|G|$. However, in some cases one can show 
that all irreducible representations have dimension exactly $|G|$.
For example, we have the following result. 

\begin{theorem}\label{ratcher}
Let $G=S_n$, $V=\h\oplus \h^*$, $\h=\Bbb C^{n}$ (the rational Cherednik algebra). 
Then for $c\ne 0$, every irreducible representation of $H_{0,c}$ 
has dimension $n!$ and is isomorphic to the regular representation of $S_n$. 
\end{theorem}

\begin{proof}
Let $E$ be an irreducible representation of $H_{0,c}$. 
Let us calculate the trace in $E$ of any permutation $\sigma\ne 1$. 
Let $j$ be an index such that $\sigma(j)=i\ne j$. 
Then $s_{ij}\sigma(j)=j$. Hence in $H_{0,c}$ we have 
$$
[y_j,x_is_{ij}\sigma]=[y_j,x_i]s_{ij}\sigma=cs_{ij}^2\sigma=c\sigma.
$$
Hence $\Tr_E(\sigma)=0$, and thus $E$ is a multiple of the regular representation of $S_n$. 
But by Theorem \ref{est}, $\dim E\le n!$, so we get that $E$ is the regular representation, as desired. 
\end{proof}

\subsection{Azumaya algebras}

Let $Z$ be a finitely generated commutative algebra over $\Bbb C$ (for simplicity without 
nilpotents), $M={\rm Spec Z}$ is the corresponding algebraic variety, and $A$ a finitely generated $Z$-algebra. 

\begin{definition}
$A$ is said to be an Azumaya algebra of degree $N$ if the completion $\hat A_\chi$ of $A$ at every maximal ideal 
$\chi$ in $Z$ is the matrix algebra of size $N$ over the completion $\hat Z_\chi$ of $Z$. 
\end{definition}

For example, if $E$ is an algebraic vector bundle on $M$ then ${\rm End}(E)$ is an Azumaya algebra. 
However, not all Azumaya algebras are of this form. 

\begin{exercise} Let $q$ be a root of unity of order $N$. 
Show that the algebra of functions on the quantum torus, 
generated by $X^{\pm 1},Y^{\pm 1}$ with defining relation 
$XY=qYX$ is an Azumaya algebra of degree $N$. 
Is this the endomorphism algebra of a vector bundle? 
\end{exercise}

It is clear that if $A$ is an Azumaya algebra then 
for every central character $\chi$ of $A$, 
$A/(\chi)$ is the algebra $Mat_N(\Bbb C)$ of complex $N$ by $N$ matrices, and 
every irreducible representation of $A$ has dimension $N$. 

The following important result is due to M. Artin.   

\begin{theorem} Let $A$ be a finitely generated (over $\Bbb C$) 
polynomial identity (PI) algebra of degree $N$ 
(i.e. all the polynomial relations of the matrix algebra of size $N$ are satisfied in $A$). 
Then $A$ is an Azumaya algebra if every irreducible representation of $A$ 
has dimension exactly $N$.  
\end{theorem}

Thus, by Theorem \ref{ratcher}, 
for $G=S_n$, $H_{0,c}$ for $c\ne 0$ is an Azumaya algebra of degree $n!$.
Indeed, this algebra is PI of degree $n!$ because the classical Dunkl representation 
embeds it into matrices of size $n!$ over $\Bbb C(x_1,...,x_n,p_1,...,p_n)^{S_n}$.  

Let us say that $\chi\in M$ is an Azumaya point 
if for some affine neighborhood $U$ of $\chi$ the localization of $A$ to $U$ 
is an Azumaya algebra. Obviously, the set $Az(M)$ of Azumaya points of $M$ is open. 

Now we come back to the study the space $M_c$ corresponding to a symplectic 
reflection algebra $H_{0,c}$. 

\begin{theorem} \label{smo}
The set $Az(M_c)$ coincides with the set of smooth points of $M_c$.
\end{theorem}

The proof of this theorem is given in the following two subsections. 

\begin{corollary} If $G=S_n$ and $V=\h\oplus \h^*$, $\h=\Bbb C^{n}$ (the rational Cherednik algebra case) then 
$M_c$ is a smooth algebraic variety for $c\ne 0$. 
\end{corollary}

\subsection{Cohen-Macaulay property and homological dimension}

To prove Theorem \ref{smo}, we will need some commutative algebra tools. 
Let $Z$ be a finitely generated commutative algebra over $\Bbb C$ without zero divisors. 
By Noether's normalization lemma, there exist elements $z_1,...,z_n\in Z$ which are 
algebraically independent, such that $Z$ is a finitely generated module over $\Bbb C[z_1,...,z_n]$.

\begin{definition} The algebra $Z$ is said to be Cohen-Macaulay if 
$Z$ is a locally free (=projective) module over $\Bbb C[z_1,...,z_n]$. 
\footnote{It was proved by Quillen that a locally free module over 
a polynomial algebra is free; this is a difficult theorem, which will not be needed here.}
\end{definition}

\begin{remark}
It was shown by Serre that if $Z$ is locally free over $\Bbb C[z_1,...,z_n]$ for {\bf some} choice
of $z_1,...,z_n$, then it happens for {\bf any} 
choice of them (such that $Z$ is finitely generated as a nodule over 
$\Bbb C[z_1,...,z_n]$).   
\end{remark}

\begin{remark} Another definition of the Cohen-Macaulay property is that 
the dualizing complex $\omega_Z^\bullet$ of $Z$ is concentrated in degree zero. 
We will not discuss this definition here. 
\end{remark}

It can be shown that the Cohen-Macaulay property is stable under 
localization. Therefore, it makes sense to make the following definition. 

\begin{definition} An algebraic variety $X$ is Cohen-Macaulay if 
the algebra of functions on every affine open set in $X$ is Cohen-Macaulay. 
\end{definition}

Let $Z$ be a finitely generated commutative algebra over $\Bbb C$ without zero divisors, and  
let $M$ be a finitely generated module over $Z$.

\begin{definition} $M$ is said to be Cohen-Macaulay
if for some algebraically independent $z_1,...,z_n\in Z$ such that $Z$ is finitely generated over 
$\Bbb C[z_1,...,z_n]$, $M$ is locally free over $\Bbb C[z_1,...,z_n]$.
\end{definition}

Again, if this happens for some $z_1,...,z_n$, then it happens for any of them. 
We also note that $M$ can be Cohen-Macaulay without $Z$ being Cohen-Macaulay, 
and that $Z$ is a Cohen-Macaulay algebra iff it is a Cohen-Macaulay module over itself. 

We will need the following standard properties of Cohen-Macaulay algebras and 
modules.

\begin{theorem}\label{cmac}
(i) Let $Z_1\subset Z_2$ be a finite extension of finitely generated commutative $\Bbb C$-algebras,
without zero divisors, and $M$ be a finitely generated module over $Z_2$. Then 
$M$ is Cohen-Macaulay over $Z_2$ iff it is Cohen-Macaulay over $Z_1$. 

(ii) Suppose that $Z$ is the algebra of functions on a smooth affine variety. 
Then a $Z$-module $M$ is Cohen-Macaulay if and only if it is projective.   
\end{theorem}

In particular, this shows that the algebra of functions 
on a smooth affine variety is Cohen-Macaulay. 
Algebras of functions on many singular varieties are also Cohen-Macaulay. 

\begin{exercise}
Show that the algebra of functions on the cone $xy=z^2$ is Cohen-Macaulay. 
\end{exercise}

Another tool we will need is {\bf homological dimension.} 
We will say that an algebra $A$ has homological dimension 
$\le d$ if any (left) $A$-module $M$ has a projective resolution of length 
$\le d$. The homological dimension of $A$ is the smallest integer having this property. 
If such an integer does not exist, $A$ is said to have infinite homological dimension. 

It is easy to show that the homological dimension of $A$ is $\le d$ if and only if 
for any $A$-modules $M,N$ one has ${\rm Ext}^i(M,N)=0$ for $i>d$.
Also, the homological dimension clearly does not decrease 
under taking associated graded of the algebra under 
a nonnegative filtration (this is clear from considering 
the spectral sequence attached to the filtration). 

It follows immediately from this definition that 
homological dimension is Morita invariant. 
Namely, recall that a Morita equivalence 
between algebras $A$ and $B$ is 
an equivalence of categories 
$A-{\rm mod}\to B-{\rm mod}$. Such an equivalence maps projective modules to projective ones, 
since projectivity is a categorical property ($P$ is projective iff
${\rm Hom}(P,?)$ is exact). This implies that if $A$ and $B$ are Morita equivalent then 
their homological dimensions are the same. 

Then we have the following important theorem. 

\begin{theorem}\label{fini} The homological dimension of a commutative 
finitely generated $\Bbb C$-algebra $Z$ is finite if and only if $Z$ is regular, i.e. is the algebra of functions 
on a smooth affine variety. 
\end{theorem}

\subsection{Proof of Theorem \ref{smo}}

First let us show that any smooth point $\chi$ of 
$M_c$ is an Azumaya point. Since $H_{0,c}={\rm End}_{B_{0,c}}H_{0,c}\e=
{\rm End}_{Z_c}(H_{0,c}\e)$, it is sufficient to 
show that the coherent sheaf on $M_c$ corresponding to the module 
$H_{0,c}\e$ is a vector bundle near $\chi$. 
By Theorem \ref{cmac} (ii), for this it suffices to show 
that $H_{0,c}\e$ is a Cohen-Macaulay $Z_c$-module.

To do so, first note that the statement is true for $c=0$. 
Indeed, in this case the claim is that $SV$ is a Cohen-Macaulay module over $(SV)^G$. 
but $SV$ is a polynomial algebra, which is Cohen-Macaulay, 
so the result follows from Theorem \ref{cmac}, (i). 

Now, we claim that if $Z,M$ are positively filtered
and ${\rm gr}(M)$ is a Cohen-Macaulay ${\rm gr}Z$-module then 
$M$ is a Cohen-Macaulay $Z$-module. Indeed, let $z_1,...,z_n$ be homogeneous algebraically independent 
elements of ${\rm gr}Z$ such that ${\rm gr}Z$ is a finite module over the subalgebra generated by them. 
Let $z_1',...,z_n'$ be their liftings to $Z$. 
Then $z_1',...,z_n'$ are algebraically independent, and 
the module $M$ over $\Bbb C[z_1',...,z_n']$ is finitely generated and (locally) free since so is 
the module ${\rm gr}M$ over $\Bbb C[z_1,...,z_n]$. 

Recall now that ${\rm gr}H_{0,c}\e=SV$, ${\rm gr}Z_c=(SV)^G$. Thus the $c=0$ case implies the general case, and we are done. 

Now let us show that any Azumaya point of $M_c$ is smooth. 
Let $U$ be an affine open set in $M_c$ consisting of Azumaya points. 
Then the localization $H_{0,c}(U):=H_{0,c}\otimes_{Z_c}\O_U$ is an Azumaya algebra. 
Moreover, for any $\chi\in U$, the unique irreducible representation of 
$H_{0,c}(U)$ with central character $\chi$ is the regular representation 
of $G$ (since this holds for generic $\chi$ by Proposition \ref{gen}). 
This means that $\e$ is a rank $1$ idempotent 
in $H_{0,c}(U)/(\chi)$ for all $\chi$. In particular, $H_{0,c}(U)\e$ is a vector bundle
on $U$. Thus the functor $F: \O_U-{\rm mod}\to H_{0,c}(U)-{\rm mod}$ given by the formula 
$F(Y)=H_{0,c}(U)\e\otimes_{\O_U}Y$ is an equivalence of categories
(the quasi-inverse functor is given by the formula $F^{-1}(N)=\e N$). 
Thus $H_{0,c}(U)$ is Morita equivalent to $\O_U$, and therefore 
their homological dimensions are the same. 

On the other hand, the homological dimension of $H_{0,c}$ is finite (in fact, it equals to $\dim V$). 
To show this, note that by the Hilbert syzygies theorem, the homological dimension of $SV$ 
is $\dim V$. Hence, so is the homological dimension of $G\ltimes SV$ (as 
${\rm Ext^*}_{G\ltimes SV}(M,N)={\rm Ext^*}_{SV}(M,N)^G$). Thus, since ${\rm gr}H_{0,c}=G\ltimes SV$, 
we get that $H_{0,c}$ has homological dimension $\le \dim V$. 
Hence, the homological dimension of $H_{0,c}(U)$ is also $\le \dim V$
(as the homological dimension clearly does not increase under the localization). 
But $H_{0,c}(U)$ is Morita equivalent to $\O_U$, so $\O_U$ has a finite homological dimension. 
By Theorem \ref{fini}, this implies that $U$ consists of smooth points. 

\begin{corollary} $Az(M_c)$ is also the set of points at which 
the Poisson structure of $M_c$ is symplectic.  
\end{corollary}

\begin{proof} The variety $M_c$ is symplectic outside of a subset
of codimension $2$, because so is $M_0$.
Thus the set $S$ of smooth points of $M_c$ where the top exterior power 
of the Poisson bivector vanishes is of codimension $\ge 2$. So By Hartogs' theorem $S$ 
is empty. Thus, every smooth point is symplectic, and the corollary follows from the theorem.  
\end{proof}

\subsection{The space $M_c$ for $G=S_n$.}

Let us consider the space $M_c$ for $G=S_n$, $V=\h\oplus \h^*$, where $\h=\Bbb C^n$. 
In this case we have only one parameter $c$ corresponding to the conjugacy class of reflections,
and there are only two essentially different cases: $c=0$ and $c\ne 0$. 

\begin{theorem} For $c\ne 0$ the space $M_c$ is isomorphic to the Calogero-Moser 
space ${\mathcal C}_n$ as a symplectic manifold. 
\end{theorem}

\begin{proof} To prove the theorem, we will first construct a map 
$f: M_c\to {\mathcal C}_n$, and then prove that $f$ is an isomorphism.  

Without loss of generality, we may assume that $c=1$. 
As we have shown before, the algebra $H_{0,c}$ is an Azumaya algebra. 
Therefore, $M_c$ can be regarded as the moduli space of irreducible representations of $H_{0,c}$. 

Let $E\in M_c$ be an irreducible representation of $H_{0,c}$. 
We have seen before that $E$ has dimension $n!$ and is isomorphic to the regular representation 
as a representation of $S_n$. Let $S_{n-1}\subset S_n$ be the subgroup which preserves the element $1$.
Then the space of invariants $E^{S_{n-1}}$ has dimension $n$. On this space we have operators 
$X,Y: E^{S_{n-1}}\to E^{S_{n-1}}$ obtained by restriction of the operators $x_1,y_1$ on $E$
to the subspace of invariants. We have 
$$
[X,Y]=T:=\sum_{i=2}^ns_{1i}.
$$
Let us now calculate the right hand side of this equation explicitly. 
Let $\p$ be the symmetrizer of $S_{n-1}$. Let us realize the regular
representation $E$ of $S_n$ as $\Bbb C[S_n]$ with action of $S_n$ by left multiplication. 
Then $v_1=\p,v_2=\p s_{12},...,v_n=\p s_{1n}$ is a basis of $E^{S_{n-1}}$. 
The element $T$ commutes with $\p$, so we have 
$$
Tv_i=\sum_{j\ne i}v_j
$$
This means that $T+1$ has rank $1$, and hence the pair $(X,Y)$ 
defines a point on the Calogero-Moser space 
$\mathcal{C}_n$. \footnote{Note that the pair $(X,Y)$ is well defined only up to conjugation, because 
the representation $E$ is well defined only up to an isomorphism.}

We now set $(X,Y)=f(E)$. It is clear that $f: M_c\to {\mathcal C}_n$ is a regular map. 
So it remains to show that $f$ is an isomorphism. 
This is equivalent to showing that the corresponding map 
of function algebras $f^*: \O({\mathcal C}_n)\to Z_c$ is an isomorphism. 

Let us calculate $f$ and $f^*$ more explicitly. 
To do so, consider the open set $U$ in $M_c$ consisting of representations 
in which $x_i-x_j$ act invertibly. These are exactly the representations that 
are obtained by restricting representations of 
$S_n\ltimes \Bbb C[x_1,...,x_n,p_1,...,p_n,\delta(x)^{-1}]$ using the classical Dunkl 
embedding. Thus representations $E\in U$ are of the form $E=E_{\lambda,\mu}$ 
($\lambda,\mu\in \Bbb C^n$, and $\lambda$ has distinct coordinates), where 
$E_{\lambda,\mu}$ is the space of complex valued functions 
on the orbit $O_{\lambda,\mu}\subset \Bbb C^{2n}$, with the following action of 
$H_{0,c}$: 
$$
(x_iF)(a,b)=a_iF(a,b), (y_iF)(a,b)=b_iF(a,b)+\sum_{j\ne i}\frac{(s_{ij}F)(a,b)}{a_i-a_j}.
$$
(the group $S_n$ acts by permutations). 

Now let us consider the space $E_{\lambda,\mu}^{S_{n-1}}$.
A basis of this space is formed by characteristic functions of $S_{n-1}$-orbits  
on $O_{\lambda,\mu}$. 
Using the above presentation, it is straightforward to calculate the matrices
of the operators $X$ and $Y$ in this basis:
$$
X={\rm diag}(\lambda_1,...,\lambda_n), 
$$
and 
$$
Y_{ij}=\mu_i
\text{ if }j=i,\ Y_{ij}=\frac{1}{\lambda_i-\lambda_j}\text{ if }j\ne i. 
$$
This shows that $f$ induces an isomorphism $f|_U: U\to U_n$, where 
$U_n$ is the subset of ${\mathcal C}_n$ consisting of pairs $(X,Y)$ for which $X$ has distinct eigenvalues. 

From this presentation, it is straigtforward that 
$f^*(\Tr(X^p))=x_1^p+...+x_n^p$ for every positive integer $p$. 
Also, $f$ commutes with the natural $SL_2(\Bbb C)$ action
on $M_c$ and ${\mathcal C}_n$ (by $(X,Y)\to (aX+bY,cX+dY)$), so we also get 
$f^*(\Tr(Y^p))=y_1^p+...+y_n^p$, and 
$$
f^*(Tr(X^pY))=\frac{1}{p+1}\sum_{m=0}^p\sum_i x_i^my_ix_i^{p-m}.
$$
Now, using the necklace bracket formula on ${\mathcal C}_n$ and 
the commutation relations of $H_{0,c}$, 
we find, by a direct computation, that 
$f^*$ preserves Poisson bracket on the elements 
$Tr(X^p)$, $Tr(X^qY)$. But these elements are a local coordinate system near a generic point, 
so it follows that $f$ is a Poisson map. Since the algebra $Z_c$ is Poisson generated by 
$\sum x_i^p$ and $\sum y_i^p$ for all $p$, we get that $f^*$ is a surjective map. 

Also, $f^*$ is injective. Indeed, by Wilson's theorem the Calogero-Moser space is connected, and hence 
the algebra ${\mathcal O}({\mathcal C}_n)$ has no zero divisors, while 
${\mathcal C}_n$ has the same dimension as $M_c$. This proves that $f^*$ is an isomorphism, 
so $f$ is an isomorphism. 
\end{proof}

\subsection{Generalizations}

Now let $L$ be a 2-dimensional symplectic vector space, 
$\Gamma\subset Sp(L)$ is a finite subgroup, and $G=S_n\ltimes \Gamma^n$. 
Let $V=L\oplus...\oplus L$ ($n$ times). In this situation 
$c$ consists of a number $k$ corresponding to the conjugacy class  
of $s_{ij}$ (for $n>1$) and numbers $c_s$ corresponding to conjugacy classes 
of nontrivial elements in $\Gamma$. In this case, the above results may be generalized. 
More specifically, we have the following result (which we will state without proof). 

\begin{theorem}
(i) $M_c$ is smooth for generic $c$. 
(It has a matrix realization in terms of McKay's correspondence). 

(ii) $M_c$ is a deformation of the Hilbert scheme of 
$n$-tuples of points on the resolution of the Kleinian singularity $\Bbb C^2/\Gamma$. 
In particular, if $n=1$ then $M_c$ is the versal deformation of the 
Kleinian singularity. 
\end{theorem}

Part (i) of this theorem is proved in \cite{EG}.
The proof of part (ii) of this theorem for $\Gamma=1$ can be
found in Nakajima's book \cite{Na}. The general case is similar.

We also have the following result about the cohomology of $M_c$.
Introduce a filtration on $\Bbb C[G]$ by setting the degree 
of $g$ to be the rank of $g-1$ on $V$. 

\begin{theorem} The cohomology ring of $M_c$ 
is $H^*(M_c,\Bbb C)={\rm gr}({\rm Center}(\Bbb C[G]))$. 
\end{theorem}

The proof is based on the following argument. 
We know that the algebra $B_{t,c}$ is a quantization of $M_c$. 
Therefore, by Kontsevich's ``compatibility with cup products'' 
theorem, the cohomology algebra of $M_c$ is the cohomology of $B_{t,c}$ (for generic $t$). 
But $B_{t,c}$ is Morita equivalent to $H_{t,c}$, so 
this cohomology is the same as the Hochschild cohomology of $H_{t,c}$. 
However, the latter can be computed by using that $H_{t,c}$ is given by generators and relations
(by producing explicit representatives of cohomology classes and
computing their product). 

\subsection{Notes} The results of this lecture are contained in
\cite{EG}. For the commutative algebraic tools (Cohen-Macaulayness,
homological dimension) we refer the reader to the textbook
\cite{Ei}. The Amitsur-Levitzki identity and 
Artin's theorem on Azumaya algebras can be found 
in \cite{Ja}. The center of the double affine Hecke algebra of
type A, which is a deformation of the rational Cherednik algebra,
was computed by Oblomkov \cite{Ob}. 

\section{Representation theory of rational Cherednik algebras}

\subsection{Rational Cherednik algebras for any finite group of
linear transformations}

Above we defined rational Cherednik algebras for reflection
groups, as a tool for understanding Olshanetsky-Perelomov
operators. Actually, it turns out that such algebras can be
defined for any finite group of linear transformations. 

Namely, let $\h$ be an $\ell$-dimensional complex vector space, 
and $W$ be a finite group of linear transformations
of $\h$. An element $s$ of $W$ is said to be 
a {\it complex reflection} if all eigenvalues of $s$ in $\h$ are
equal to $1$ except one eigenvalue $\lambda\ne 1$ (which will of
course be a root of unity). For each $s\in S$, let $\alpha_s$
be an eigenvector of $s$ in $\h^*$ with nontrivial eigenvalue,
and $\alpha_s^\vee$ the eigenvector of $s$ in $\h$ such that 
$(\alpha_s,\alpha_s^\vee)=2$. Let $S$ be the set of conjugacy
classes of complex reflections in $W$, and $c: S\to \Bbb C$ be a
function invariant under conjugation. 

\begin{definition} The rational Cherednik algebra 
$H_{1,c}=H_{1,c}(\h,W)$ attached to $\h,W,c$ is 
the quotient of the algebra \linebreak
$\Bbb CW\ltimes \bold T(\h\oplus \h^*)$ 
(where $\bold T$ denotes the tensor algebra) 
by the ideal generated by the relations
$$
[x,x']=0,\ [y,y']=0,\ [y,x]=(y,x)-\sum_{s\in S}
c_s(y,\alpha_s)(x,\alpha_s^\vee)s.
$$
\end{definition}

Obviously, this is a generalization of the definition 
of the rational Cherednik algebra for a reflection group, 
and a special case of symplectic reflection algebras. 
In particular, the PBW theorem holds for $H_{1,c}$, and thus we
have a tensor product decomposition $H_{1,c}=\Bbb C[\h]\otimes 
\Bbb CW\otimes \Bbb C[\h^*]$ (as a vector space). 

\subsection{Verma and irreducible lowest weight modules over $H_{1,c}$}

We see that the structure of $H_{1,c}$ is similar to the
structure of the universal enveloping algebra of a simple Lie
algebra: $U(\g)=U(\n_-)\otimes U(\h)\otimes U(\n_+)$. 
Namely, the subalgebra $\Bbb CW$ plays the role of the Cartan
subalgebra, and the subalgebras $\Bbb C[\h^*]$ and $\Bbb C[\h]$ 
play the role of the positive and negative nilpotent
subalgebras. This similarity allows one to define 
and study the category $\O$ analogous to the
Bernstein-Gelfand-Gelfand category $\O$ for simple Lie algebras,
in which lowest weights will be representations of
$W$. 

Let us first define the simplest objects of this category -- the
standard, or Verma modules. 

For any irreducible representation $\tau$ of $W$,
let $M_c(\tau)$ be the standard representation
of $H_{1,c}$ with lowest weight
$\tau$, i.e. \linebreak $M_c(\tau)=H_{1,c}\otimes_{\Bbb C[W]\ltimes \Bbb
C[\h^*]}\tau$,
where $\tau$ is the representation of
$\Bbb CW\ltimes \Bbb C[\h^*]$, in which $y\in \h$
act by $0$. Thus, as a vector space,
$M_c(\tau)$ is naturally identified with $\Bbb C[\h]\otimes
\tau$.

An important special case of $M_c(\tau)$ is
$M_c=M_c(\Bbb C)$, the polynomial representation,
corresponding to the case when $\tau=\Bbb C$ is trivial.
The polynomial representation can thus be
naturally identified with $\Bbb C[\h]$.
Elements of $W$ and $\h^*$ act in this space in the obvious way,
while elements $y\in \h$ act by Dunkl operators
$$
\partial_y-\sum_{s\in
S}\frac{2c_s}{1-\lambda_s}\frac{(\alpha_s,y)}{\alpha_s}
(1-s),
$$
where $\lambda_s$ is the nontrivial eigenvalue of $s$ in 
$\h^*$. 

Introduce 
an important element $\bold h\in H_{1,c}$:
$$
\bold h=\sum_i x_iy_i+\frac{\ell}{2}-\sum_{s\in
S}\frac{2c_s}{1-\lambda_s}s,
$$
where $y_i$ is a basis of $\h$ and $x_i$ the dual basis of
$\h^*$.

\begin{exercise}
Show that if $W$ preserves an inner product in $\h$
then the element $\bold h$ can be included in an $sl_2$ triple
$\bold h,\bold E=\frac{1}{2}\sum x_i^2,\bold F=-\frac{1}{2}\sum
y_i^2$, where $x_i$, $y_i$ are orthonormal bases of $\h^*$ and
$\h$, respectively.
\end{exercise}

\begin{lemma}\label{elemh}
The element $\bold h$ is $W$-invariant and satisfies 
the equations $[\bold h,x]=x$ and
$[\bold h,y]=-y$.
\end{lemma}

\begin{exercise} Prove Lemma \ref{elemh}.
\end{exercise}

It is easy to show that the lowest eigenvalue of
$\bold h$ in $M_c(\tau)$ is $h(\tau):=
\ell/2-\sum_s \frac{2c_s}
{1-\lambda_s}s|_\tau$. 

\begin{corollary} \label{irrquo}
The sum $J_c(\tau)$ of all proper submodules of $M_c(\tau)$ 
is a proper submodule of $M_c(\tau)$, and the quotient module 
$L_c(\tau):=M_c(\tau)/J_c(\tau)$ is irreducible. 
\end{corollary}

\begin{proof} Every eigenvalue $\mu$ of $\bold h$ in a proper submodule
of $M_c(\tau)$, and hence in $J_c(\tau)$, satisfies the
inequality $\mu-h(\tau)>0$. Thus $J_c(\tau)\ne M_c(\tau)$. 
\end{proof} 

The module $L_c(\tau)$ can be characterized in terms of the
contragredient standard modules.
Namely, let $\hat M_c(\tau)=\tau^*\otimes_{\Bbb C[W]\ltimes \Bbb
C[\h]}H_{1,c}(W)$
be a right $H_{1,c}(W)$-module, and $M_c(\tau)^\vee=\hat M_c(\tau)^*$
its restricted dual, which
may be called the contragredient standard module.
Clearly, there is a natural morphism
$\phi: M_c(\tau)\to M_c(\tau)^\vee$.

\begin{lemma}\label{ima}
The module $L_c(\tau)$ is the image of $\phi$.
\end{lemma}

\begin{exercise}
Prove Lemma \ref{ima}.
\end{exercise}

Note that the map $\phi$ can be viewed as a
bilinear form $B: \hat M_c(\tau)\otimes M_c(\tau)\to \Bbb C$.
This form is analogous to the Shapovalov form in Lie theory.

If $W$ preserves an inner product on $\h$,
we can define an antiinvolution of $H_c(W)$
by $x_i\to y_i$, $y_i\to x_i$, $g\to g^{-1}$ for
$g\in W$ (where $x_i,y_i$ are orthonormal bases of $\h^*$ and
$\h$ dual to each other). Under this antiinvolution, the right
module $\hat M_c(\tau)$ turns into the left module
$M_c(\tau^*)$, so the form $B$ is a (possibly degenerate) pairing
$M_c(\tau^*)\otimes M_c(\tau)\to \Bbb C$.
Moreover, it is clear that if $Y,Y'$ are any quotients
of $M_c(\tau),M_c(\tau^*)$ respectively, then $B$ descends to
a pairing $Y'\otimes Y\to \Bbb C$
(nondegenerate iff $Y,Y'$ are irreducible). This pairing
satisfies the contravariance equations $B(a,x_ib)=B(y_ia,b)$,
$B(a,y_ib)=B(x_ia,b)$, and $B(ga,gb)=B(a,b)$ for $g\in W$.

\subsection{Category $\O$.}

\begin{definition} 
The category $\mathcal O=\mathcal O_c$ of $H_{1,c}$-modules is the category
of $H_{1,c}$-modules $V$, such that $V$ is the direct sum of
finite dimensional generalized eigenspaces of $\bold h$, and the
real part of the spectrum of $\bold h$ is bounded below. 
\end{definition}

Obviously, the standard representations $M_c(\tau)$ and
their irreducible quotients $L_c(\tau)$ belong to ${\mathcal O}$.

\begin{exercise} Show that $\O$ 
is an abelian subcategory of the category of
all $H_{1,c}$-modules, which is closed under extensions, 
and its set of isomorphism classes of simple objects
is $\lbrace{ L_c(\tau)\rbrace}$.  
Moreover, show that every object of $M_c(\tau)$ has finite length. 
\end{exercise}

\begin{definition}
The character of a module $V\in {\mathcal O}$ is
$\chi_V(g,t)=\Tr|_V(g t^\bold h)$, $g\in W$ (this is a series in
$t$).
\end{definition}

For example, the character of $M_c(\tau)$ is
$$
\chi_{M_c(\tau)}(g,t)=\frac{\chi_\tau(g)t^{h(\tau)}}{\det|_{\h^*}(1-gt)}.
$$

On the other hand, determining the characters 
of $L_c(\tau)$ is in general a hard problem. 
An equivalent problem is determining the multiplicity 
of $L_c(\sigma)$ in the Jordan-H\"older series of $M_c(\tau)$. 

\begin{exercise} 1. Show that for generic $c$ (outside of countably
many hyperplanes), $M_c(\tau)=L_c(\tau)$, so $\O$ is a semisimple
category. 

2. Determine the characters of $L_c(\tau)$ for all $c$ in
the case $\ell=1$.
\end{exercise}

\subsection{The Frobenius property}

Let $A$ be a $\Bbb Z_+$-graded commutative algebra. 
The algebra $A$ is called Frobenius if the top degree $A[d]$ of
$A$ is 1-dimensional, and the multiplication map 
$A[m]\times A[d-m]\to A[d]$ is a nondegenerate pairing for any
$0\le m\le d$. In particular, the Hilbert polynomial 
of a Frobenius algebra $A$ is palindromic.  

Now, let us go back to considering modules 
over the rational Cherednik algebra $H_{1,c}$. 
Any submodule $J$ of the polynomial 
representation $M_c(\Bbb C)=M_c=\Bbb C[\h]$ is an ideal in $\Bbb C[\h]$, 
so the quotient $A=M_c/J$ is a $\Bbb Z_+$-graded commutative algebra. 

Now suppose that $W$ preserves an inner product in $\h$. 

\begin{theorem} \label{Gorirr} If 
$A=M_c/J$ is finite dimensional, then $A=L_c=L_c(\Bbb C)$
(i.e. $A$ is irreducible) if and only if it is a Frobenius
algebra. 
\end{theorem}

\begin{proof}
(i) Suppose $A$ is an irreducible $H_{1,c}$-module,
i.e. $A=L_c(\Bbb C)$.  
Recall that $H_{1,c}$ has an ${\frak sl}_2$ subalgebra 
generated by $\bold h$, $\bold E=\sum x_i^2$, $\bold F=\sum
y_i^2$. Thus the ${\frak {sl}}_2$-representation theory shows that
 the top degree homogeneous component of $A$ is 1-dimensional. 
Let $\phi\in A^*$ denote a nonzero linear functional on $A$ which
factors through the projection to the top component; it is unique
up to scaling. 

Further, since $A$ is finite dimensional, the action of ${\frak
sl}_2$ in $A$ integrates to an action of the group 
$SL_2(\Bbb C)$. Also, as we have explained, this module has a canonical 
contravariant nondegenerate bilinear form $B(,)$. 
 Let $F$ be the `Fourier transform'
endomorphism of $A$ corresponding to the
action of the matrix $\left(\begin{matrix} \,0 & 1\\ -1 & 0
\end{matrix}\right)\in SL_2(\Bbb C)$.
Since the Fourier transform interchanges 
$\h$ and $\h^*$, the bilinear form: $E: (v_1,v_2)\longmapsto B(v_1\,,\,Fv_2)$
on  $A$ is compatible with the algebra
structure, i.e. for any $x\in \h^*$, we have
$E(xv_1,v_2)=E(v_1,xv_2)$. 
This means that for any polynomials 
$p,q\in \Bbb C[\h]$, one has $E(p(x),q(x))=\phi(p(x)q(x))$
(for a suitable normalization of $\phi$). 
Hence, since the form $E$ is nondegenerate,
$A$ is a Frobenius algebra. 

(ii) Suppose $A$ is Frobenius. Then the highest degree component of $A$ is 
1-dimensional, and the pairing $E: A\otimes A\to \Bbb C$ given by 
$E(a,b):= \phi(ab)$ (where, as before, $\phi$ stands for the 
highest degree coefficient) is nondegenerate. 
This pairing obviously satisfies the condition 
$E(xa,b)=E(a,xb), x\in \h^*$. Set  
$\tilde B(a,b):=E(a,Fb)$. Then $B$ satisfies the 
equations $\tilde B(a,x_ib)=\tilde B(y_ia,b)$.
So for any $f_1,f_2\in \Bbb C[\h]$, 
one has $\tilde B(p(x)v,q(x)v)=
\tilde B(q(y)p(x)v,v)$, where $v=1$ is the lowest weight vector
of $A$. This shows that $\tilde B$ coincides with the Shapovalov form $B$ 
on $A$, up to scaling. Thus $A$ is an irreducible representation of $H_{1,c}$. 
\end{proof} 

{\bf Remark.} It easy to see by considering the rank 1 case that 
for complex reflection groups Theorem \ref{Gorirr} is, in
general, false.  

Let us now consider the Frobenius property 
of $L_c$ in the case of a general group $W\subset GL(\h)$. 

\begin{theorem} \label{fingor}
Let $U\subset M_c=\Bbb C[\h]$ 
be a $W$-subrepresentation of dimension $\ell=\dim(\h)$ 
sitting in degree $r$, consisting of singular vectors 
(i.e. those killed by $y\in \h$). 
Let $J$ be the ideal (=submodule) generated by $U$. Assume that the
quotient representation $A=M_c/J$ is finite dimensional. 
Then 

(i) The algebra $A$ is Frobenius. 

(ii) The representation $A$ admits a BGG type resolution 
$$
A\leftarrow M_c(\Bbb C)\leftarrow M_c(U)\leftarrow M_c(\wedge^2U)
\leftarrow...\leftarrow M_c(\wedge^\ell U)\leftarrow 0. 
$$

(iii) The character of $A$ is given by the formula 
$$
\chi_A(g,t)=t^{\frac{\ell}{2}-\sum_s \frac{2c_s}{1-\lambda_s}}
\frac{\det|_U(1-gt^r)}{\det_{\h^*}(1-gt)}.
$$
In particular, the dimension of $A$ is $r^\ell$. 

(iv) If $W$ preserves an inner product on $\h$, then $A$ is irreducible. 
\end{theorem}

\begin{proof}
(i) Since ${\rm Spec}(A)$ is a complete intersection,
$A$ is Frobenius (=Gorenstein of dimension 0) (see \cite{Ei}, p.541).

(ii) Consider the subring $SU$ in $\Bbb C[\h]$.  
Then $\Bbb C[\h]$ is a finitely generated $SU$-module. 
A standard theorem of Serre which we have discussed above
(Section 10.5) says that 
if $B=\Bbb C[t_1,...,t_n]$, $f_1,...,f_n\in B$ are homogeneous, 
$A=\Bbb C[f_1,...,f_n]\subset B$, and 
$B$ is a finitely generated module over $A$, then
$B$ is a free module over $A$. Applying this in our situation,  
we see that $\Bbb C[\h]$ is a free $SU$-module. 
It is easy to see by computing the 
Hilbert series that the rank of this free module 
is $r^\ell$. 

Consider the Koszul complex attached to the module $\Bbb C[\h]$
over $SU$. Since this module is free, 
the Koszul complex is exact (i.e. it is a resolution of the zero-fiber).  
At the level of $\Bbb C[\h]$ modules, this resolution 
looks exactly as we want in (ii). So we need to show 
that the maps of the resolution are in fact morphisms 
of $H_{1,c}$-modules and not only $\Bbb C[\h]$-modules.
This is easily established by induction (going from left to right),
cf. proof of Proposition 2.2 in \cite{BEG}. 

(iii) Follows from (ii) by the Euler-Poincare principle. 

(iv) Follows from Theorem \ref{Gorirr}.
\end{proof}
 
\subsection{Representations of the rational Cherednik algebra of type $A$}

\subsubsection{The results} 

Let $W=S_n$, and $\h$ be its reflection representation. 
In this case the function $c$ reduces to one number $k$. 
We will denote the rational Cherednik algebra $H_{1,k}(S_n)$ by $H_k(n)$. 
The polynomial representation $M_k$ of this algebra is the space of 
$\Bbb C[x_1,...,x_n]^T$ of polynomials of $x_1,...,x_n$, which are 
invariant under simultaneous translation $x_i\mapsto x_i+a$. 
In other words, it is the space of regular functions
on $\h=\Bbb C^n/\Delta$, where $\Delta$ is the diagonal.

\begin{proposition} \label{Asing} (C. Dunkl)  
Let $r$ be a positive integer not divisible by $n$, and $k=r/n$. 
Then $M_k$ contains a copy of the reflection representation 
$\frak h$ of $S_n$, which consists of singular vectors
(i.e. those killed by $y\in \h$). 
This copy sits in degree $r$ and is spanned 
by the functions 
$$
f_i(x_1,...,x_n)={\rm Res}_\infty [(z-x_1)...(z-x_n)]^{\frac{r}{n}}
\frac{dz}{z-x_i}.
$$
(the symbol ${\rm Res}_\infty$ denotes the residue 
at infinity). 
\end{proposition} 

{\bf Remark.} The space spanned by $f_i$ is $n-1$-dimensional, since 
$\sum_i f_i=0$ (this sum is the residue of an exact differential). 

\begin{proof}
This proposition can be proved by a straightforward computation. 
The functions $f_i$ are a special case of 
Jack polynomials.
\end{proof}
       
\begin{exercise} Do this computation.
\end{exercise}
       
\bigskip
	
Let $I_k$ be the submodule of $M_k$ 
generated by $f_i$. 
Consider the $H_k(n)$-module $V_k=M_k/I_k$, and 
regard it as a $\Bbb C[\h]$-module.

Our result is 

\begin{theorem}\label{main} Let $d=(r,n)$ denote the greatest common 
divisor of $r$ and $n$. 
Then the (set-theoretical) 
support of $V_k$ is the union of $S_n$-translates 
of the subspaces of $\Bbb C^n/\Delta$,
defined by the equations 
$$
\begin{array}{c}
x_1=x_2 = \dots =x_{\frac{n}{d}};\\
x_{\frac{n}{d}+1}= \dots =x_{2\frac{n}{d}};\\
 \dots  \\
x_{(d-1)\frac{n}{d}+1}= \dots = x_n.
\end{array}
$$
In particular, the Gelfand-Kirillov dimension of $V_k$ 
is $d-1$. 
\end{theorem}

The theorem allows us to give a simple proof 
of the following result of Berest, Etingof, Ginzburg (without the use of 
the KZ functor and Hecke algebras). 

\begin{corollary} If $d=1$ 
then the module $V_k:=M_k/I_k$ is finite dimensional, irreducible, 
admits a BGG type resolution, and its character is 
$$
\chi_{V_k}(g,t)=t^{(1-r)(n-1)/2}\frac{\det|_\h(1-gt^r)}{\det|_\h(1-gt)}.
$$
\end{corollary}

\begin{proof} 
For $d=1$ 
Theorem \ref{main} says that the support of 
$M_k/I_k$ is $\{0\}$. This implies that 
$M_k/I_k$ is finite dimensional. The rest follows from  
Theorem \ref{fingor}.
\end{proof} 

\subsubsection{Proof of Theorem \ref{main}}
The support of $V_k$ is the zero-set of $I_k$, i.e. 
the common zero set of $f_i$. 
Fix $x_1,...,x_n\in \Bbb C$. Then $f_i(x_1,...,x_n)=0$ for all
$i$ iff 
$\displaystyle \sum_{i=1}^n \lambda_i f_i=0$ 
for all $\lambda_i$, i.e.
$$
 \operatorname{Res}_{\infty}\left( 
\prod_{j=1}^n(z-x_j)^{\frac{r}{n}}\sum_{i=1}^n\frac{\lambda_i}{z-x_i} 
\right)dz=0.  
$$

Assume that $x_1, \dots x_n$ take distinct values $y_1, \dots, y_p$ with 
positive multiplicities $m_1, \dots, m_p$.
The previous equation implies that the point 
$(x_1,...,x_n)$ is in the zero set iff
$$
{\rm Res}_\infty 
\prod_{j=1}^p (z-y_j)^{m_j\frac{r}{n}-1} \left(\sum_{i=1}^p
\nu_i(z-y_1)\dots \widehat{(z-y_i)} \dots (z-y_p)\right)dz=0 \quad
\forall \nu_i.
$$
Since $\nu_i$ are arbitrary, this is equivalent to the 
condition
$$
{\rm Res}_\infty 
\prod_{j=1}^p (z-y_j)^{m_j\frac{r}{n}-1}z^i dz=0, \quad i=0, \dots, 
p-1. 
$$

We will now need the following lemma. 

\begin{lemma}\label{mainlemma}
Let $\displaystyle a(z)=\prod_{j=1}^p(z-y_j)^{\mu_j}$, where 
$\mu_j\in \Bbb C$, $\sum_j \mu_j 
\in \mathbb Z$ 
and $\sum_j \mu_j > -p$. Suppose 
$${\rm Res}_\infty a(z)z^idz=0, \quad i=0,1,\dots, p-2.$$
Then $a(z)$ is a polynomial.
\end{lemma}

\medskip

\begin{proof}
Let $g$ be a polynomial.
Then 
$$0={\rm Res}_\infty d(g(z)\cdot a(z))=
{\rm Res}_\infty(g^{\prime}(z)a(z)+a^{\prime}(z)g(z))dz$$

and hence

$${\rm Res}_\infty \left(g^{\prime}(z)+ \sum_i 
\frac{\mu_j}{z-y_j}g(z)\right)a(z)dz=0.$$

\smallskip

Let $\displaystyle g(z)=z^l \prod_j(z-y_j)$. Then $ 
\displaystyle g^{\prime}(z)+ \sum 
\frac{\mu_j}{z-y_j}g(z)$ is a polynomial of degree $l+p-1$ with highest 
coefficient $l+p+\sum \mu_j \ne 0$ (as $\sum \mu_j>-p$). 
This means that for every $l \ge 0$,
$\displaystyle {\rm Res}_\infty z^{l+p-1}a(z)dz$ is a linear combination of 
residues of $z^qa(z)dz$ 
with $q<l+p-1$. By the assumption of the lemma, 
this implies by induction in $l$ 
that all such residues are $0$ and hence $a$ is a polynomial.
\end{proof}

In our case $\sum (m_j \frac{r}{n}-1)=r-p$ (since $\sum m_j=n$) and 
the conditions of the lemma are satisfied. Hence 
$(x_1,...,x_n)$ is in the zero set of $I_k$ iff
$\displaystyle\prod_{j=1}^p(z-y_j)^{m_j\frac{r}{n}-1}$ is a polynomial.
This is equivalent to saying that all  $m_j$ are divisible by 
$\frac{n}{d}$. 

We have proved that $(x_1, \dots x_n)$ is in the zero set of $I_k$ iff
$(z-x_1) \dots (z-x_n)$ is the $n/d$-th power of a polynomial of degree $d$.
This implies the theorem. 

\subsection{Notes} In this lecture we have followed the paper 
\cite{CE}. The classification of finite dimensional 
representations of the rational Cherednik algebra of type A 
was obtained in \cite{BEG} (the corresponding result 
for double affine Hecke algebra is due to Cherednik \cite{Ch2}).
The general theory of category ${\mathcal O}$ for rational Cherednik algebras 
is worked out in \cite{GGOR}.

\end{document}